\documentclass[EJP]{ejpecp} 

\usepackage{enumerate}  
\usepackage{mathrsfs}

\newcounter{epsilon}
  \newcommand{\nep}[1]{\refstepcounter{epsilon}\label{#1}}
  \newcommand{\uep}[1]{\varepsilon_{\textnormal{\tiny \ref{#1}}}}
\newcounter{constant}
  \newcommand{\nc}[1]{\refstepcounter{constant}\label{#1}}
  \newcommand{\uc}[1]{c_{\textnormal{\tiny \ref{#1}}}}
\newcounter{Konstant}
  \newcommand{\nC}[1]{\refstepcounter{Konstant}\label{#1}}
  \newcommand{\uC}[1]{C_{\textnormal{\tiny \ref{#1}}}}

\SHORTTITLE{FPPHE is not monotone}

\TITLE{First passage percolation in hostile environment is not monotone
    \thanks{E.C.\ was supported by the project ``Programma per Giovani Ricercatori Rita Levi Montalcini'' awarded by the Italian Ministry of Education. E.C.\ also acknowledges partial support by ``INdAM--GNAMPA Project 2019'' and ``INdAM--GNAMPA Project 2020''. A.S.\ acknowledges support from EPSRC Early Career Fellowship EP/N004566/1.}
      } 

\AUTHORS{%
Elisabetta Candellero\footnote{elisabetta.candellero@uniroma3.it;  Universit\`a Roma Tre, Dipartimento di Matematica e Fisica, Largo S.\ Murialdo 1, 00146, Rome, Italy.} \and Alexandre Stauffer\footnote{a.stauffer@bath.ac.uk;  University of Bath, Dept of Mathematical Sciences, BA2 7AY Bath, UK \ \& \ Universit\`a Roma Tre, Dipartimento di Matematica e Fisica, Largo S.\ Murialdo 1, 00146, Rome, Italy.}  
}


\AMSSUBJ{60K35; 
82C43
} 

\SUBMITTED{August 19, 2022} 
\ACCEPTED{May 14, 2024} 

\VOLUME{29}
\YEAR{2024}
\PAPERNUM{85}
\DOI{10.1214/24-EJP1145}

\ABSTRACT{We study a natural growth process with competition, modeled by two first passage percolation processes, $\Fo$ and $\Fl$, spreading on a graph. 
$\Fo$ starts at the origin and spreads at rate $1$, 
whereas $\Fl$ starts from a random set of \emph{inactive seeds} distributed as Bernoulli percolation of parameter $\mu\in (0,1)$.
A seed of $\Fl$ gets activated when one of the two processes attempts to occupy its location, and from this moment onwards spreads at some fixed rate $\lambda>0$.
In previous works~\cite{Stauffer-Sidoravicius-MDLA, Candellero-Stauffer-FppCompetition, Finn-Stauffer-FPPHE} it has been shown that when 
both $\mu$ or $\lambda$ are small enough, then $\Fo$ \emph{survives} (i.e., it occupies an infinite set of vertices) with positive probability.
It might seem intuitive that decreasing $\mu$ or $\lambda$ is beneficial to $\Fo$. However, we prove that, in general, 
this is indeed false by constructing a graph for which the probability that $\Fo$ survives is not a monotone function of $\mu$ or $\lambda$, implying the existence of multiple phase transitions.
This behavior contrasts with other natural growth processes such as the $2$-type Richardson model.}

\newcommand{\E}{\mathbb{E}}
\newcommand{\T}{\mathbb{T}}

\newcommand{\N}{\mathbb{N}}
\newcommand{\Z}{\mathbb{Z}}
\newcommand{\G}{\mathbb{G}}
\newcommand{\s}{\mathbf{s}}

\newcommand{\cz}{\mathbf{c_0}}
\newcommand{\co}{\mathbf{c_1}}
\newcommand{\bt}{\mathbf{t}}
\newcommand{\p}{\mu}
\newcommand{\den}{\mu}
\newcommand{\ns}{(1-\mu)}

\newcommand{\A}{\mathcal{A}}
\newcommand{\M}{\mathcal{M}}
\newcommand{\cE}{\mathcal{E}}
\newcommand{\oT}{\overline{\mathbb{T}}}

\newcommand{\Fo}{\operatorname{FPP}_1}
\newcommand{\Fl}{\operatorname{FPP}_\lambda}
\newcommand{\cin}{\mathbf{c}_{\operatorname{in}}}
\newcommand{\cout}{\mathbf{c}_{\operatorname{out}}}
\newcommand{\Wup}{W_{\text{up}}}
\newcommand{\Wlow}{W_{\text{low}}}
\newcommand{\Oup}{O_{\text{up}}}
\newcommand{\Olow}{O_{\text{low}}}
\renewcommand{\P}{\mathbb{P}}

\newcommand{\f}{\mathbf{f}}

\newcommand{\bseta}{\boldsymbol{\eta}}
\newcommand{\wtP}{\widetilde{\mathbb{P}}}
\newcommand{\wtE}{\widetilde{\mathbb{E}}}
\newcommand{\kun}{\underline{\kappa}}
\newcommand{\kov}{\overline{\kappa}}

\begin{document}



\section{Introduction}

We consider a natural growth process known as \emph{First passage percolation in a hostile environment} (FPPHE), which  
consists of two first passage percolation processes that spread at different rates, which we shall call $\Fo$ and $\Fl$.
Such processes are competing against each other to infect vertices of the underlying graph, according to the following dynamics.
Consider a connected infinite graph of bounded degree, choose a reference point on it (call it \emph{the origin}), and place a $\Fo$ particle there.
Then fix a parameter $\mu\in (0,1)$, and at each vertex other than the origin independently place a $\Fl$ particle with probability $\mu$ and nothing with probability $1-\mu$.
We shall call such initial $\Fl$ particles \emph{seeds}.

At time zero, all seeds of $\Fl$ are \emph{inactive}, and $\Fo$ starts spreading from the origin as a first passage percolation process of rate $1$. 
Whenever there is an attempt to occupy a site that is currently occupied by an inactive seed, then the seed gets activated and starts to spread as a first passage percolation of 
rate $\lambda$. Whenever a site gets occupied by one of the two processes, it will remain so forever. Thus the two processes compete with one another for space.
We see $\Fo$ and $\Fl$ as infection processes, so we will say that a vertex is \emph{infected} by $\Fo$ or $\Fl$ to mean that the vertex is occupied by the given process.

FPPHE is a natural model for the growth of two processes with competition and, as we will see below, it has many interesting aspects and phase transitions.
Nonetheless, it was introduced by Sidoravicius and Stauffer~\cite{Stauffer-Sidoravicius-MDLA} as an auxiliary model for analyzing a notoriously challenging model called Multiparticle Diffusion Limited Aggregation (MDLA).
The use of FPPHE as a tool of analysis has been further explored in Dauvergne and Sly~\cite{Dauvergne-Sly} to analyze a non-equilibrium process for the spread of an infection. 
The way in which FPPHE is used to analyze such processes with non-equilibrium dynamics is by making the growth of $\Fo$ represent the propagation of the front as it spreads through ``typically good regions'', 
while the activation of a $\Fl$ seed represents places where 
the front behaves ``atypically bad''. The influence of discovering such bad events is represented by the spread of $\Fl$ from the activated seed. 
In this context, one is interested to understand whether $\Fo$ grows indefinitely and leaves behind any cluster of activated seeds of $\Fl$.

There are several challenging questions about FPPHE.
One example is the \emph{survival} of $\Fo$ (resp.\ $\Fl$), 
meaning that $\Fo$ (resp.\ $\Fl$) infects an infinite connected region of the graph.
Note that it is natural to require connectivity of the infected region in the definition of survival since $\Fl$ already starts from an infinite set of seeds.
Another fundamental question is that of \emph{coexistence} of $\Fo$ and $\Fl$, which corresponds to both processes surviving simultaneously.

The first paper about FPPHE~\cite{Stauffer-Sidoravicius-MDLA} established that on $\Z^d$ if $\lambda<1$, then for all $\mu$ small enough $\Fo$ survives with positive probability and produces a \emph{bulky} shape (cf.\ \cite[Theorems 1.3]{Stauffer-Sidoravicius-MDLA}). 
Simple arguments show that, on $\Z^d$, $\Fo$ dies out almost surely if $\lambda\geq 1$ or $\mu$ large enough.
The justification of the first case can be found in Finn's PhD thesis \cite[Chapter 4]{Finn-PhD}.
The latter case follows from an percolation argument.
Very recently, \cite{Finn-Stauffer-FPPHE} showed the opposite direction that $\Fo$ survives with positive probability for any $\mu$ for which non-seeds percolate on $\Z^d$ provided $\lambda$ is small enough.
This seems analogous to the work of \cite{Stauffer-Sidoravicius-MDLA}, but actually FPPHE behaves rather differently in this setting. For example,  
the result of~\cite{Finn-Stauffer-FPPHE} implies a coexistence regime on $\Z^d$ 
for $d\geq 3$, and even a regime where both types occupy a set of positive density.
Coexistence is a challenging phenomenon to establish. 
Even in more classical models, such as the two-type Richardson model, the regime of coexistence has not been fully determined (cf.~\cite{Haggstrom-Pemantle}, and \cite{GM2008} for a more general analysis of the model, not restricted to exponential passage times).
For FPPHE on $\Z^d$, it remains entirely open whether coexistence occurs when $d=2$. 
Before that, coexistence for FPPHE was established on hyperbolic non-amenable graphs by~\cite{Candellero-Stauffer-FppCompetition}, which also showed that FPPHE behaves quite different in such graphs when compared to $\Z^d$. 
In fact, \cite{Candellero-Stauffer-FppCompetition} shows that $\Fl$ survives almost surely, for any fixed pair of values for the parameters $0<\mu\leq 1$ and $\lambda>0$.

The core of this work will focus on the behavior of FPPHE on a type of graph that is hyperbolic, non-amenable, quasi-transitive and is endowed with a ``tree-like structure''. 
However, our result does \emph{not} apply to trees, because of the lack of loops in their structure.
Although this construction naturally leads to a specific kind of graph, we are confident that the main ideas developed in this paper will be applicable to other contexts, at least to other growth models.
A work that follows a similar flavor is \cite{Deijfen-nonmonotonic}, in which the authors investigate competing Richardson processes on graphs.
They construct specific graphs where the set of values for $\lambda$ that allow coexistence has very counter-intuitive properties.
However, a crucial difference is that the Richardson model is a monotone model, in the sense that increasing the rate of one type can only increase the probability that said type survives. 
The focus of \cite{Deijfen-nonmonotonic} is on the event of coexistence, while we are interested in studying monotonicity properties of the FPPHE model on its own.

It is natural to believe that by increasing the density $\mu$ or the rate $\lambda$ one would favor $\Fl$.
Such monotonicity property would be quite useful as it plays a crucial role in the analysis of several other related models, 
such as the contact process and the two-type Richardson model \cite{liggett-SISys,Haggstrom-Pemantle,HP,GM2005,Hoffman}.
However, it is possible to show that there are choices for the locations of the seeds and for the passage times such that the standard coupling fails to be monotone.
In particular, on graphs whose structure present a lot of loops of different lengths, such configurations can be engineered so that adding a seed is beneficial to $\Fo$. The reason is that when the seed is activated, 
the process gets slowed down locally, and this can delay the activation of seeds that are further away in that direction. 
This could be beneficial for $\Fo$ as it grows along nearby directions.

The lack of knowledge of such monotonicity gives rise to serious challenges; the papers \cite{Stauffer-Sidoravicius-MDLA, Candellero-Stauffer-FppCompetition, Finn-Stauffer-FPPHE} develop quite involved multi-scale arguments to analyze FPPHE without using 
monotonicity.
The aim of the present article is to show that FPPHE is indeed \emph{not} monotone in general, meaning that there are graphs for which the probability of survival of $\Fo$ is not a monotone function of the parameters.
The graph that we obtain is rather different from a $d$-dimensional lattice; roughly speaking, it looks like a homogeneous tree of large degree, where all edges are replaced by identical copies of a finite graph that we call \emph{tile}.
(A formal definition is given in Section \ref{sect:ideas-proof}).

We believe that on regular graphs such as $\Z^d$ ($d\geq 1$) the process should be monotone.
However, the challenges in this type of analysis stem from several facts.
For example, 
geodesics on $\Z^d$ can go very far apart from each other, making it more difficult to control the spread of FPPHE.
Moreover, there are arbitrarily long loops, allowing $\Fo$ (resp.\ $\Fl$) to ``surround'' arbitrarily large clusters of $\Fl$ (resp.\ $\Fo$) within a short time.
At the moment a general result still eludes us, but we hope that this work will trigger a more general analysis of this question.

%
%


For any event $E$ and considering $\lambda$ fixed, let $\P_{\mu}(E)$ denote the probability of $E$, given that the initial density of seeds of $\Fl$ equals $\mu$.

\begin{theorem}\label{thm:main1}
There is a connected, infinite graph of bounded degree and values $\mu_1,\mu_2$ with $0<\mu_1<\mu_2<1$, such that, whenever the rate $\lambda$ is small enough,
\[
\P_{\mu_1}(\Fo \text{ survives})=0 \quad \text{and} \quad \P_{\mu_2}(\Fo \text{ survives})>0.
\]
\end{theorem}

\begin{remark}
   In our proof of Theorem~\ref{thm:main1}, we will construct a graph (which we denote $\G^\infty$) with the above properties.
   We note that $\G^\infty$ is even a \emph{quasi transitive} graph. Recall that a graph is quasi transitive if its vertex set can be partitioned into a finite number of sets $V_1,V_2,\ldots,V_\ell$ such 
   that for any $k$ and any two vertices $u,v\in V_k$, there exists an automorphism of the graph that maps $u$ onto $v$. Theorem~\ref{thm:main1} then implies at least three phase transitions as $\lambda$ is fixed 
   and $\mu$ varies. It is possible to get more involved constructions that lead to a larger number of phase transitions, and a similar non-monotonicity result can be obtained for $\mu$ fixed and $\lambda$ varying (see Theorem~\ref{thm:main2} and 
   Remark~\ref{rem:main2}).
\end{remark}

In our analysis of whether $\Fo$ survives we will use the following classical fact regarding percolation on regular trees (cf.\ \cite[Chapter 5]{LyonsPeres-book} for more details).

\begin{theorem}{\cite[Theorem 5.15]{LyonsPeres-book}.}\label{thm:lyons-peres-percolation} 
Let $T$ denote a homogeneous tree of degree $b\geq 3$. 
Then, 
\[
p_c(T)=\frac{1}{b-1},
\]
where $p_c(T)$ denotes the percolation threshold for the appearance of an infinite cluster for Bernoulli percolation on $T$. 
\end{theorem}
Note that in the case of trees the threshold for bond and for site percolation is the same.


\subsection{Idea of the proof}\label{sect:ideas-proof}
The proof proceeds through various steps, which we outline here.

\paragraph{Construction of the graph.}
The graph $\G^\infty$ should be thought of as a sort of infinite tree where each vertex has degree $\Phi+1$, where $\Phi$ is a large integer.
We replace the edges of the tree by identical copies of a finite graph that we call ``\emph{tile}''.
The structure of a tile plays a key role in determining which process gets (with reasonable probability) to advance more quickly. 
Below we give an overall view on the construction of the tiles, leaving the details to Section \ref{sect:construction-tile}.

The origin of the tile is a vertex from where two edges depart, we shall say that they begin the \emph{upper part} and the \emph{lower part} of the tile (refer to Figure \ref{fig:graph-lenghts}).
Following the lower part, we find that the endpoint of the first edge is connected with a finite binary tree of some (large) height $H$. 
The last generation of this tree is then identified into one vertex, which represents the start of a path of length $R$ (for some large integer $R$). 
%
%
Following the upper part, we see that the endpoint of the edge departing from the origin is connected with a finite $D$-ary tree (for some large integer $D>2$) of some height $L$.
Again, the last generation of the $D$-ary tree is identified into one vertex, to which one edge is attached.
Finally, the two loose ends that we are left with (the end of the path in the lower part and the endpoint of the edge in the upper part) are identified into one vertex.

\paragraph{Analysis of FPPHE for $\mu=\mu_1$ and $\mu=\mu_2$.}
The next step consists in investigating the typical behavior of FPPHE on a tile, hence we split the study of the process into two parts. 

In Section \ref{sect:calculations-mu1} we set the initial density of seeds to be some properly chosen value $\mu=\mu_1\in (0,1/2)$, and we show that the first process to infect the opposite endpoint of the tile is likely to be $\Fl$.
The idea is that for small values of $\mu$, the seeds do not block the advance of $\Fo$ through the two trees (the $D$-ary tree in the upper part and the binary tree in the lower part).
Therefore it is likely that $\Fo$ is the first process to reach the end of both trees.
By letting $H$ much smaller than $L$, we obtain that $\Fo$ reaches the end of the binary tree first.
At this point though, in the lower part it is likely that $\Fo$ will encounter a seed along the path of length $R$.
Hence we will have that the endpoint of the tile can be infected either by $\Fo$ (which is making its way along the upper part of the tile) or by $\Fl$ (started by a seed activated in the lower part).
The construction of the tile, in particular by properly choosing the height of the $D$-ary tree and the binary tree, 
makes it more likely that $\Fl$ (through the lower part) arrives first.

We note that implementing the idea above requires a delicate analysis of FPPHE on a finite tree. The main challenge is that even if $\Fo$ spreads quickly along a tree it may so happen that, near the last 
levels of the tree, seeds are located on the fastest route of $\Fo$. Therefore, it is not enough to look at the geodesics of first passage percolation on trees: we also need to control paths 
that are ``near optimal'' for $\Fo$. We do this via a comparison with branching random walks, for which strong concentration results are available.
This is the most technical part of the paper. It is established in Proposition~\ref{prop:Fo-wins}, whose proof is given in Section \ref{sect:Fo-wins}.

The second part is carried out in Section \ref{sect:calculations-mu2}, where we set the initial density to some properly chosen value $\mu=\mu_2 \in (1/2, 1-1/D)$.
This will imply that while in the upper part $\Fo$ continues to be the process that reaches the last generation of the $D$-ary tree, 
in the lower part it will be likely that $\Fo$  will be blocked by the seeds in the binary tree (as a straightforward consequence of Theorem \ref{thm:lyons-peres-percolation}).
This delays the progress of the process through the lower part, making it likely that the first process to reach the opposite endpoint of the tile is $\Fo$ (through the upper part).

\paragraph{Conclusion of the proof.}
To conclude, recall that $\Phi$ was a large integer appearing in the construction of $\G^\infty$ as the ``degree'' of each vertex.
Then we show that the event 
\[
\mathcal{E}:=\{\text{the first process to infect the endpoint of the tile is } \Fo \}
\]
is such that  $\P_{\mu_1}(\mathcal{E})<\frac{1}{\Phi}$  whereas $ \P_{\mu_2}(\mathcal{E})>\frac{1}{\Phi}$.
Then, by viewing each tile where $\mathcal{E}$ holds as an open edge of a suitable Bernoulli percolation process on a $\Phi$-ary tree and using Theorem \ref{thm:lyons-peres-percolation}, we conclude the proof.

\paragraph{Organization of the paper.}
The paper is organized as follows.
Section \ref{sect:constructionANDproof} deals with the construction of $ \G^\infty$ with vertex set $V(\G^\infty)$ and edge set $E(\G^\infty)$ and states a result (Proposition \ref{prop:aux}) which is the main step of the proof, giving the typical behavior of $\Fo$ on a tile. 
Using Proposition \ref{prop:aux}, we conclude Section \ref{sect:constructionANDproof} with the proof of Theorem \ref{thm:main1}.
In Section \ref{sect:proofsOnTile} we prove auxiliary lemmas, and in Section \ref{sect:proof-prop-aux} we prove Proposition \ref{prop:aux} using these lemmas and Proposition \ref{prop:Fo-wins}. 
Section \ref{sect:Fo-wins} is  devoted to the proof of Proposition \ref{prop:Fo-wins}.

\subsection{Passage times}

A path $\pi$ of length $n$ is a sequence of $n$ consecutive edges $\{e_{\{0,1\}}, e_{\{1,2\}}, \ldots, e_{\{n-1,n\}}\} $.
Fix a value $\gamma>0$.
For every path $\pi=\{e_{\{0,1\}}, \ldots, e_{\{n-1,n\}}\}$ define its \emph{passage time} by
\[
T^{(\gamma)}(\pi):=\sum_{j=0}^{n-1} t^{(\gamma)}_{\{j, j+1\}},
\]
where $t^{(\gamma)}_{\{j, j+1\}}$ is 
a random variable associated to edge $e_{\{j,j+1\}}$ and distributed according to Exp$(\gamma)$.
For any given $\gamma$ the variables $t^{(\gamma)}_{\{j, j+1\}}$ are independent and are often called \emph{passage times} of the corresponding edges.
In particular, when we write $\Fo$ we mean first passage percolation with $\gamma=1$ whereas when we write $\Fl$ we mean first passage percolation with $\gamma=\lambda$.

\begin{definition}[First Passage Times]\label{def:FPT}
At this point we can define the \emph{passage time} between two vertices $x$ and $y$ of the graph as
\[
T^{(\gamma)}(x\to y):= \min_{\pi : \, x \to y} T^{(\gamma)}(\pi),
\]
where the minimum is taken over all paths $\pi$ that connect $x$ to $y$.
Throughout the paper we will also need a ``mixed'' first passage time. 
Given a path $\pi$ from a vertex $x$ to a vertex $y$, if there exists a  vertex of $\pi$ that hosts a seed, let $w$ be the first such vertex and 
let $k\leq n$ denote the distance between $x$ and $w$ along $\pi$; note that for fixed $x,y$ and $\pi$, the value $k$ is a random variable as it depends on whether vertices of $\pi$ host a seed.
Then the passage time of $\pi$ is defined as
\[
T^{(1, \lambda)}(\pi):=\sum_{j=0}^{k-1} t^{(1)}_{\{j, j+1\}} + \sum_{j=k}^{n-1} t^{(\lambda)}_{\{j, j+1\}};
\]
otherwise, if no vertex of $\pi$ hosts a seed simply define $T^{(1,\lambda)}(\pi):=T^{(1)}(\pi)$.
Furthermore, for any two vertices $x,y\in V(G)$ that are \emph{cut-points}\footnote{A vertex is a cut-point when, if removed, it disconnects the graph.} in $\G^\infty$, we also define the ``mixed'' passage time between two vertices $x$ and $y$ as
\[
T^{(1, \lambda)}(x\to y):= \min_{\pi : \, x \to y} T^{(1, \lambda)}(\pi),
\]
where the minimum is taken over all paths $\pi$ that connect $x$ to $y$.
Since $x$ and $y$ are cut-points in $\G^\infty$, there are only finitely many paths between $x$ and $y$.
\end{definition}

\begin{remark}\label{rem:coupling}
%
%
We can couple the passage times for different values of $\gamma$ as follows. 
For each edge $e$, let $t^{(1)}_e$ be independent exponential random variables of rate $1$. Then, for any $\gamma>0$, $\gamma\neq 1$ and any edge $e$, we set $t^{(\gamma)}_e := t^{(1)}_e/\gamma$.
\end{remark}

We will need a result that is based on the known correspondence between the Gamma distribution and a Poisson process.
For a complete proof the reader is referred to \cite[Lemma 2.5]{Candellero-Stauffer-FppCompetition}.
In the following statement we denote by $\P$ the probability measure defined on the considered graph where seeds are not present.
More specifically, we consider a FPP process with (generic) rate $\gamma>0$ started at a unique source.
We write the statement for a first passage percolation process of rate $\gamma>0$, as we shall be interested in applying this result to the values $\gamma=1$ for $\Fo$ and $\gamma=\lambda$ for $\Fl$.
For a more clear statement, in the following lemma we set $\gamma=1$, but the statement hold for arbitrary $\gamma>0$.
In fact, it suffices to observe that for any path $\pi$, we have that $T^{(1)}(\pi)$ has the same distribution as $\gamma T^{(\gamma)}(\pi)$.

\begin{lemma}\label{lemma:typical-FPP} {\cite[Lemma 2.5 + Remark 2.6]{Candellero-Stauffer-FppCompetition}}
Consider an infinite graph with bounded degree, and let $o$ denote a reference vertex.
For any constant $\co  >0$, there is a constant $\cout =\cout ( \co  ) >1$ so that for all $k\geq 1$
\begin{equation}\label{eq:Eq1}
\P\left [\exists \text{ a path }\pi_k \text{ of length $k$ started at $o$, such that } T^{(1)}(\pi_k)\leq \frac{k}{ \cout}\right ]\leq e^{-\co   k }.
\end{equation}
For any constant $0<\cz <1$ there is a constant $\cin=\cin(\cz )<1$ so that for all $k\geq 1$
\begin{equation}\label{eq:Eq2}
\P\left [ \exists \text{ a path }\pi_k\text{ of length $k$ started at $o$,  such that } T^{(1)}(\pi_k)\geq \frac{k}{ \cin}\right ]\leq e^{-\cz   k }.
\end{equation}
\end{lemma}

\section{Construction of the graph}\label{sect:constructionANDproof}


The construction of the sought graph is done by steps, we start by describing the ``tiles'' that we will use.
The procedure about how to ``attach'' them together in a suitable manner is explained in Section \ref{sect:construction-G-infinity}.

%
%

\subsection{The tile $G$: construction, notation and some facts}\label{sect:construction-tile}

\begin{figure}[h!]
\begin{center}
\includegraphics[width=11cm]{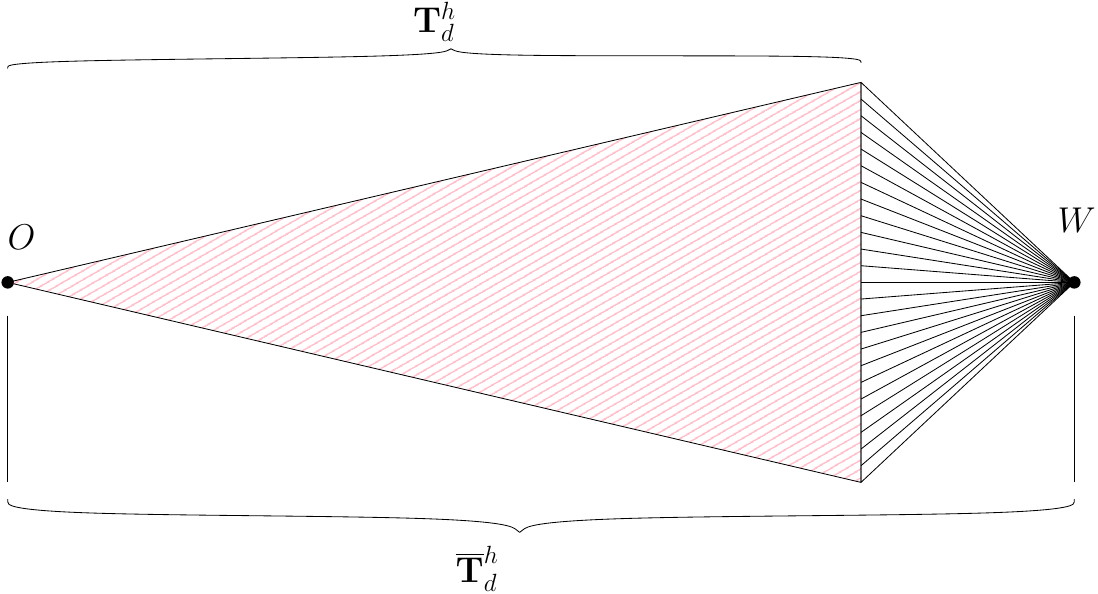}
\end{center}
\caption{Sketch of $\oT_d^h$: the triangle on the left (pink) represents the tree $\T_d^h$, while the one on the right (black) represents the set of edges connecting the $h$-th generation of $\T_d^h$ to $W$.
}\label{fig:treeTdh}
\end{figure}

Roughly speaking, a ``tile'' $G$ with vertex set $V(G)$ and edge set $E(G)$ is a connected graph consisting of different types of trees.
For any integer $d\geq 1$ let $\T_d$ be the infinite tree where all its nodes (including the root) have exactly $d$ children.
For any $h\in \mathbb{N}$ we let $\T_d^{h}$  denote the tree $\T_d$ truncated at generation $h$.
Finally, the notation $\oT_d^{h}$ denotes the tree $\T_d^{h+1}$ where all nodes in generation $h+1$ have been identified into one vertex (refer to Figure \ref{fig:treeTdh}, where such vertex is denoted by $W$).

\begin{figure}[h!]
\begin{center}
\includegraphics[width=10cm]{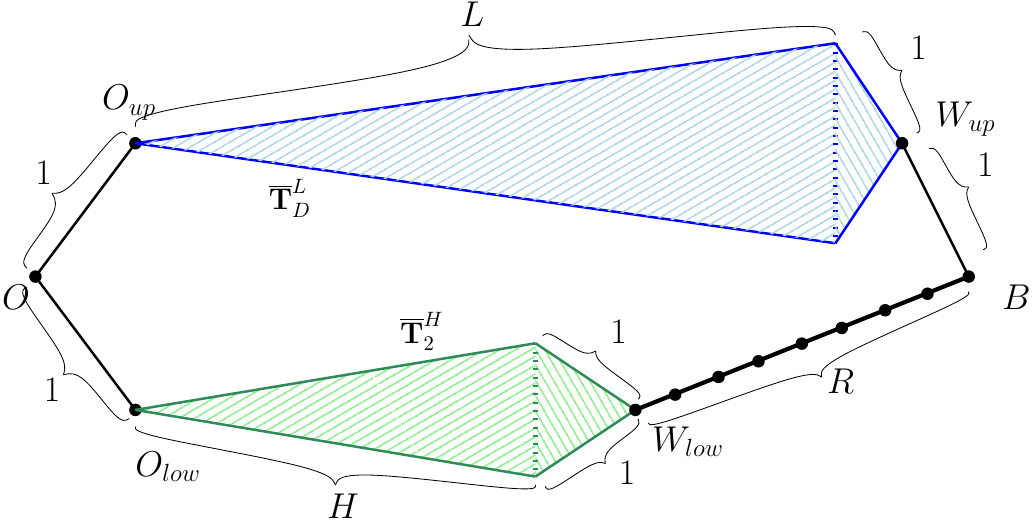}
\end{center}
\caption{Sketch of the tile $G$ and the lengths of its components.}\label{fig:graph-lenghts}
\end{figure}

For ease of readability, the reader can refer to Figure \ref{fig:graph-lenghts} for a sketch of $G$.
Start from a vertex $O$ (which we also call \emph{the origin}) and attach two edges to it, that is, identify the starting point of these two edges into $O$.
Call the other two endpoints $O_{\text{up}}$ and $O_{\text{low}}$.

For the moment we shall simply say ``fix a sufficiently large integer'', meaning that we shall be more precise later on by fixing the conditions that such value must satisfy.
Fix two sufficiently large integers $D>2$ and $L>0$, and identify the root of the graph $\oT_D^L$ with $O_{\text{up}}$.
Similarly, for another large integer $H>0$ identify the root of $\oT_2^H$ with $O_{\text{low}}$.
Basically, we are ``attaching'' $\oT_D^L$ to one edge incident to the origin, and $\oT_2^H$ to the other one.
The last vertex of $\oT_D^L$ will be denoted by $\Wup$ and the last vertex of $\oT_2^H$ will be denoted by $\Wlow$.

Finally, attach a single edge from $\Wup$ to a new vertex $B$ and a path of length $R$ (for a large integer $R$) from $\Wlow$ to $B$.
(Refer to Figure \ref{fig:graph-lenghts}.)


Recall the constants $\cz ,\co ,\cin$ and $\cout$ from Lemma \ref{lemma:typical-FPP}.
Later on, for $d\in \mathbb{N}$, we shall emphasize this dependence by writing
\begin{equation}\label{eq:cin-cout-d}
\begin{split}
& \cin^{(d)}:=\cin^{(d)}(\cz ):= \text{value of }\cin\text{ satisfying \eqref{eq:Eq2} when }\cz \text{ is given and the graph is } \T_d,\\
& \cout^{(d)}:=\cout^{(d)}(\co ):= \text{value of }\cout\text{ satisfying  \eqref{eq:Eq1} when }\co \text{ is given and the graph is }\T_d.
\end{split}
\end{equation}
For the particular case where $d=1$, we write $\cin^{(1)}$ and $\cout^{(1)}$ to denote the values that on a large interval of $\Z$ satisfy \eqref{eq:Eq2} and \eqref{eq:Eq1} respectively. 

Before proceeding, we introduce some notation.
%
For any fixed $\mu\in (0,1)$, any fixed $d\in \mathbb{N}$ and any infinite tree $\T_d$ define the event
\[
\cE_{\text{perc}}(\T_d):= \{ \text{cluster of non-seeds containing  the root percolates on }\T_d\},
\]
then set $\f_d(\mu)$ to be
\begin{equation}\label{eq:defBFf}
\f_d(\mu):=1-\P_\mu \left ( \cE_{\text{perc}}(\T_d)\right ).
\end{equation}
By monotonicity of percolation, we know that $\f_d(\mu)$ is increasing with $\mu$.
%
To simplify the notation throughout, for any vertex $x\in V(G)$ and $\mu \in (0,1)$ fixed, set
\begin{equation}\label{eq:def-Ex-mu}
\cE_x:=\{\text{vertex }x\text{ does not host a seed}\},
\end{equation}
for which we have 
\begin{equation}\label{eq:prob(1-mu)}
\P_\mu(\cE_x)=1-\mu. 
\end{equation}
On $ \oT_d^{h}$ with the notation borrowed from Figure \ref{fig:treeTdh} we set
\begin{equation}\label{eq:def-perc}
\cE_{\text{perc}}(\oT_d^{h}):=  \{\text{there is a cluster of non-seeds in $\oT_d^{h}$ connecting }O\text{ to }W\},
\end{equation}
for which we have
\begin{equation}\label{eq:prob(1-f)}
\P_\mu \left [ \cE_{\text{perc}}(\oT_d^{h})\right ]\geq \bigl ( 1-\f_d(\mu)\bigr )(1-\mu).
\end{equation}
In particular, the first term of the product on the R.H.S.\ follows from \eqref{eq:defBFf} and the second is the probability that vertex $W$ does not host a seed.
(Recall that in the definition of the model the origin is defined to be free of seeds.)

\subsection{Behavior of FPPHE on a tile}
%


Fix $\lambda>0$ and consider a copy of the tile $G$ as constructed above, and start FPPHE by placing a $\Fo$ particle at $O$, and a seed independently at each vertex with probability $\mu$.
From now on, the value of $\lambda>0$ will be considered fixed and \emph{small enough} for our purposes.


Now we are able to state a result that will be crucial for the whole analysis.
The considered graph is $\oT_d^{m}$, and the notation that we use is the same as in Figure \ref{fig:treeTdh}.


\begin{proposition}\label{prop:Fo-wins}
Fix an integer $d\geq 2$, take $\lambda>0$ small enough, and any $\mu \in (0,1)$.
For all $m\in \N$ large enough, consider $\oT_d^{m}$ with the following initial setting:
place a particle of $\Fo$ at $o $ and at every $x\in V(\oT_d^{m})\setminus \{o \}$  place seeds independently with probability $\mu$. 
Then there is a function $\bseta_d=\bseta_d(\mu, \lambda)>0$ so that 
the following hold.
\begin{itemize}
\item[(i)] For all $\mu$ small enough, 
\[
\P_{\mu}
\Bigl [ W\text{ is infected by }\Fo \Bigr ]\geq \bigl (1 - \f_d(\mu)-\bseta_d\bigr )(1-\mu).
\]
\item[(ii)] If $\mu\in (1/2, 1)$ so that 
\begin{equation}\label{eq:techCondMu}
d\ns>1 \quad \text{ and } \quad d^2(1-\mu)^2\mu^{d-1}<1,
\end{equation}
then,
\[
\P_{\mu} \Bigl [ W\text{ is infected by }\Fo \Bigr ]\geq \bigl (1 - \f_d(\mu)-\bseta_d\bigr )(1-\mu) .
\]
\item[(iii)]  For any $\xi>0$, if $\lambda$ is small enough and either $\mu$ is small enough or $d$ is large enough, then $\bseta_d\leq \xi$.
\end{itemize}
\end{proposition}

The factor $1-\mu$ corresponds to the probability that $W$ does not host a seed. 

\begin{remark}
Observe that, since $\lambda>0$ is fixed throughout, the quantity $\bseta_d (\mu, \lambda)$ can be seen as a function of $\mu$ only.
Furthermore, whenever $\lambda$ and $\mu$ are fixed, $\bseta_d (\mu, \lambda)$ is decreasing in $d$.
In order to avoid making the notation even more cumbersome, we chose that whenever $\mu$ is clear from the context we shall omit its explicit dependence from $\bseta_d$.
\end{remark}
Since the proof of Proposition \ref{prop:Fo-wins} is quite technical, we postpone it to Section \ref{sect:Fo-wins}, thus we will proceed while assuming the validity of the result.
The next statement provides a crucial fact regarding the behavior of the process on a tile $G$ with origin $O$.
\begin{proposition}\label{prop:aux}
Fix a large integer $D$, and let $\lambda$ be small enough.
Then there are values 
$0<\mu_{1,\text{low}}^{(\text{crit})}< \mu_{1,\text{high}}^{(\text{crit})} <1/2 < \mu_{2,\text{low}}^{(\text{crit})} <  \mu_{2,\text{high}}^{(\text{crit})} <1 $, and values of $H, L, R,\Phi$
so that in such a tile $G$ we obtain
\[
\P_{\mu} \left ( B\text{ is infected by }\Fo \right )<\frac{1}{\Phi} , \quad \text{ for all }\mu \in\left ( \mu_{1,\text{low}}^{(\text{crit})}, \mu_{1,\text{high}}^{(\text{crit})}\right )
\]
and
\[
\P_{\mu} 
\left ( B\text{ is infected by }\Fo \right )>\frac{1}{\Phi}, \quad \text{ for all }\mu \in \left ( \mu_{2,\text{low}}^{(\text{crit})},  \mu_{2,\text{high}}^{(\text{crit})}\right ).
\]
\end{proposition}
This statement will be proven in Section \ref{sect:proof-prop-aux}, admitting the result of Proposition \ref{prop:Fo-wins} (which will be proven in Section \ref{sect:Fo-wins}).

\subsection{The infinite graph $\G^\infty$}\label{sect:construction-G-infinity}
We proceed with the construction that provides the sought counterexample to monotonicity. 
We view each tile with an orientation, where the vertex $O$ is called the \emph{origin} of the tile, and the vertex $B$ is called the \emph{tail} of the tile.
Then $\G^\infty$ is constructed as follows.
Fix a large integer $\Phi$ and consider an infinite $(\Phi+1)$-regular tree. 
Consider an infinite (self-avoiding) path starting at the root, this determines an end of the tree, and thus a specific direction.
Now, for each vertex $v$ of the tree, let $u_0,u_1,\ldots,u_\Phi$ be the neighbors of $v$ with $u_0$ being the unique neighbor that is on the direction of the chosen end. Now we replace each edge $(v,u_i)$ of the tree, for $i\in\{1,2,\ldots,\Phi\}$, with a copy of the tile, 
identifying the origin of the tile (i.e., the vertex $O$) with $v$ and the tail of the tile (i.e., the vertex $B$) with $u_i$. For the edge $(v,u_0)$, the orientation of the tile is reversed; i.e., 
we replace $(v,u_0)$ with a copy of the tile, but identify the origin of the tile with $u_0$ and the tail of the tile with $v$. 

A graphical representation is shown in Figure \ref{fig:G-infty} (for simplicity we used $\Phi=2$, although it follows from the computations in Section \ref{sect:proof-prop-aux} that $\Phi$ is much larger than $2$).
Note that the graph $\G^\infty$ is quasi transitive.


\subsection{Proof of Theorem \ref{thm:main1}}

We state our result for FPPHE on the graph $\G^\infty$ below, from which Theorem \ref{thm:main1} follows directly.
\begin{theorem}\label{thm:main2}
Consider FPPHE starting from the origin of $\G^\infty$.
Fix a large integer $D$, let $\lambda$ be small enough, and let the values $\mu_{1,\text{low}}^{(\text{crit})}$, $\mu_{1,\text{high}}^{(\text{crit})}$, $ \mu_{2,\text{low}}^{(\text{crit})} $, $ \mu_{2,\text{high}}^{(\text{crit})}$ be as in Proposition \ref{prop:aux}.
Then there are values for $H, L, R$ and $\Phi$ so that, 
\[
\P_{\mu}(\Fo \text{ survives})=0 \quad \text{ for all }
\mu \in  \left (\mu_{1,\text{low}}^{(\text{crit})}, \, \mu_{1,\text{high}}^{(\text{crit})} \right ) 
\]
and
\[
\P_{\mu}(\Fo \text{ survives})>0 \quad \text{ for all }
\mu \in 
\left (\mu_{2,\text{low}}^{(\text{crit})}, \mu_{2,\text{high}}^{(\text{crit})} \right ).
\]
\end{theorem}

\begin{figure}[h!]
   \begin{center}
   \includegraphics[width=11cm]{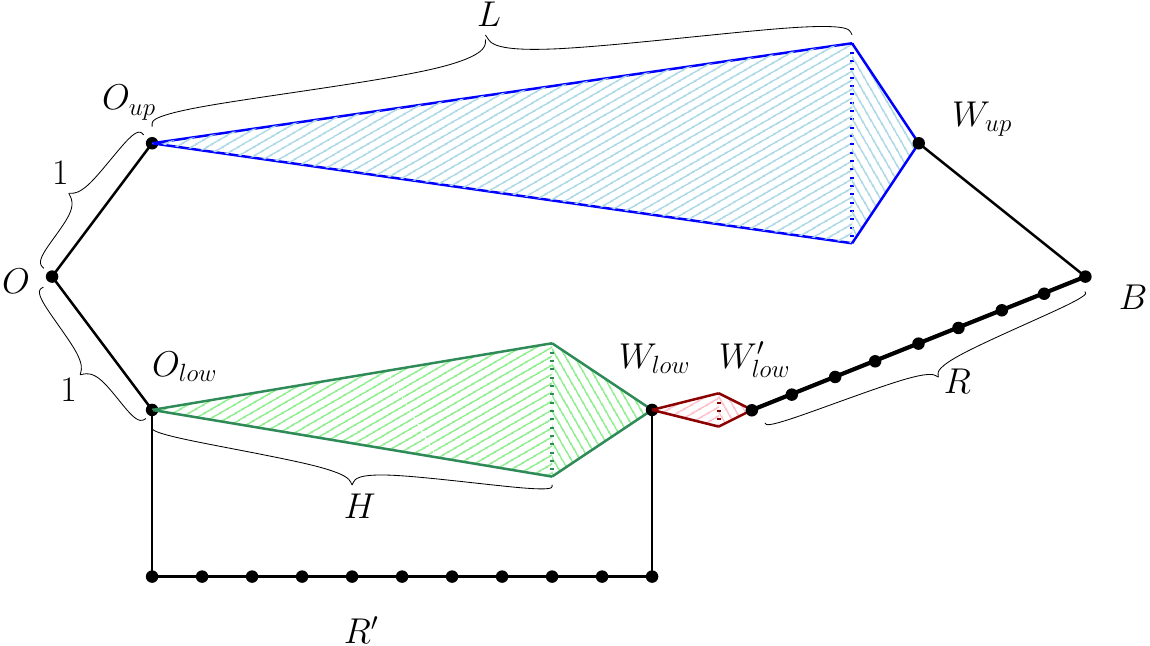}
   \end{center}
   \caption{The tile $G$ that gives non-monotonicity in $\lambda$ as in Remark~\ref{rem:main2} part \ref{it:main2}.}\label{fig:graph-lenghts3}
\end{figure}

\begin{remark}\label{rem:main2}
   We make the following observations.
   \begin{enumerate}
      \item It is immediate that $\Fo$ survives if $\mu$ is small enough (much smaller than $\mu_{1,\text{low}}^{(\text{crit})}$) since seed-free tiles percolate on $\G^\infty$ for all small enough $\mu$. 
      Moreover, it is immediate that $\Fo$ dies out if $\mu$ is large enough (much larger than $\mu_{2,\text{high}}^{(\text{crit})}$) so that the origin is blocked by seeds with probability one. 
      Therefore, this reasoning and Theorem~\ref{thm:main2} establish at least three phase transitions on $\G^\infty$: from survival of $\Fo$ for small enough $\mu$, then non-survival for some $\mu$ in $(0,1/2)$, 
      then a survival regime for some $\mu$ in $(1/2,1)$ to a final regime of non-survival for all large $\mu$. 
      
      \item By further working the construction of the tiles we can make as many phase transitions as we desire. 
         For example, we could have a tile composed of $\kappa$ copies of the upper and lower parts so that $O$ would have degree 
         $2\kappa$, with each copy of the upper and lower parts of different sizes and different degrees for the corresponding trees. Then by tuning the parameters properly, as we increase $\mu$ 
         we could make $B$ be occupied by $\Fl$ through the first lower part, then by $\Fo$ through the first upper part, then by $\Fl$ through the second lower part, then 
         by $\Fo$ through the second upper part, and so on and so forth.
      
      \item\label{it:main2} We choose to fix $\lambda$ small enough and show non-monotonicity when varying $\mu$. A similar result could be obtained by fixing $\mu$ small enough first and then varying $\lambda$ inside the interval $(0,1)$. 
      For example, in the lower part of the tile we could split the vertex $\Wlow$ into two vertices ($\Wlow$ and $\Wlow'$), connect them with a small finite tree, 
      and then create an additional path of length $R'$ connecting $\Olow$ to $\Wlow$ (see Figure~\ref{fig:graph-lenghts3}). 
      Then, from $\Olow$ there would be two ways to $\Wlow$, one through
      $\oT_2^H$ which is favorable to $\Fo$ by having $\mu$ small and $\lambda<1$, and the new path which will be likely to have a seed in it. 
      Even for not-too-small values of $\lambda$, we can make the path of length $R'$ slow enough to let $\Fo$ get to $\Wlow$ through $\oT_2^H$, leading to a fast passage through the small tree from $\Wlow$ to $\Wlow'$ 
      and the vertex $B$ to be 
      occupied by $\Fl$ from the lower part through the activation of a seed in the path of length $R$. 
      On the other direction, by increasing $\lambda$ we can make the path of length $R'$ to be the fastest
      way to $\Wlow$, but since $\lambda<1$ this makes the way through the small tree from $\Wlow$ to $\Wlow'$ slow enough to let the upper part of the tile be the fastest. 
      This leads to $\Fo$ occupying $B$ from the upper part. By tuning the parameters we can then establish the non-monotonicity. For conciseness, we only do the details for the phase transition on $\mu$.
   \end{enumerate}
\end{remark}


\begin{figure}[h!]
\begin{center}
\includegraphics[width=8cm]{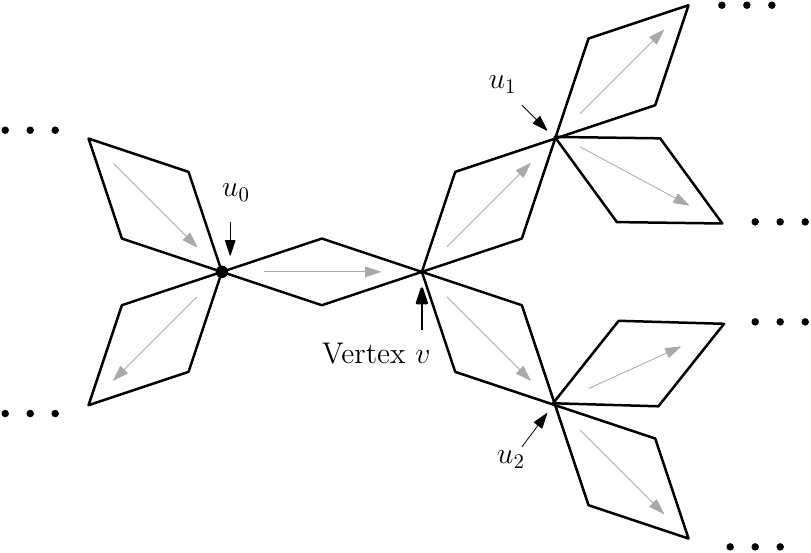}
\end{center}
\caption{The construction of $\G^\infty$, every diamond represents a copy of the tile $G$, while arrows inside each tile point from the origin of the tile to the tail. Considering vertex $v$ as the origin of $\G^\infty$, 
we obtain that $\vec\G^\infty$  (as formally defined in the proof of Theorems \ref{thm:main1} and~\ref{thm:main2}) as the part of the figure to the right of $v$.}\label{fig:G-infty}
\end{figure}

\begin{remark}
In the following, we set the parameters in the following order.
First, we choose $\mu_{1,\text{low}}^{(\text{crit})}, \mu_{1,\text{high}}^{(\text{crit})}$ and $\mu_{2,\text{low}}^{(\text{crit})}, \mu_{2,\text{high}}^{(\text{crit})}$, then we choose the values of $L,R,H, D$ uniformly over $\mu_1\in \left ( \mu_{1,\text{low}}^{(\text{crit})}, \mu_{1,\text{high}}^{(\text{crit})}\right )$ and $\mu_2\in \left ( \mu_{2,\text{low}}^{(\text{crit})}, \mu_{2,\text{high}}^{(\text{crit})}\right )$.
\end{remark}

\begin{proof}[Proof of Theorems \ref{thm:main1} and~\ref{thm:main2}]
We now prove Theorem~\ref{thm:main2} (Theorem~\ref{thm:main1} is just a direct consequence of it).
The proof of Theorem~\ref{thm:main2} is a consequence of Proposition \ref{prop:aux} and classical facts about percolation.
Proposition \ref{prop:aux} gives values for the quantities $D,L,H,R$ so that all calculations work for all $\mu_1 \in \left ( \mu_{1,\text{low}}^{(\text{crit})}, \mu_{1,\text{high}}^{(\text{crit})}\right )$ and $\mu_2 \in \left ( \mu_{2,\text{low}}^{(\text{crit})}, \mu_{2,\text{high}}^{(\text{crit})} \right )$.
Thus, we can fix $\mu_1$ and $\mu_2$ in the above intervals.
We can assume, without loss of generality, that the origin of $\G^\infty$ is a vertex that has been identified with the origin of a tile (i.e., it is not an internal vertex of a tile). 
The reason why we can assume this is that, if the origin were an internal vertex of a tile, then with positive probability $\Fo$ would spread to both the origin of the tile and the tail of the tile 
(since there is a positive probability
that the tile contains no seed). This means that $\Fo$ has a positive probability of surviving in this case if and only if it has a positive probability of survival when the origin of $\G^\infty$ is not internal to a tile.

From now on, let $o$ denote the origin of $\G^\infty$, which is identified with the origin of $\Phi$ tiles and with the tail of one tile, we denote this last tile $G_\mathrm{parent}$. 
Referring to Figure \ref{fig:G-infty} with the vertex $v$ being $o$, we obtain that $G_\mathrm{parent}$ is the tile between $u_0$ and $v$.
We shall denote the corresponding percolation probability measure by $\P^p$, in order to avoid confusion with the previously defined measures.
Now let $\vec\G^\infty$ be the subgraph where we remove from $\G^\infty$ all vertices and edges whose geodesic to $o$ intersects $G_\mathrm{parent}$; in Figure \ref{fig:G-infty} this is the part of the graph to the 
right of $v$. 
Note that by collapsing each tile of $\vec\G^\infty$ to an edge, we obtain a $\Phi$-ary tree rooted at $o$.
Then, by performing (edge) Bernoulli percolation with parameter $p$ on an infinite $\Phi$-ary tree we have (cf., e.g., \cite{LyonsPeres-book})
\[
\P^p\bigl (\text{cluster containing the root is infinite}\bigr )>0 \quad \Leftrightarrow \quad p>\frac{1}{\Phi}. 
\]
In a similar fashion, consider a percolation process on the graph $\vec\G^\infty$ defined as follows.
Fix a value $p\in (0,1)$.
For every copy of $G$ in $\vec\G^\infty$, toss an independent coin which gives heads with probability $p$ and tails with probability $1-p$.
If the coin gives heads, then declare the corresponding copy of $G$ to be \emph{open}, otherwise declare it \emph{closed}.
This process on $\vec\G^\infty$ is a percolation process of parameter $p$ on a $\Phi$-ary tree, and thus 
\begin{equation}\label{eq:Pp-cluster}
\P^p(\text{cluster containing vertex $o$ is infinite} )>0\quad \Leftrightarrow \quad p>\frac{1}{\Phi}.
\end{equation}
%
%
From Proposition \ref{prop:aux} it follows immediately that, 
for a large integer $\Phi$ (which will be specified later in Equation \eqref{eq:defPhi}) we have 
\begin{equation}\label{eq:p1p2}
\P_{\mu_1} \left (  B\text{ is infected by }\Fo \right )>\frac{1}{\Phi} \quad \text{ and } \quad \P_{\mu_2} \left (  B\text{ is infected by }\Fo \right )<\frac{1}{\Phi}.
\end{equation}
Therefore, it suffices to realize that the event $ \{ \Fo $ survives$\}$ corresponds to the event
\[
\cE := \left \{ \exists \text{ infinite cluster of tiles so that, for each such tile, vertex $B$ is infected by }\Fo\right \}.
\]
Hence, by \eqref{eq:Pp-cluster}
\[
\P^p\bigl (\cE \bigr )>0\quad \Leftrightarrow \quad 
p>\frac{1}{\Phi}.
\]

Now we turn to analyzing $\G^\infty$. When $\Fo \text{ survives}$ in $\vec\G^\infty$, then it also survives on $\G^\infty$. We are then left to show that when $\Fo$ dies out almost surely in $\vec\G^\infty$ then 
it also dies out almost surely in $\G^\infty$. To see this, let $o_1,o_2,\ldots$ be a sequence of origins of tiles, starting with $o=o_1$, such that for each $i\geq 1$ 
there is a tile for which $o_{i+1}$ is its origin and $o_i$ its tail.
In other words, $o_1,o_2,\ldots$ is the unique sequence of origins of tiles that are traversed by a path that starts from $o$ and can only traverse tiles from the tail to the origin. Since this sequence is unique and the corresponding 
tiles form a one-dimensional chain of tiles, almost surely $\Fo$ occupies only a finite number of elements in this sequence (for example, it will happen almost surely that there exists a $k$ for which 
$o_k$ hosts a seed). 
Let $\ell<\infty$ be such that $o_\ell$ is occupied by $\Fo$ and $o_{\ell+1}$ is occupied by $\Fl$. For $\Fo$ to survive, it has to survive by spreading through a path from $o_1,o_2,\ldots,o_\ell$ which first traverses a tile from the 
origin to the tail. But such paths correspond to spreading of $\Fo$ in $\vec\G^\infty$. Since $\Fo$ dies out almost surely in $\vec\G^\infty$, it will also die out when spreading from each one of $o_1,o_2,\ldots,o_\ell$, 
concluding the proof. 
\end{proof}

\section{Analysis of the process on the tile}\label{sect:proofsOnTile}

This section is devoted to the analysis of the process on each tile.
In particular, we will state and prove two results (Lemmas \ref{lemma:toProp-mu1} and \ref{lemma:toProp-mu2} below) which we will use in the proof of Proposition \ref{prop:aux}. 
In order to prove the lemmas, we will apply the result from Proposition \ref{prop:Fo-wins}. 
So $\mu_1$ and $\mu_2$ are given as specified in Proposition \ref{prop:Fo-wins}.
In the following, we shall refer to the \emph{upper part} of $G$ as the subgraph of $G$ induced by the vertices of $\oT_D^L$ together with $O$ and $B$.
Symmetrically, the \emph{lower part} of $G$ is the subgraph of $G$ induced by all vertices of $G$ except for those of $\oT_D^L $; cf.\ Figure \ref{fig:graph-lenghts}.

Informally speaking, we study the spread of the process along two main streams, the lower one and the upper one.
Since they only communicate at $O$ and at $B$, they are independent until the moment at which vertex $B$ gets infected.
Once this happens, it will be crucial to understand whether $B$ got infected from the upper or from the lower part.

In order to avoid confusion, we define another notation for passage times that are only local to $G$.
More precisely, if the edges of a path $\pi=(e_{\{0,1\}}\ldots e_{\{n-1,n\}})$ are completely contained in $E(G)$ and is such that when exploring the path from vertex $0$ to vertex $n$ the first seed is encountered at $k$, then (similarly to the reasoning in  Definition \ref{def:FPT}) we define the quantity  
\[
\Upsilon^{(1,\lambda)}(\pi):=\sum_{j=0}^{k-1} t^{(1)}_{\{j, j+1\}} + \sum_{j=k}^{n-1} t^{(\lambda)}_{\{j, j+1\}}.
\]
Consequently, similarly to what discussed in Definition \ref{def:FPT}, the ``mixed'' passage time between two vertices $x,y \in V(G)$ along a path completely contained in $G$
is denoted by $\Upsilon^{(1,\lambda)}(x \to y)$.
When $k=0$ or $k=n$ in the above formula, we shall simplify the notation writing $\Upsilon^{(\lambda)}(x \to y)$ or $\Upsilon^{(1)}(x \to y)$, respectively.
Note that this is the first passage time between $x$ and $y$ restricted to paths contained in $G$, and thus in general we would have $\Upsilon^{(1,\lambda)}(x \to y)\geq T^{(1,\lambda)}(x \to y)$.

To avoid confusion, we denote by  $\Upsilon^{(1,\lambda)}_{\text{up}}(x \to y)$ the first passage time between $x$ and $y$ restricted to paths in the \emph{upper part} of $G$, and by $\Upsilon^{(1,\lambda)}_{\text{low}}(x \to  y)$ the first passage time from $x$ to $y$ restricted to paths in the \emph{lower part}.
Moreover, when we write the event $\left \{ \Wup \text{ is infected by }\Fo \right \}$, we mean that this happens when we restrict to the upper part of the tile, so we do not account for the possibility that 
$\Wup$ is infected by $B$. 

We start from the fact that for any $\varepsilon>0$ fixed and any spread rate $\gamma>0$ there is a constant $C(\varepsilon, \gamma)>0$ such that 
uniformly in $e\in E(G)$ we have $\P\bigl (t_e^{(\gamma)}  < C(\varepsilon, \gamma) \bigr )>(1-\varepsilon)$.
In what follows, set $\varepsilon>0$ as an arbitrarily small fixed constant, and we will set
\begin{equation}\label{eq:def-cC}
\mathfrak{C}:=\max \bigl \{ C(\varepsilon, 1), \, C(\varepsilon, \lambda )\bigr \}.
\end{equation}
The next two subsections are devoted to prove two results.
More specifically, in Lemma \ref{lemma:toProp-mu1} we show that when $\mu$ is small, the process is very likely to be ``quick'' in the upper part and ``slow'' in the lower part (where ``quick'' and ``slow'' are to be made precise below).
Subsequently, in Lemma \ref{lemma:toProp-mu2} we show that when $\mu$ is large (but not too large), it is very likely that the situation is reversed.
The proof of Proposition \ref{prop:aux} will strongly rely on these results.


\begin{lemma}\label{lemma:toProp-mu1}
%
Let $\mathfrak{C}$ be as in \eqref{eq:def-cC}.
There are critical values $0<\mu_{1,\text{low}}^{(\text{crit})}< \mu_{1,\text{high}}^{(\text{crit})} <1/2 $ and large enough values of $L$, $H$, $R$ and $D$ so that for any fixed $\varepsilon>0$ there are 
$\cout^{(D)}= \cout^{(D)}(\varepsilon)$,  $\cin^{(1)}=\cin^{(1)}(\varepsilon)$, $\cin^{(2)}=\cin^{(2)}(\varepsilon)$ 
such that for all $\mu_1 \in \left (\mu_{1,\text{low}}^{(\text{crit})}, \, \mu_{1,\text{high}}^{(\text{crit})} \right )$ the following hold.
\begin{itemize}
\item[(i)] In the lower part,
\[
\begin{split}
\P_{\mu_1} & \left [ \Upsilon_{\text{low}}^{(1,\lambda)}(O\to B)\leq \mathfrak{C} + \frac{H+1}{\cin^{(2)}} + \frac{R}{\lambda \cin^{(1)}}\right ]\\
& \geq 
(1-\varepsilon)^3 (1-\mu_1)^2 \left [ (1-\f_2(\mu_1)-\bseta_2)(1-\mu_1)^2-\varepsilon \right ],
\end{split}
\]
where the factor $\bseta_2$ is as  in Proposition \ref{prop:Fo-wins} and $\f_2(\mu)$ as in \eqref{eq:defBFf} when $d=2$.
\item[(ii)] In the upper part,
\[
\P_{\mu_1} \left [ \Upsilon_{\text{up}}^{(1,\lambda)}(O\to B) \geq \frac{L+1}{\cout^{(D)}} \right ] \geq 1-\varepsilon.
\]
\end{itemize}
\end{lemma}

\paragraph*{Remark.}
The critical value $0<\mu_{1,\text{high}}^{(\text{crit})} <1/2$ depends  on $\lambda$ and has to satisfy some technical conditions that will be specified later on (expressed in details in conditions \eqref{eq:def-mu1-bar}, \eqref{eq:def-mu-bar-prime}, \eqref{eq:def-mu-doubleprime}, \eqref{eq:def-mu-triple}).
Moreover, since the result holds for all $\mu_1$ small enough, we can simply set $\mu_{1,\text{low}}^{(\text{crit})}:= \mu_{1,\text{high}}^{(\text{crit})} /2$.
Note that 
such values are independent of $L,H,R,D$.

\begin{lemma}\label{lemma:toProp-mu2}

Let $\mathfrak{C}$ be as in \eqref{eq:def-cC}.
There are critical values $1/2<\mu_{2,\text{low}}^{(\text{crit})}<\mu_{2,\text{high}}^{(\text{crit})}<1$ and 
large enough values of $L$, $H$, $R$ and $D$ so that for any fixed $\varepsilon>0$ there are 
$\cout^{(D)}= \cout^{(D)}(\varepsilon)$,  $\cin^{(1)}=\cin^{(1)}(\varepsilon)$, $\cin^{(2)}=\cin^{(2)}(\varepsilon)$ such that for all $\mu_2\in \left (\mu_{2,\text{low}}^{(\text{crit})}, \, \mu_{2,\text{high}}^{(\text{crit})} \right )$, 
%
%
\begin{itemize}
\item[(i)] In the lower part,
\[
\P_{\mu_2} \left [ \Upsilon_{\text{low}}^{(1,\lambda)}(O\to B)\geq \frac{H/2}{\cout^{(2)}}+\frac{H/2+1}{\lambda \cout^{(2)}}+ \frac{R}{\lambda \cout^{(1)}}\right ] \geq (1-3\varepsilon).
\]
\item[(ii)] In the upper part,
\[
\begin{split}
& \P_{\mu_2} \left [ \left \{ \Upsilon_{\text{up}}^{(1,\lambda)}(O\to B)\leq 2\mathfrak{C}+\frac{L+1}{\cin^{(D)}}\right \} \cap \left \{ \Wup \text{ is infected by }\Fo \right \} \right ]  \\
& \geq (1-\varepsilon)^2 (1-\mu_2)^2\left [ \bigl (1-\bseta_D-\f_D(\mu_2)\bigr )(1-\mu_2)^2-\varepsilon \right ] .
\end{split}
\]
where $\bseta_D$ is as in Proposition \ref{prop:Fo-wins} and $\f_D(\mu)$ as in \eqref{eq:defBFf} when $d=D$.
\end{itemize}
\end{lemma}

\subsection{Proof of Lemma \ref{lemma:toProp-mu1}}\label{sect:calculations-mu1}
The aim of this section is to understand the ``typical'' behavior of the process on the tile $G$ when $\mu=\mu_1$ is small enough.
%
Sections \ref{sect:mu_1_low} and \ref{sect:mu_1_up} take care of the analysis in the lower and upper part of $G$ respectively.

\subsubsection{Proof of Lemma \ref{lemma:toProp-mu1}(i), ``lower part'' of $G$.}\label{sect:mu_1_low}

Consider the variable $\Upsilon_{\text{low}}^{(1,\lambda)}(O\to B)$.
Fix any pair of vertices $x,y$ in the lower part of $G$ and an initial density of seeds $\mu_1\in \left (\mu_{1,\text{low}}^{(\text{crit})}, \, \mu_{1,\text{high}}^{(\text{crit})} \right )$, where $\mu_{1,\text{high}}^{(\text{crit})}$ is a small enough value (cf.\ Section \ref{sect:Fo-wins}).
For this part of the proof consider only the lower part of $G$, more precisely define the passage times $\{ \overline{t}_e^{(\gamma)}\}$ as follows.
Let $\overline{t}_e^{(\gamma)}:=t_e^{(\gamma)}$ if the edge $e$ is in the lower part of $G$, and set $\overline{t}_e^{(\gamma)}:=\infty$, if $e$ is in the upper part.
With respect to the passage times $\{\overline{t}_e^{(\gamma)}\}_{e\in E(G)}$, recalling the definition of $\cE_{y}$ from \eqref{eq:def-Ex-mu}, define the event
\begin{equation}\label{eq:eventB}
\mathfrak{B}_{\text{low}}(x,y):= \cE_{y}\cap \{ \Fo \text{ started at }x\text{ infects }y\text{ before }\Fl\}.
\end{equation}
\begin{remark}
Note that the event $\mathfrak{B}_{\text{low}}(x,y)$ does not require $x$ to be seed-less.
\end{remark}
Using the fact that $\Olow $ and $\Wlow$ are cutpoints, under the event $ \cE_{O} \cap \mathfrak{B}_{\text{low}}(O,\Olow)\cap \mathfrak{B}_{\text{low}}(\Olow,\Wlow) $,
\[
\Upsilon_{\text{low}}^{(1,\lambda)}(O\to B) = \Upsilon_{\text{low}}^{(1,\lambda)}(O\to \Olow) + \Upsilon_{\text{low}}^{(1,\lambda)}(\Olow\to \Wlow)+\Upsilon_{\text{low}}^{(1,\lambda)}(\Wlow \to B).
\]
Thus,
\begin{equation}\label{eq:aux1-OB-low1}
\begin{split}
& \P_{\mu_1}  \left [ \Upsilon_{\text{low}}^{(1,\lambda)}(O\to B)\leq \mathfrak{C} + \frac{H+1}{\cin^{(2)}} + \frac{R}{\lambda \cin^{(1)}} \right ]\\
&
\geq \P_{\mu_1}  \left [ \left \{ \Upsilon_{\text{low}}^{(1,\lambda)}(O\to B)\leq \mathfrak{C} + \frac{H+1}{\cin^{(2)}} + \frac{R}{\lambda \cin^{(1)}} \right \} \cap \cE_{O} \cap \mathfrak{B}_{\text{low}}(O,\Olow)\cap \mathfrak{B}_{\text{low}}(\Olow,\Wlow) \right ] \\
& \geq 
\P_{\mu_1} \left [  \cE_{O} \cap \mathfrak{B}_{\text{low}}(O,\Olow)\cap \left \{ \Upsilon_{\text{low}}^{(1,\lambda)}(O\to O_{\text{low}})\leq \mathfrak{C} \right \} \right  ] \times \\
& \quad  \times \P_{\mu_1} \left [ \mathfrak{B}_{\text{low}}(\Olow,\Wlow)\cap \left \{ \Upsilon_{\text{low}}^{(1,\lambda)}(O_{\text{low}} \to \Wlow)\leq \frac{H+1}{\cin^{(2)}} \right \} \right ] \times \\
& \quad \times  \P_{\mu_1} \left [ \Upsilon_{\text{low}}^{(1,\lambda)}(\Wlow \to B)\leq \frac{R}{\lambda \cin^{(1)}} \right ],
\end{split}
\end{equation}
where we used independence of the passage times in the second bound.

We proceed our analysis by bounding each term in \eqref{eq:aux1-OB-low1}.

\paragraph*{First term.} Since the graph $\{O, \Olow\}$ is actually an edge, then
\[
\begin{split}
& \P_{\mu_1} \bigl [ \cE_{O} \cap \mathfrak{B}_{\text{low}}(O,\Olow) \cap \left \{  \Upsilon_{\text{low}}^{(1,\lambda)}(O\to O_{\text{low}})\leq \mathfrak{C} \right \} \bigr ] \\
& = \P_{\mu_1} \left [ \bigl \{\Upsilon_{\text{low}}^{(1)}(O\to O_{\text{low}})\leq \mathfrak{C}\bigr \} \ \cap \ \cE_{O} \ \cap \ \cE_{O_{\text{low}}} \right ] \stackrel{ \eqref{eq:prob(1-mu)} , \, \eqref{eq:def-cC} }{>} (1-\varepsilon) (1-\mu_1)^2.
\end{split}
\]

\paragraph*{Second term.}

Because of Proposition \ref{prop:Fo-wins}, by taking $o\equiv O_{\text{low}}$, $m=H$, $d=2$ and $W\equiv \Wlow$, the event 
$ \mathfrak{B}_{\text{low}}(\Olow,\Wlow) $ has probability bounded from below by $(1-\f_2(\mu_1)-\bseta_2)(1-\mu_1)$.
For the following part we need a technical construction which will be developed in detail later on, therefore at this point we present a result whose proof is deferred to Section \ref{sect:proof-claim-later}, since it requires a comparison with branching random walks that will be described in Section \ref{sect:Fo-wins}.
At this point we emphasize that the bound obtained in the following result is independent of $\lambda$ (provided that $\lambda$ is chosen to be small).

\begin{claim}\label{claim:to-prove-later}
Suppose that $\mu_1$ and $\mu_2$, as well as $L,H,R,D$, are as in Lemmas \ref{lemma:toProp-mu1} and \ref{lemma:toProp-mu2}.
Let $\varepsilon>0$ be a small enough constant so that $\cin^{(2)}<1$ and $\cin^{(D)}<1$ are small enough.
Then there are values $C,\delta,\beta>0$ independent of $H$ and of $\cin^{(2)}$ such that 
\[
\begin{split}
\P_{\mu_1} & \left [\cE_{\Olow} \cap  \mathfrak{B}_{\text{low}}(\Olow,\Wlow)\cap \left \{  \Upsilon_{\text{low}}^{(1,\lambda)}(O_{\text{low}} \to \Wlow)\leq \frac{H+1}{\cin^{(2)}} \right \} \right ]\\
& \geq (1-\f_2(\mu_1)-\bseta_2)(1-\mu_1)^2- e^{-\beta H}-C e^{-\delta (H+1)/(2\cin^{(2)})}.
\end{split}
\]
Similarly, 
there are constants $\overline{C},\overline{\delta}, \overline{\beta}>0 $ that are independent of $L$ and of $\cin^{(D)}$, so that 
\[
\begin{split}
& \P_{\mu_2}\left [ \cE_{\Oup} \cap  \mathfrak{B}_{\text{up}}(\Oup, \Wup)  \cap \left \{ \Upsilon_\text{up}^{(1,\lambda)}(O_{\text{up}}\to \Wup)\leq \frac{L+1}{\cin^{(D)}} \right \} \right ] \\
& \geq (1-\bseta_D-\f_D(\mu_2))(1-\mu_2)^2 - e^{-\overline{\beta} H}-\overline{C} e^{-\overline{\delta} (L+1)/(2\cin^{(D)})}.
\end{split}
\]
\end{claim}
Hence in this context we choose $H$ large enough, so that $e^{-\beta H}+C e^{-\delta (H+1)/(2\cin^{(2)})}<\varepsilon$.

\paragraph*{Third term.}
By the \emph{geodesic path} connecting $\Wlow$ to $B$ we mean the path of length $R$ as described in Section \ref{sect:construction-tile}.
Define the event
\begin{equation}\label{eq:seeds-AB-1}
\cE_{\text{seed}}:=\{\text{there is a seed on a vertex along the (geodesic) path going from }\Wlow\text{ to }B\}.
\end{equation}
Since 
\begin{equation}\label{eq:prob-seeds-AB-1}
\P_{\mu_1} \left ( \cE_{\text{seed}}\right )=1-(1-\mu_1)^R,
\end{equation}
then, by choosing $R$ large enough uniformly over $\mu_1\in \left (\mu_{1,\text{low}}^{(\text{crit})}, \, \mu_{1,\text{high}}^{(\text{crit})} \right )$, we can make the above probability arbitrarily close to $1$.
We emphasize that if the first seed appears at the beginning of the path, then the passage time of $\Fl$ along the path is higher than if it originated towards the end of the path.
This can be obtained by a simple coupling argument.
In fact, since $\lambda$ is small, then when adding a seed the passage time for each edge $e$ of the path $(t_e^{(\gamma)}) $ is not decreasing.
More precisely, if $e$ is preceding the seed, then $t_e^{(\gamma)}=t_e^{(1)}$, whereas if $e$ is after the seed, then $t_e^{(\gamma)}=t_e^{(\lambda)}>t_e^{(1)}$.
From this observation together with Lemma \ref{lemma:typical-FPP}, applied to $\gamma=\lambda<1$, it follows that for any $\cz ^{(1)}\in (0,1)$ fixed, there is $\cin^{(1)}<1$ such that
\begin{equation}\label{eq:time-AB-1}
\begin{split}
\P_{\mu_1} \left [  \Upsilon_{\text{low}}^{(1,\lambda)}(\Wlow \to B)\leq \frac{R}{\lambda \cin^{(1)}} \right ] & \geq \P_{\mu_1} \left [ \left \{ \Upsilon_{\text{low}}^{(\lambda)}(\Wlow \to B)\leq \frac{R}{\lambda \cin^{(1)}} \right \} \ \cap \ \cE_{\text{seed}} \right ]\\
& \stackrel{\eqref{eq:prob-seeds-AB-1} }{\geq } \left ( 1-e^{-\cz ^{(1)}R }\right ) \left ( 1-(1-\mu_1)^R\right ).
\end{split}
\end{equation}

We can now conclude the proof of the first part of the lemma.

\begin{proof}[Lemma \ref{lemma:toProp-mu1}(i)]
By plugging all these lower bounds into \eqref{eq:aux1-OB-low1} we obtain
\[
\begin{split}
& \P_{\mu_1}  \left [ \Upsilon_{\text{low}}^{(1)}(O\to B)\leq \mathfrak{C} + \frac{H+1}{\cin^{(2)}} + \frac{R}{\lambda \cin^{(1)}}\right ] \\
& \geq 
(1-\varepsilon) (1-\mu_1)^2 \left [ (1-\f_2(\mu_1)-\bseta_2)(1-\mu_1)^2-\varepsilon \right ] 
\left [ 1-e^{-\cz ^{(1)}R }\right ] \left [ 1-(1-\mu_1)^R\right ].
\end{split}
\]
The statement follows by taking 
$R$ large enough, uniformly over $\mu_1\in \left (\mu_{1,\text{low}}^{(\text{crit})}, \, \mu_{1,\text{high}}^{(\text{crit})} \right )$, so that $e^{-\cz ^{(1)}R} \leq \varepsilon$ and $(1-\mu_1)^R < \varepsilon$.
\end{proof}

\subsubsection{Proof of Lemma \ref{lemma:toProp-mu1}(ii), ``upper part'' of $G$.}\label{sect:mu_1_up}
Here we investigate what is likely to happen on the upper part of $G$.
%
%
%
%
%
Here $D>5$ denotes a large integer that will be specified later.
Note that the tree structure implies that the geodesics have length exactly $L+1$.
Also, note the following simple coupling argument: 
for any $\lambda<1$, if $X$ is an exponential random variable of parameter $1$, then $X/\lambda$ is an exponential random variable of parameter $\lambda$.
Therefore, the sum of independent exponential random variables of rate $\lambda$ stochastically dominates the sum of the independent exponential random variables of rate $1$.
Then, by Lemma \ref{lemma:typical-FPP}, 
for any fixed $\co ^{(D)}>0$ there is a $\cout^{(D)}>1$ so that
\begin{equation}\label{eq:time-o1W-1}
\begin{split}
& \P_{\mu_1} \left [ \Upsilon_{\text{up}}^{(1,\lambda)}(O\to B) \geq \frac{L+1}{\cout^{(D)}} \right ] 
\geq \P_{\mu_1} \left [ \Upsilon_{\text{up}}^{(1,\lambda)}(O_{\text{up}}\to \Wup) \geq \frac{L+1}{\cout^{(D)}} \right ]\\
& 
%
\stackrel{\text{coupling}}{\geq} \P_{\mu_1} \left [ \Upsilon_{\text{up}}^{(1)}(O_{\text{up}}\to \Wup) \geq \frac{L+1}{\cout^{(D)}} \right ] 
\geq 1-e^{-\co ^{(D)}(L+1)}.
\end{split}
\end{equation}
Now it suffices to take $L$ so large that $e^{-\co ^{(D)}(L+1)}\leq \varepsilon$ and the proof of the lemma is concluded.

\subsection{Proof of Lemma \ref{lemma:toProp-mu2}}\label{sect:calculations-mu2}
As required in Proposition \ref{prop:Fo-wins} (and in the hypotheses of the Lemma), we will assume that $D$ is large enough and $\mu_2 \in \left (\mu_{2,\text{low}}^{(\text{crit})}, \, \mu_{2,\text{high}}^{(\text{crit})} \right )$ satisfies \eqref{eq:techCondMu}.
Since $D(1-\mu_2)>1 $, then
\begin{equation}\label{eq:P(T2)=0}
\P_{\mu_2} \left (\cE_{\text{perc}}(\T_2) \right )=0 \quad \text{ and } \quad \P_{\mu_2} \left (\cE_{\text{perc}}(\T_{D}) \right )>0.
\end{equation}
This is clearly not the case for $\mu=\mu_1< \frac{1}{2}$, and this difference will turn out to be crucial.

\subsubsection{Proof of Lemma \ref{lemma:toProp-mu2}(i), ``lower part'' of $G$}

Consider the connected component containing $\Olow $ consisting of vertices $\oT_2^H$ which are not hosting seeds (placed with density $\mu_2$).
This set has the same distribution as the first $H$ generations of a \emph{subcritical} Galton-Watson tree (GWT) whose reproduction rule is a Bin$(2, (1-\mu_2))$.

By letting $\mathbf{P}_{\mu_2} $ denote the probability measure defined on the space of Galton-Watson trees with the above-mentioned distribution, then it is a standard fact that
\[
\mathbf{P}_{\mu_2}\left ( \text{subcritical GWT survives for more than }  H/2\, \text{ generations}\right )\leq \bigl ( 2(1-\mu_2)\bigr )^{H/2}.
\]
Thus, the event
\begin{equation}\label{eq:def-subcr}
\cE_{\text{subcr}}(H):=\left \{\text{cluster of non-seeds on $\oT_2^H $ is a tree of height at most } H/2\right \}
\end{equation}
has probability bounded from below by
\begin{equation}\label{eq:prob-subcr}
\mathbf{P}_{\mu_2}\left (\cE_{\text{subcr}}(H)\right )\geq 1-\bigl ( 2(1-\mu_2)\bigr )^{H/2},
\end{equation}
which can be made as close to $1$ as we need, uniformly over all $\mu_2 \in \left ( \mu_{2,\text{low}}^{(\text{crit})}, \mu_{2,\text{high}}^{(\text{crit})}\right )$, by taking a sufficiently large value of $H$. 

Recall the definition $\mathfrak{B}_{\text{low}}(x,y)$ from \eqref{eq:eventB}.
We will show that it is very likely that, when restricting the analysis to the lower part, $\Fl$ takes a very long time before it can potentially infect $\Wlow$.
In order to proceed, we need an auxiliary result.
\begin{lemma}\label{lemma:time-o2A-2}
For the quantities defined above we have
\[
\begin{split}
\P_{\mu_2} & \left ( \mathfrak{B}_{\text{low}}^c (\Olow, \Wlow) \cap  \left \{  \Upsilon_\text{low}^{(1,\lambda)}(O_{\text{low}}\to \Wlow)\geq \frac{H/2}{\cout^{(2)}}+\frac{H/2+1}{\lambda \cout^{(2)}} \right \} \right ) \\
& \geq 
1- 2^{H/2+1} e^{-\co ^{(2)} H/2} -\bigl ( 2(1-\mu_2)\bigr )^{H/2}.
\end{split}
\]
\end{lemma}
\begin{proof}
For every $k\in \N$ define $\mathcal{L}_{2,k}$ to be the set of vertices at generation $k$ of $\oT_2^H$, that is,
\[
\mathcal{L}_{2,k} := \left \{ x\in V(\oT_2^H) \ : \ d_{\oT_2^H}\left ( O_{\text{low}}, x\right )=k
\right \},
\]
where $d_{\oT_2^H}\left ( O_{\text{low}}, x\right )$ is the distance induced by the graph metric on $\oT_2^H$.
We observe that since $\lambda<1$, via a simple coupling argument (similar to the one outlined in Section \ref{sect:mu_1_up}), we see that the probability of the event
\[
\mathfrak{B}_{\text{low}}^c (\Olow, \Wlow) \cap \left \{ \Upsilon_\text{low}^{(1,\lambda)}(O_{\text{low}}\to \Wlow)\geq \frac{H/2}{\cout^{(2)}}+\frac{H/2+1}{\lambda \cout^{(2)}} \right \}
\]
is bounded from below by the probability of the event
\begin{equation}\label{eq:aux-3bits}
\left \{ \Upsilon_\text{low}^{(1)}(O_{\text{low}}\to \mathcal{L}_{2,\frac{H}{2}} ) \geq \frac{ H/2}{\cout^{(2)}} \right \} \, \cap \,  \cE_{\text{subcr}}(H) \, \cap \, \left \{ \Upsilon_\text{low}^{(\lambda)}(\mathcal{L}_{2,\frac{H}{2}} \to \Wlow)\geq \frac{H/2+1}{\lambda \cout^{(2)}}\right \}.
\end{equation}
To see this fact more clearly, observe that the event $\cE_{\text{subcr}}(H)$ only depends on the vertices up to generation $H/2$, and that the function $\Upsilon_\text{low}^{(1)}$ does not take into account the presence of seeds.
Lemma \ref{lemma:typical-FPP} 
yields
\begin{equation}\label{eq:time-Olow-x}
\begin{split}
\P_{\mu_2}\left [ \left \{ \Upsilon_\text{low}^{(1)}(O_{\text{low}}\to \mathcal{L}_{2,\frac{H}{2}} ) \geq \frac{ H/2}{\cout^{(2)}} \right \} \cap  \cE_{\text{subcr}}(H)\right ] 
& \geq  \left [1-e^{-\co ^{(2)} H/2}\right ] \left [ 1-\bigl [ 2(1-\mu_2)\bigr ]^{H/2}\right ] .
\end{split}
\end{equation}
%
%
Furthermore, by Lemma \ref{lemma:typical-FPP} with $\gamma=\lambda$, 
\begin{equation}\label{eq:time-x-Wlow}
\P_{\mu_2}\left [ \Upsilon_\text{low}^{(\lambda)}(\mathcal{L}_{2,\frac{H}{2}}\to \Wlow)\geq \frac{H/2+1}{\lambda \cout^{(2)}} \right ] 
\geq 1-2^{H/2}e^{-\co ^{(2)}\bigl (H/2+1\bigr ) } .
\end{equation}
Thus, by \eqref{eq:time-Olow-x} and \eqref{eq:time-x-Wlow} and using independence of the passage times before and after level $\mathcal{L}_{2,\frac{H}{2}}$,
\[
\begin{split}
& \P_{\mu_2} \left [ \left \{ \Upsilon_\text{low}^{(1)}(O_{\text{low}}\to \mathcal{L}_{2,\frac{H}{2}} ) \geq \frac{ H/2}{\cout^{(2)}} \right \} \, \cap \,  \cE_{\text{subcr}}(H) \, \cap \, \left \{ \Upsilon_\text{low}^{(\lambda)}(\mathcal{L}_{2,\frac{H}{2}} \to \Wlow)\geq \frac{\frac{H}{2}+1}{\lambda \cout^{(2)}}\right \}\right ] \\
& = \P_{\mu_2} \left [ \left \{ \Upsilon_\text{low}^{(1)}(O_{\text{low}}\to \mathcal{L}_{2,\frac{H}{2}} ) \geq \frac{ H/2}{\cout^{(2)}} \right \} \cap   \cE_{\text{subcr}}(H) \right ]  \P_{\mu_2} \left [ \Upsilon_\text{low}^{(\lambda)}(\mathcal{L}_{2,\frac{H}{2}} \to \Wlow)\geq \frac{\frac{H}{2}+1}{\lambda \cout^{(2)}}\right ]
\\
& \geq \left ( 1-e^{-\co ^{(2)} H/2} -\bigl ( 2(1-\mu_2)\bigr )^{H/2}\right ) \left ( 1-2^{H/2}e^{-\co ^{(2)}\bigl (H/2+1\bigr ) }\right )\\
& \geq 1- 2^{H/2+1} e^{-\co ^{(2)} H/2} -\bigl ( 2(1-\mu_2)\bigr )^{H/2},
\end{split}
\] 
as claimed.
\end{proof}
Furthermore, Lemma \ref{lemma:typical-FPP} applied to the path from $\Wlow$ to $B$ when $\gamma=\lambda<1$ gives
\begin{equation}\label{eq:time-AB-2}
%
\P_{\mu_2}\left [\Upsilon_\text{low}^{(\lambda)}(\Wlow \to B) \geq \frac{R}{\lambda \cout^{(1)}} \right ]\geq 1-e^{-\co ^{(1)}R}.
\end{equation}
Now we are ready to proceed with the proof of the first part of  Lemma \ref{lemma:toProp-mu2}.
\begin{proof}[Proof of Lemma \ref{lemma:toProp-mu2}(i)]
%
Since $\Wlow$ is a cutpoint we can bound the sought probability as 
\[
\begin{split}
& \P_{\mu_2} \left [  \Upsilon_\text{low}^{(1,\lambda)}(O\to B)\geq  \frac{H/2}{\cout^{(2)}}+\frac{H/2+1}{\lambda \cout^{(2)}}+ \frac{R}{\lambda \cout^{(1)}}  \right ]\\
& \geq 
\P_{\mu_2}\Biggl [ \mathfrak{B}_{\text{low}}^c (\Olow, \Wlow) \cap \left \{\Upsilon_\text{low}^{(1,\lambda)}(O_{\text{low}}\to \Wlow)\geq \frac{H/2}{\cout^{(2)}}+\frac{H/2+1}{\lambda \cout^{(2)}} \right \}  \cap \\
& \qquad \cap \left \{  \Upsilon_\text{low}^{(\lambda)}(\Wlow \to B)\geq \frac{R}{\lambda \cout^{(1)}}\right \}  \Biggr ]\\
& =
\P_{\mu_2}\Biggl [ \mathfrak{B}_{\text{low}}^c (\Olow, \Wlow) \cap  \left \{ \Upsilon_\text{low}^{(1,\lambda)}(O_{\text{low}}\to \Wlow)\geq \frac{H/2}{\cout^{(2)}}+\frac{H/2+1}{\lambda \cout^{(2)}} \right \} \Biggr ] \times \\
& \qquad \times  \P_{\mu_2}\left [ \Upsilon_\text{low}^{(\lambda)}(\Wlow \to B)\geq \frac{R}{\lambda \cout^{(1)}} \right ].
\end{split}
\]
By Lemma  \ref{lemma:time-o2A-2} the above is bounded from below by
\[
\begin{split}
& \left ( 1- 2^{H/2+1} e^{-\co ^{(2)} H/2} -\bigl ( 2(1-\mu_2)\bigr )^{H/2} \right ) \P_{\mu_2}\left [ \Upsilon_\text{low}^{(\lambda)}(\Wlow \to B)\geq \frac{R}{\lambda \cout^{(1)}} \right ]\\
& \stackrel{\eqref{eq:time-AB-2} }{\geq } \left (1-2^{H/2+1} e^{-\co ^{(2)} H/2}- \bigl ( 2(1-\mu_2)\bigr )^{H/2}\right ) \left ( 1-e^{-\co ^{(1)}R} \right ).
\end{split}
\]
At this point we can choose $R$ and $H$ so large that, for all values of $\mu_2 \in \left (\mu_{2,\text{low}}^{(\text{crit})}, \, \mu_{2,\text{high}}^{(\text{crit})} \right ) $ one has $e^{-\co ^{(1)}R}< \varepsilon$, $\bigl( 2(1-\mu_2)\bigr )^{H/2}< \varepsilon$ and $2^{H/2+1}e^{-\co ^{(2)} H/2}<\varepsilon$, hence the statement follows.
\end{proof}

\subsubsection{Proof of Lemma \ref{lemma:toProp-mu2}(ii), ``upper part'' of $G$}

Similarly to the previous section, consider the variable $\Upsilon_{\text{up}}^{(1,\lambda)}(O\to B)$, fix any pair of vertices $x,y$ in the upper part of $G$ and an initial density of seeds $\mu$.
For this part of the proof consider only the upper  part of $G$, more precisely define the passage times $\{ \overline{\overline{t}}_e^{(\gamma)}\}$ as follows.
Set $ \overline{\overline{t}}_e^{(\gamma)}:=t_e^{(\gamma)}$ if the edge $e$ is in the upper part of $G$, and set $\overline{\overline{t}}_e^{(\gamma)}:=\infty$, if $e$ is in the lower part.
With respect to the passage times $\{\overline{\overline{t}}_e^{(\gamma)}\}_{e\in E(G)}$, recalling the definition of $\cE_{y}$ from \eqref{eq:def-Ex-mu}, define the event
\begin{equation}\label{eq:eventBup}
\mathfrak{B}_{\text{up}}(x,y):=\cE_{y}\cap \{ \Fo \text{ started at }x\text{ infects }y\text{ before }\Fl\}.
\end{equation}
By Proposition \ref{prop:Fo-wins}, 
\begin{equation}\label{eq:Bup}
\P_{\mu_2} \left [ \mathfrak{B}_{\text{up}}(\Oup ,\Wup)\right ]\geq 
\bigl (1 - \f_D(\mu_2)-\bseta_D\bigr )(1-\mu_2) .
\end{equation}
Now we want to show that, with high probability, when the initial density of seeds is $\mu_2\in \left (\mu_{2,\text{low}}^{(\text{crit})}, \, \mu_{2,\text{high}}^{(\text{crit})} \right )$ the time needed by $\Fo$ started at $O$ to reach $B$ is smaller than that needed by $\Fl$.

Now recall the definition of $\mathfrak{C}$ from \eqref{eq:def-cC} and the bounds obtained in Lemma \ref{lemma:time-o2A-2} and \eqref{eq:time-AB-2}. 
At this point we are ready to prove the second bound of Lemma \ref{lemma:toProp-mu2}.

\begin{proof}[Proof of Lemma \ref{lemma:toProp-mu2}(ii)]
We can bound the sought probability by
\begin{equation}\label{eq:aux1-OB-up2}
\begin{split}
& \P_{\mu_2} \left [ 
\left \{ \Upsilon_{\text{up}}^{(1,\lambda)}(O\to B)\leq 2\mathfrak{C}+\frac{L+1}{\cin^{(D)}}\right \} \cap \left \{ \Wup \text{ is infected by }\Fo \right \} \right ] \\
& \geq  \P_{\mu_2} \left [\left \{ \Upsilon_\text{up}^{(1,\lambda)}(O\to B) \leq 2\mathfrak{C}+\frac{L+1}{\cin^{(D)}} \right \} \cap \cE_{O}\cap \cE_{\Oup} \cap  \mathfrak{B}_{\text{up}}(\Oup, \Wup) \cap \cE_{\Wup} \right ].
\end{split}
\end{equation}
Using the fact that $\Oup$ and $\Wup$ are cutpoints, under the event $\cE_{O}\cap \cE_{\Oup} \cap  \mathfrak{B}_{\text{up}}(\Oup, \Wup) $ we have
\[
\Upsilon_\text{up}^{(1,\lambda)}(O\to B) = \Upsilon_\text{up}^{(1,\lambda)}(O \to O_{\text{up}}) + \Upsilon_\text{up}^{(1,\lambda)}(O_{\text{up}}\to \Wup) + \Upsilon_\text{up}^{(1,\lambda)}(\Wup \to B) ,
\]
hence the last line in \eqref{eq:aux1-OB-up2} is bounded from below by
\begin{equation}\label{eq:aux2-up}
\begin{split}
& \P_{\mu_2}\Bigl [ 
\left \{ \Upsilon_\text{up}^{(1,\lambda)}(O \to O_{\text{up}})\leq \mathfrak{C}\right \} \cap 
\left \{ \Upsilon_\text{up}^{(1,\lambda)}(O_{\text{up}}\to \Wup)\leq \frac{L+1}{\cin^{(D)}} \right \} \cap \\
& \qquad \cap \left \{ \Upsilon_\text{up}^{(1,\lambda)}(\Wup \to B)\leq \mathfrak{C} \right \} \cap \cE_{O} \cap \cE_{\Oup} \cap  \mathfrak{B}_{\text{up}}(\Oup, \Wup) \cap \cE_{\Wup} \Bigr ] .
\end{split}
\end{equation}
Under the event $\cE_{O}\cap \cE_{\Oup} \cap  \mathfrak{B}_{\text{up}}(\Oup, \Wup) \cap \cE_{\Wup}$ we have that 
\begin{itemize}
\item $\Upsilon_\text{up}^{(1,\lambda)}(O \to O_{\text{up}})$ is equal to $\Upsilon_\text{up}^{(1)}(O \to O_{\text{up}})$;
\item $\Upsilon_\text{up}^{(1,\lambda)}(\Wup \to B)$ is equal to $\Upsilon_\text{up}^{(1)}(\Wup \to B)$.
\end{itemize}
Furthermore, the distributions of such random variables only depend on the passage times of $\Fo$, which are independent of the position of the seeds.
Therefore, the quantity in \eqref{eq:aux2-up} is bounded from below by
\[
\begin{split}
& \P_{\mu_2} \left [ \left \{ \Upsilon_\text{up}^{(1)}(O \to O_{\text{up}})\leq \mathfrak{C} \right \} \cap \cE_{O} \right ] 
\P_{\mu_2} \left [ \Upsilon_\text{up}^{(1)}(\Wup \to B)\leq \mathfrak{C} \right ] \times \\
& \quad \times 
\P_{\mu_2}\left [ \left \{ \Upsilon_\text{up}^{(1,\lambda)}(O_{\text{up}}\to \Wup)\leq \frac{L+1}{\cin^{(D)}} \right \} \cap \cE_{\Oup} \cap  \mathfrak{B}_{\text{up}}(\Oup, \Wup) \cap \cE_{\Wup} \right ] 
\\
& \stackrel{\eqref{eq:def-cC} }{\geq } 
(1-\varepsilon)^2 (1-\mu_2)^2
\P_{\mu_2}\left [ \left \{ \Upsilon_\text{up}^{(1,\lambda)}(O_{\text{up}}\to \Wup)\leq \frac{L+1}{\cin^{(D)}} \right \} \cap  \cE_{\Oup} \cap  \mathfrak{B}_{\text{up}}(\Oup, \Wup)  \right ]
\end{split}
\]
For the last term of the last line, we proceed analogously to the proof of Lemma \ref{lemma:toProp-mu1}.
Consider now a sufficiently small value of $\cin^{(D)}$.
By using similar arguments to the ones mentioned in Section \ref{sect:mu_1_low}, in particular by Claim \ref{claim:to-prove-later}, one finds that there are constants $\overline{C},\overline{\delta}, \overline{\beta}>0 $ that are independent of $L$ and of $\cin^{(D)}$, so that 
\[
\begin{split}
& \P_{\mu_2}\left [ \left \{ \Upsilon_\text{up}^{(1,\lambda)}(O_{\text{up}}\to \Wup)\leq \frac{L+1}{\cin^{(D)}} \right \} \cap \cE_{\Oup} \cap  \mathfrak{B}_{\text{up}}(\Oup, \Wup)  \right ] \\
& \geq (1-\bseta_D-\f_D(\mu_2))(1-\mu_2)^2 - e^{-\overline{\beta} H}-\overline{C} e^{-\overline{\delta} (L+1)/(2\cin^{(D)})}.
\end{split}
\]
%
%
Now it suffices to pick $L, H$ large enough so that $ e^{-\overline{\beta} H}+\overline{C} e^{-\overline{\delta} (L+1)/(2\cin^{(D)})} <\varepsilon $.
%
By putting all these bounds together into \eqref{eq:aux1-OB-up2} and using \eqref{eq:Bup}, we obtain the claim.
\end{proof}

\section{Proof of Proposition \ref{prop:aux}}\label{sect:proof-prop-aux}
We start with a technical lemma which we use to appropriately set $H,R$ and $L$.

\begin{lemma}\label{lemma:to-change}
As before, we consider $\lambda$ to be an arbitrarily small constant fixed beforehand; let $D$ be large.
For $\mathfrak{C}$ fixed as in \eqref{eq:def-cC}, for any large integer $R$, there is a choice of values  $H,L$, so that
\begin{equation}\label{eq:1-red-infects-B}
\frac{H+1}{\cin^{(2)}}+\frac{R}{\lambda \cin^{(1)}}+ \mathfrak{C} < \frac{L+1}{\cout^{(D)}},
\end{equation}
and 
\begin{equation}\label{eq:2-white-infects-B}
\frac{H/2}{\cout^{(2)}}+\frac{H/2+1}{\lambda \cout^{(2)}}+ \frac{R}{\lambda \cout^{(1)}} > 2\mathfrak{C}+\frac{L+1}{\cin^{(D)}}
\end{equation}
hold simultaneously.
\end{lemma}
\begin{proof}
   By solving \eqref{eq:1-red-infects-B} and \eqref{eq:2-white-infects-B} for $(L+1)$ and comparing the results, we see that it suffices to prove the following (stronger) condition, 
   \begin{equation}\label{eq:1+2}
   \cout^{(D)}\left ( \frac{H+1}{\cin^{(2)}}+\frac{R}{\lambda \cin^{(1)}}+ 2\mathfrak{C} \right ) +1 \ < \ \cin^{(D)}\left (\frac{H/2}{\cout^{(2)}}+\frac{H/2}{\lambda \cout^{(2)}}+ \frac{R}{\lambda \cout^{(1)}} -2\mathfrak{C} \right ).
   \end{equation}
   In other words, \eqref{eq:1+2} ensures that we can use \eqref{eq:1-red-infects-B} and \eqref{eq:2-white-infects-B} at the same time.
   The value ``$+1$'' on the left-hand side guarantees we can pick a large integer $L$ depending on $\lambda, H,R, \mathfrak{C}, \cin^{(1)}, \cin^{(2)}, \cin^{(D)}, \cout^{(1)}, \cout^{(2)}, \cout^{(D)}$ that satisfies both \eqref{eq:1-red-infects-B} and \eqref{eq:2-white-infects-B}.
   Relation \eqref{eq:1+2} is equivalent to
   \begin{equation}\label{eq:second-last18}
   \frac{\cout^{(D)}}{\cin^{(2)}} +\frac{R}{\lambda} \left ( \frac{\cout^{(D)}}{\cin^{(1)}} -\frac{\cin^{(D)}}{\cout^{(1)}} \right ) + 2\mathfrak{C} \left ( \cout^{(D)}+\cin^{(D)}\right )+1 
   < H\left (\frac{\cin^{(D)}}{2\cout^{(2)}}\Bigl (1+\frac{1}{\lambda}\Bigr ) - \frac{\cout^{(D)}}{\cin^{(2)}} \right ).
   \end{equation}
   Recall from Lemma~\ref{lemma:typical-FPP} that $\cout^{(d)}>1>\cin^{(d)}$ for all $d$.
Suppose now that $R$ is a large value chosen in a convenient way.
   Then we need to check that the coefficient of $H$ is positive, that is,
   \[
   \frac{\cin^{(D)}}{2\cout^{(2)}}\Bigl (1+\frac{1}{\lambda}\Bigr ) - \frac{\cout^{(D)}}{\cin^{(2)}} >0,
   \]
   which is guaranteed for all values $0<\lambda<\lambda_0$ where we have set 
   \[
   \lambda_0:= 
   \frac{\cin^{(D)}\cin^{(2)}}{2\cout^{(2)}\cout^{(D)}-\cin^{(D)}\cin^{(2)}}>0.
   \]
   The proof is now completed since, after fixing $R$, it suffices to choose $H$ large enough so that \eqref{eq:second-last18} is satisfied. 
\end{proof}


Now we proceed by showing that for all $\mu_1\in \left (\mu_{1,\text{low}}^{(\text{crit})}, \, \mu_{1,\text{high}}^{(\text{crit})} \right )$, it is likely that $B$ is infected by the process spreading in the lower part of the tile, whereas for all $\mu_2 \in \left (\mu_{2,\text{low}}^{(\text{crit})}, \, \mu_{2,\text{high}}^{(\text{crit})} \right )$ it is likely to be the opposite.
We recall that, as in Lemma \ref{lemma:toProp-mu2}, the event $\left \{ \Wup \text{ is infected by }\Fo \right \}$ is restricted to the upper part of the tile, so we do not consider the possibility that 
$\Wup$ is infected by $B$.

\begin{lemma}\label{lemma:to-changeNext}
Let $G$ be constructed as in Section \ref{sect:construction-tile}. 
Fix a large integer $D$ and a small enough $\lambda>0$. 
Then there are values $0<\mu_{1,\text{low}}^{(\text{crit})}<\mu_{1,\text{high}}^{(\text{crit})} <1/2 <\mu_{2,\text{low}}^{(\text{crit})}<\mu_{2,\text{high}}^{(\text{crit})} <1$ satisfying 
\[
D\left (1-\mu_{2,\text{high}}^{(\text{crit})} \right )>1 \quad \text{ and }\quad  D^2 \left (1-\mu_{2,\text{low}}^{(\text{crit})}\right )^2 \left (\mu_{2,\text{high}}^{(\text{crit})} \right )^{D-1}<1,
\]
as well as values $H$, $L$ and $R$ satisfying Lemma \ref{lemma:to-change} together with \eqref{eq:time-AB-1} and \eqref{eq:time-AB-2}, so that
\[
\begin{split}
&  
\P_{\mu_1} \left (  \Upsilon_{\text{low}}^{(1,\lambda)}(O\to B)<\Upsilon_{\text{up}}^{(1, \lambda)}(O\to B)\right )\geq 1-\bseta_2-\f_2(\mu_1)-5\varepsilon-4\mu_1,
\\
& 
\P_{\mu_2} \left ( \left \{ \Upsilon_{\text{low}}^{(1, \lambda)}(O\to B)>\Upsilon_{\text{up}}^{(1, \lambda)}(O\to B) \right \} \cap \left \{ \Wup \text{ is infected by }\Fo \right \} \right )\\
& \qquad \geq 
(1-5\varepsilon) (1-\mu_2)^2\left [ \bigl (1-\bseta_D-\f_D(\mu_2)\bigr )(1-\mu_2)^2-\varepsilon \right ] ,
\end{split}
\]
where $\bseta_2$ and $\bseta_D$ are obtained as in Proposition \ref{prop:Fo-wins}.
\end{lemma}
\begin{proof}
In order to bound the first formula we use Lemma \ref{lemma:toProp-mu1} together with relation \eqref{eq:1-red-infects-B} and the fact that FPPHE spreads in the upper and lower part of the tile independently (until $B$ gets infected).
We find that
\[
\begin{split}
& 
\P_{\mu_1} \left ( \Upsilon_{\text{low}}^{(1,\lambda)}(O\to B)<\Upsilon_{\text{up}}^{(1, \lambda)}(O\to B)\right )\\
& 
\geq  
\P_{\mu_1} \left [  \left \{ \Upsilon_{\text{low}}^{(1, \lambda)}(O\to B)\leq \mathfrak{C} + \frac{H+1}{\cin^{(2)}} + \frac{R}{\lambda \cin^{(1)}}\right \}  \cap \left \{ \Upsilon_{\text{up}}^{(1, \lambda)}(O\to B)\geq \frac{L+1}{\cout^{(D)}}\right \} \right ]   \\
& \geq  \bigl (1-\varepsilon\bigr )^3 \bigl (1-\mu_1\bigr )^2
\bigl[ (1 - \f_2(\mu_1)-\bseta_2 )(1-\mu_1)^2-\varepsilon \bigr ]
( 1-\varepsilon ),
\end{split}
\]
implying the first part of the claim.
Analogously, for the second formula of the statement we use Lemma \ref{lemma:toProp-mu2}, relation \eqref{eq:2-white-infects-B} and again independence of FPPHE in the upper and lower part, finding
\[
\begin{split}
& 
\P_{\mu_2} \left ( \left \{ \Upsilon_{\text{low}}^{(1, \lambda)}(O\to B)>\Upsilon_{\text{up}}^{(1, \lambda)}(O\to B) \right \} \cap \left \{ \Wup \text{ is infected by }\Fo \right \} \right )\\
& 
\geq  
\P_{\mu_2} \Bigl [ \left \{ \Upsilon_{\text{low}}^{(1, \lambda)}(O\to B)\geq \frac{H/2}{\cout^{(2)}}+\frac{H/2+1}{\lambda \cout^{(2)}}+ \frac{R}{\lambda \cout^{(1)}}\right \}  \cap \left \{ \Upsilon_{\text{up}}^{(1, \lambda)}(O\to B)\leq 2\mathfrak{C}+\frac{L+1}{\cin^{(D)}}\right \} \cap \\
& \qquad \cap \left \{ \Wup \text{ is infected by }\Fo \right \} \Bigr ]\\
& \geq  (1-3\varepsilon)(1-\varepsilon)^2 (1-\mu_2)^2\left [ \bigl (1-\bseta_D-\f_D(\mu_2)\bigr )(1-\mu_2)^2-\varepsilon \right ] ,
\end{split}
\]
which implies the claim.
\end{proof}

Now recall that $\lambda$ is small enough and note that, given $\mu_{2,\text{high}}^{(\text{crit})} $,  the function $\f_D( \mu_{2,\text{high}}^{(\text{crit})} ) $ is decreasing in $D$.
Moreover, 
by Proposition \ref{prop:Fo-wins}(iii),
the quantity $\bigl (1-\bseta_D-\f_D(\mu_{2,\text{high}}^{(\text{crit})})\bigr )(1-\mu_{2,\text{high}}^{(\text{crit})})^2$ can be made positive by choosing $D$ large enough.
Hence, let $D$ be so large that, for all $\mu_2\in \left (\mu_{2,\text{low}}^{(\text{crit})}, \, \mu_{2,\text{high}}^{(\text{crit})} \right ) $ one has $\bigl (1-\bseta_D-\f_D(\mu_2)\bigr )(1-\mu_2)^2>0$.
Hence we assume that $\varepsilon$ is so that, uniformly in $\mu_2 \in \left (\mu_{2,\text{low}}^{(\text{crit})}, \, \mu_{2,\text{high}}^{(\text{crit})} \right )$,
\begin{equation}\label{eq:DaDefinire}
\varepsilon<\min \left \{ \frac{1}{700}, \, (1-\bseta_D-\f_D(\mu_2) )(1-\mu_2)^2 \right \}.
\end{equation}
Since $\bigl (1-\bseta_D-\f_D(\mu_2)\bigr )(1-\mu_2)^2-\varepsilon >0$, define
\begin{equation}\label{eq:defPhi}
\Phi:= \sup_{\mu_2 \in \left (\mu_{2,\text{low}}^{(\text{crit})}, \, \mu_{2,\text{high}}^{(\text{crit})} \right ) }\left \lceil \frac{2D^3}{\frac{99}{100}
\left [ \bigl (1-\bseta_D-\f_D(\mu_2)\bigr )(1-\mu_2)^2-\varepsilon \right ]} \right \rceil;
\end{equation}
this quantity will be left constant throughout.
Now we are ready to prove Proposition \ref{prop:aux}.
Here we shall use that the critical value $\mu_{1,\text{high}}^{(\text{crit})}$ is small enough; more precisely, with the notation above, let $\overline{\mu_1}:=\sup\{ \mu \in (0,1) \ : \ 7\varepsilon+5\mu+\bseta_2+\f_2(\mu) < 1/\Phi\} $.
Then we need
\begin{equation}\label{eq:def-mu1-bar}
0<\mu_{1,\text{high}}^{(\text{crit})}\leq \overline{\mu_1}.
\end{equation}

\begin{proof}[Proof of Proposition \ref{prop:aux}]
Start by considering $\mu_1 \in \left (\mu_{1,\text{low}}^{(\text{crit})}, \, \mu_{1,\text{high}}^{(\text{crit})} \right )$ and recall the definition of $\cE_y$ from \eqref{eq:def-Ex-mu}.
Then, since $R$ was so that $(1-\mu_1)^{R}<\varepsilon$ (cf.\ \eqref{eq:time-AB-1} and \eqref{eq:time-AB-2}). 
\[
\begin{split}
& 
\P_{\mu_1} \left (  B\text{ is infected by }\Fo \right ) \\
& \leq 
\P_{\mu_1} \left (  \left \{ \Upsilon_{\text{low}}^{(1,\lambda)}(O\to B)> \Upsilon_{\text{up}}^{(1,\lambda)}(O\to B) \right \}\cap \cE_B \right ) 
+ \P_{\mu_1} \left ( \text{there is no seed on the path to }B\right )\\
& \stackrel{\text{Lem.\ \ref{lemma:to-changeNext} } }{\leq }\bigl ( 1- \left (1-5\varepsilon-4\mu_1-\bseta_2-\f_2(\mu_1)\right ) \bigr  )(1-\mu_1) + (1-\mu_1)^{R-1} \\
& <  5\varepsilon+5\mu_1+\bseta_2+\f_2(\mu_1) +2\varepsilon.
\end{split}
\]
Using the fact that $\varepsilon$ is small (cf.\ \eqref{eq:DaDefinire}), as well as $\lambda$, 
together with \eqref{eq:def-mu1-bar}, then, for all $\mu_1\in \left (\mu_{1,\text{low}}^{(\text{crit})}, \, \mu_{1,\text{high}}^{(\text{crit})} \right )$ we have  $7\varepsilon+5\mu_1+\bseta_2+\f_2(\mu_1) < 1/\Phi$, concluding the first part.
Analogously,
\[
\begin{split}
& 
\P_{\mu_2} \left ( B\text{ is infected by }\Fo \right ) \\
& \geq 
\P_{\mu_2} \left [ \left \{  \Upsilon_{\text{low}}^{(1,\lambda)}(O\to B)> \Upsilon_{\text{up}}^{(1,\lambda)}(O\to B)\right \} \cap \cE_B \cap \{\Wup \text{ is infected by }\Fo \} \right ]\\
& \stackrel{\text{Lem.\ \ref{lemma:to-changeNext} } }{\geq }
(1-5\varepsilon) \left [ \bigl (1-\bseta_D-\f_D(\mu_2)\bigr )(1-\mu_2)^2-\varepsilon \right ] (1-\mu_2)^3 \\
& \stackrel{\eqref{eq:techCondMu} }{> }
(1-5\varepsilon) \left [ \bigl (1-\bseta_D-\f_D(\mu_2)\bigr )(1-\mu_2)^2-\varepsilon \right ] \frac{1}{D^3}.
\end{split}
\]
By our assumption on $\varepsilon$ we have that $1-5\varepsilon>99/100$.
By \eqref{eq:defPhi},
\[
(1-5\varepsilon) \left [ \bigl (1-\bseta_D-\f_D(\mu_2)\bigr )(1-\mu_2)^2-\varepsilon \right ] \frac{1}{D^3}\geq \frac{2}{\Phi}> \frac{1}{\Phi}, 
\]
as claimed.
\end{proof}

\section{Analysis of FPPHE on a finite tree: proof of Proposition \ref{prop:Fo-wins}}\label{sect:Fo-wins}
In this section we analyze FPPHE on the tree $\oT_d^m$ and prove Proposition \ref{prop:Fo-wins}. Since the proof is articulated into several intermediate results, we provide a sketch.
The idea is the following.
If $\Fl$ occupies $W$, then we can find an index $i\in \{1, \ldots , m\}$ corresponding to a generation, such that there is a \emph{seed} at generation $i$ that is activated by $\Fo$ and which belongs to the geodesic that reaches $W$.
Note that $i$ is the generation of the first seed in the geodesic between the root and $W$.

Once we fix an arbitrarily small constant $\uep{e0}>0$, ``Step 1'' and ``Step 2'' below deal with the cases $i\leq \lceil (1-\uep{e0})m \rceil $  and $i> \lceil (1-\uep{e0})m \rceil $ respectively.
From now on, in order to avoid cumbersome notation, instead of writing $\lceil (1-\uep{e0})m \rceil$ we only write $ (1-\uep{e0})m$.

\nep{e0}
\textbf{Step 1. }
Fix an arbitrarily small constant $\uep{e0}>0$.
We start by showing that the probability that there is a seed at some generation $i\leq  (1-\uep{e0})m $ that originates a $\Fl$ process which will reach level $m$ much quicker than expected, decreases exponentially fast in $m$.
The underlying idea is that 
if $i\leq  (1-\uep{e0})m $, then, in order to reach generation $m$, $\Fl$ has to cross at least $\uep{e0}m$ edges, and since $m$ is large, this is unlikely.
This step then will show that if $\Fl$ does occupy $W$, then it is more likely that it did so from a seed within distance $\uep{e0}m$ of $W$.
%
%

\textbf{Step 2. }
Subsequently, we take care of what happens at levels $i $ for $  (1-\uep{e0})m  \leq i \leq m$.
Since the seed is close to $W$, we cannot argue that the path traversed by $\Fl$ from the seed to $W$ has a large delay with respect to what $\Fo$ would have had. 
Note that with probability at least $\mu$ the fastest path to $W$ will have a seed in its last vertex before reaching $W$.
Therefore, in order to give an upper bound to the probability that $W$ is occupied by $\Fl$, we will consider what happened on the tree before level $(1-\uep{e0})m$.
For simplicity, set $j:=m-i$; here is the idea behind this procedure.
\begin{itemize}
\item[(i)] We show that for all $j\in \{1, \ldots , \lfloor \uep{e0}m \rfloor \} $, it is likely that \emph{all} vertices of type $\Fo$ at level $j$ are infected within some time proportional to $j$. 
\item[(ii)] Fix $j\in \{1, \ldots , \lfloor \uep{e0}m \rfloor \} $.
We consider the (unique, if any) geodesic of vertices infected by $\Fo$ started at each vertex at level $j$ described in (i), and ended at the first infected vertex at level $m-j$.
Note that up to level $m-j-1$ such line consists of vertices infected by $\Fo$, whereas the one at level $m-j$ can be either of type $\Fo$ or a seed.
\item[(iii)] We use a relation between first passage percolation and a continuous-time branching random walk, as well as concentration results for the latter, 
    to show that all vertices at level $m-j$ described in (ii) are infected within a time that is very concentrated around its expected value (which is of order $m-j$).
\item[(iv)] Next, we show that the probability that among the first-infected particles described in (iii) there is a non-seed originating a $\Fo$ process that infects level $m$ before any other process does is high.
\end{itemize}

Before proceeding we make a brief reflection about the steps above. 
Part (i) is used to get several independent vertices at level $j$ to be occupied by $\Fo$, for which the arrival time at level $j$ is within 
distance of order $j$ of one another. 
Parts (ii) and (iii) give that each such vertex will give rise to a geodesic path to level $m-j$, whose passage time is very concentrated around a given value of order $m-j$. 
Call the collection of such geodesic paths $\mathcal{P}$. 
Since $m-j$ is much bigger than $j$, 
the paths in $\mathcal{P}$ all arrive at level $m-j$ within a small time interval. 
Part (iv) gives that, since $\mathcal{P}$ reaches several vertices at level $m-j$,
it is unlikely that too many of them are occupied by seeds. Similarly, since $\mathcal{P}$ is large, it is unlikely that a seed from level $m-j$ that is activated by a non-geodesic path (that is, by a path not in $\mathcal{P}$) 
arrives at level $m$ faster than the paths in $\mathcal{P}$.


Recall that $d\geq 2$ is the degree of the tree, and consider $\p$ 
as in the statement of the proposition.
In this proof we will deal with the first passage times restricted to $\oT_d^{m}$, thus to make this more explicit we shall write $\Upsilon^{(1,\lambda)}_{\oT_d^{m}}$.
%
%
Recall that $m$ is chosen to be large.
For $d\geq 2$ and $0\leq k\leq m$ let 
\[
\mathcal{L}_k= \mathcal{L}_{d,k}:= \{x \in V(\oT_d^{m} ) \ : \ d_{\oT_d^{m}}(o , x)=k\},
\]
that is, the set of all vertices of $\oT_d^{m}$ at distance $k$ from the root.

\paragraph{Auxiliary facts.}
For $d\geq 2$ and $\p$ as in \eqref{eq:techCondMu}, for all levels $i\in \{1, \ldots , m\}$ define the events 
\begin{equation}\label{eq:def-event-Bi}
\begin{split}
B_i :=& \{\exists\text{ seed in } \mathcal{L}_i\text{ which is infected by }\Fo \text{ and}\\
& \text{ gives rise to a }\Fl \text{ process that infects } \mathcal{L}_m\text{ first}\}.
\end{split}
\end{equation}
These are ``bad'' events, thus we need to show that the probability that there is a level $i$ for which $B_i$ occurs is small.
To shorten the notation, sometimes we shall write ``vertex $x$ is of type $1$ (resp.\ $\lambda$)'' as an equivalent formulation of ``vertex $x$ is of type $\Fo$ (resp.\ $\Fl$)''.

A fundamental tool is provided by \cite[Theorem 3]{AddarioBerry-Reed} which we recall below.
Let $\gamma\in \{1, \lambda\}$.
For $d\geq 2$ and each $i,n\in\{1, \ldots, m\}$ with $i<n$, and for every fixed $x\in \mathcal{L}_i$, $x$ of type $\gamma$, set
\begin{equation}\label{eq:defMx}
M_{i,n}^{(\gamma)}(x):=
\inf_{y\in \mathcal{L}_n, ~ y\text{ of type } \gamma, ~ y \text{ descendant of }x} \{\Upsilon_{\oT_d^{m}}^{(\gamma)}(x \to y)\},
\end{equation}
with the standard convention that $\inf \{\emptyset\}:= +\infty$, and thus define
\[
M_{i,n}^{(\gamma)}:= \inf_{x\in \mathcal{L}_i, ~ x \text{ of type }
\gamma }M_{i,n}^{(\gamma)}(x).
\]
In words, $M_{i,n}^{(\gamma)}$ is the time needed by the quickest among the infections started at any $x\in\mathcal{L}_i $ to reach $\mathcal{L}_n$ when the infection spreads at rate $\gamma$.
%
%
Now we consider a Galton-Watson branching process (where individuals have i.i.d.\ offspring) whose branching distribution is given by Bin($d, \, \ns$).
Define $\mathscr{S}$ to be the event
\[
\mathscr{S}:=\{\text{branching process with distribution given by Bin($d, \, \ns$) survives}\}.
\]
By \eqref{eq:defBFf},
\begin{equation}\label{eq:PS}
\P_\p (\mathscr{S})= 1-\f_d(\p).
\end{equation}

\paragraph{Description of a Branching Random Walk in Continuous Time (BRWCT).}
We will couple the process $\Fo $ inside $\oT_d^m$ with initial density of seeds $\mu$ with a BRWCT on $\mathbb{R}$ as follows. 
\begin{itemize}
\item
Each vertex of type $\Fo $ corresponds to an individual in the BRWCT.
The origin $o$ plays the role of the starter of the progeny, and for $j\in \{0,\ldots, m-1\}$ we say that individual $u \in \mathcal{L}_j$ is the \emph{father} of $x\in \mathcal{L}_{j+1}$ if and only if \textbf{(a)} $u$ and $x$ are not seeds, and \textbf{(b)} $\{u,x\}\in E(\oT_d^m)$.
For $x,y\in V(\oT_d^m)$, we say that $y$ is a \emph{descendant} of $x$ if there is a line of fathers from $y$ to $x$.
\item
The distribution governing the underlying random walk has (independent) increments distributed as exponentials of parameter $\gamma=1$.
\item
The offspring distribution is binomial with parameters $d $ and $1-\den$.
\item The above defines the BRWCT up to generation $m$; from generation $m+1$ onwards, we let the BRWCT evolve independently of $\Fo$ on $\oT_d^m$.
\end{itemize}
Whenever we will write ``GWT'' we will mean a Galton-Watson Tree having this distribution.
Also, from the coupling above, we see that the event that the BRWCT survives does not depend on the passage times of $\Fo$, but just on the location of the seeds.

Similarly, we couple the evolution of $\Fl $ on $\T_d$ from each seed activated by $\Fo$ as independent BRWCT on $\mathbb{R}$ by setting the offspring distribution to be the constant $d$, and the distribution governing the underlying random walk to have (independent) increments as exponentials of parameter $\gamma=\lambda$.
Clearly, individuals here correspond to vertices of type $\Fl $, and each progeny is started by a seed.
The theorem below is the concentration result we shall use. 
The measure $\P$ below is defined on the space of the trajectories of BRWCT on $\mathbb{R}$, where the branching mechanism and the step size are governed by the prescribed distributions.

\begin{theorem}[{\cite[Theorem 3, adapted to our setting]{AddarioBerry-Reed}}]\label{thm:ABR}
Consider a supercritical BRWCT on $\mathbb{R}$ with branching distribution Bin($d,\, \ns$) and step size $X\sim $ Exp($\gamma$) (for a fixed $\gamma>0$) coupled with FPPHE as specified above. 
Then, for each fixed 
%
$x\in \mathcal{L}_i$ there exist constants $C_x > 0, \delta_x > 0$ which depend on $\gamma$, such that for all $\alpha>0 $ and all $n\geq1$
\[
\P \left ( \left | M_{i,n}^{(\gamma)}(x)-\E\left (  M_{i,n}^{(\gamma)} (x)\mid \mathscr{S} \right ) \right | > \alpha \mid \mathscr{S} \right )\leq C_x e^{-\delta_x \alpha },
\]
where $M_{i,n}^{(\gamma)}(x)$ was defined in \eqref{eq:defMx}.
\end{theorem}
\begin{remark}
It follows from the definition that, whenever the branching distribution is Bin($d,\, \ns$), the measure $\P$ above can be coupled with $\P_\mu$.
From now on, for every measurable event $E$ define
$
\wtP_{\den}[E]:=\P_{\den} [E \mid \mathscr{S}]
$,
that is, the probability measure $\wtP_{\den}$ coincides with ${\P}_{\den}$ conditioned on the event of survival $\mathscr{S}$.
Consistently, the corresponding expectation will be denoted by $\widetilde{\E}_{\den}$.
\end{remark}

\subsection{Proof of Claim \ref{claim:to-prove-later}}\label{sect:proof-claim-later}

Before continuing with the proof of Proposition \ref{prop:Fo-wins}, in this short section we use the auxiliary results stated above to establish  Claim \ref{claim:to-prove-later}, in particular we will only write the proof for the lower part, as it is completely analogous for the upper part.

\begin{proof}[Proof of Claim \ref{claim:to-prove-later}]
Start by defining the following auxiliary event
\[
\begin{split}
& \overline{\mathscr{E}}:=\{\exists \text{ a path $\pi$ from }\Olow \text{ to } \Wlow \text{ that is free of seeds and s.t.\ }\Upsilon_{\text{low}}^{(1)}(\pi)\leq  \frac{H+1}{\cin^{(2)}} \}.
\end{split}
\]
We can bound the sought probability from below as follows,
\[
\begin{split}
& 
\P_{\mu_1} \left [  \cE_{\Olow} \cap  \mathfrak{B}_{\text{low}}(\Olow,\Wlow)\cap \left \{  \Upsilon_{\text{low}}^{(1,\lambda)}(O_{\text{low}} \to \Wlow)\leq \frac{H+1}{\cin^{(2)}} \right \} \right ] \\
&  \geq \P_{\mu_1} \left [  \cE_{\Olow} \cap  \mathfrak{B}_{\text{low}}(\Olow,\Wlow)\cap \overline{\mathscr{E}} \cap \mathscr{S}\right ]\\
& = \P_{\mu_1} \left [  \cE_{\Olow} \cap  \mathfrak{B}_{\text{low}}(\Olow,\Wlow) \cap \mathscr{S} \right ] - \P_{\mu_1} \left [  \cE_{\Olow} \cap \mathfrak{B}_{\text{low}}(\Olow,\Wlow)\cap \overline{\mathscr{E}}^c \cap \mathscr{S} \right ].
\end{split}
\]
Now let $\mathscr{S}_H^{\text{not}}$ denote the event that the GWT survives for $H+1$ generations but it does not survive forever.
At this point, the above expression is at least
\[
\begin{split}
& \P_{\mu_1} \left [  \cE_{\Olow} \cap \mathfrak{B}_{\text{low}}(\Olow,\Wlow) \cap \mathscr{S} \right ] - \P_{\mu_1} \left [ \overline{\mathscr{E}}^c \cap \mathscr{S} \right ] \\
& \geq  \P_{\mu_1} \left [  \cE_{\Olow} \cap \mathfrak{B}_{\text{low}}(\Olow,\Wlow) \cap \mathscr{S} \right ] - \P_{\mu_1} \left [ \overline{\mathscr{E}}^c \mid \mathscr{S} \right ] \\
& =  \P_{\mu_1} \left [  \cE_{\Olow} \cap \mathfrak{B}_{\text{low}}(\Olow,\Wlow) \right ]-  \P_{\mu_1} \left [  \cE_{\Olow} \cap  \mathfrak{B}_{\text{low}}(\Olow,\Wlow) \cap \mathscr{S}^c \right ] - \wtP_{\mu_1} \left [ \overline{\mathscr{E}}^c  \right ] \\
& \geq \P_{\mu_1} \left [  \cE_{\Olow} \cap \mathfrak{B}_{\text{low}}(\Olow,\Wlow) \right ]-  \P_{\mu_1} \left [ \mathscr{S}_H^{\text{not}} \right ] - \wtP_{\mu_1} \left [ \overline{\mathscr{E}}^c  \right ].
\end{split}
\]
By Proposition \ref{prop:Fo-wins} the first term in the last line is bounded from below by $(1-\f_2(\mu_1)-\bseta_2)(1-\mu_1)^2 $.
Regarding the second term, it suffices to recall that a super-critical GWT conditioned on dying out has the same distribution of a \emph{sub-critical} Galton-Watson tree (cf.\ \cite[Section I.12]{athreya-ney-BP}).
Thus, the event $\mathscr{S}_H^{\text{not}}$ is so that there is a  constant $\beta>0$ dependent on $\mu_1$ but independent of $H$, so that
\[
\P_{\mu_1} \left [ \mathscr{S}_H^{\text{not}} \right ]\leq e^{-\beta H}.
\]
In order to bound the third term we proceed as follows.
Since the event $\Upsilon_{\text{low}}^{(1)}(\pi)> \frac{H+1}{\cin^{(2)}}$ must hold for all paths of length $H+1$ starting at $\Olow$ that are free of seeds, it holds for the one with minimum passage time.
By construction, it is easy to see that there is a constant $C'>0$ (independent of $H$) so that
\[
\wtE_\mu \left [  M_{0,H+1}^{(1)}\right ]\leq C' (H+1).
\]
Furthermore, we shall assume that the constant $\cin^{(2)}$ in the statement is so that $ 1/\cin^{(2)}\geq 2 C'$.
Thus,
\[
\begin{split}
\wtP_{\mu_1} \left [ \overline{\mathscr{E}}^c  \right ]
& \leq \wtP \left ( \left | M_{0,H+1}^{(1)}(O_{\text{low}})-\wtE\left (  M_{0,H+1}^{(1)} (O_{\text{low}})  \right ) \right | > \left (\frac{1}{\cin^{(2)}}-C'\right )(H+1)  \right )\\
& \leq \wtP \left ( \left | M_{0,H+1}^{(1)}(O_{\text{low}})-\wtE\left (  M_{0,H+1}^{(1)} (O_{\text{low}}) \right ) \right | >  \frac{H+1}{2\cin^{(2)}}  \right )\\
& \stackrel{\text{Thm \ref{thm:ABR} } }{\leq }C_{O_{\text{low}}} e^{-\delta_{O_{\text{low}}} (H+1)/(2\cin^{(2)}) },
\end{split}
\]
concluding the proof.
\end{proof}

\subsection{Step 1, taking care of levels up to $ (1-\uep{e0})m $}\label{sect:step1}
We start by showing that it is unlikely that $\Fl$ originated from any seed located far away from $\mathcal{L}_m$ reaches level $m$ too quickly.
To achieve this we show that the probability of the occurrence of a bad event $B_i$ (cf.\ \eqref{eq:def-event-Bi}) at some level much earlier than $m$ decreases exponentially fast.

Now we proceed with our analysis as follows.
Fix an arbitrarily small constant $\uep{e0}>0$; for any fixed $i\leq  (1-\uep{e0})m $ the event $B_i$  necessarily implies that there must be a vertex $\mathbf{x}_i\in \mathcal{L}_i$ hosting a seed (and whose tree-geodesic to the origin is free of seeds), which starts a $\Fl$ process such that
\[
M_{i,m}^{(\lambda)}(\mathbf{x}_i)< M_{0,m}^{(1)} - \Upsilon_{\oT_d^m}^{(1)} (o \to \mathbf{x}_i).
\]
For simplicity, let $\M(\mathbf{x}_i)$ denote this event, then
\begin{equation}\label{eq:Bi-in-Mx}
B_i \subseteq \left \{ \exists \text{ a seed } \mathbf{x}_i \in \mathcal{L}_i \ : \ \M(\mathbf{x}_i) \text{ is realized} \right \} .
\end{equation}
We shall proceed with bounding the probability of the RHS of the above.
We shall use the obvious fact that $\M(\mathbf{x}_i) \Rightarrow M_{i,m}^{(\lambda)}(\mathbf{x}_i)< M_{0,m}^{(1)}$.
%
Fix a large constant $\uc{ch}>0$, and define the following auxiliary event
\[
\A_{(1)}:= \left \{ M_{0,m}^{(1)}<(1+\uc{ch}) \wtE_\mu \left [ M_{0,m}^{(1)}\right ] \right \} .
\]
The value $ M_{0,m}^{(1)}$ is a minimum over all paths starting at the origin as $\Fo$ and arriving to $\mathcal{L}_m$ as $\Fo$. Thus it is easy to see that there is a constant $C'>0$ (independent of $m$) so that
\begin{equation}\label{eq:E-M0m1}
\wtE_\mu \left [  M_{0,m}^{(1)}\right ]\leq C' m.
\end{equation}
Similarly, there must be a constant $\mathcal{C}^\ast >0$ independent of $\lambda$, such that for all $0\leq i<n\leq m$
\begin{equation}\label{eq:definitionClambda}
\E_\mu \left [ M_{i,n}^{(\lambda)}  \right ] \geq \frac{\mathcal{C}^\ast }{\lambda} (n-i).
\end{equation}
We shall proceed by showing that 
$ \wtP_\mu \left [ \exists \, i\leq  (1-\uep{e0})m \text{ s.t. } \exists \text{ seed }\mathbf{x}_i \, : \, \M(\mathbf{x}_i) \text{ occurs} \right ] $
decreases exponentially fast in $m$.

\nc{ch}
\nc{c}

\begin{claim}\label{claim:bd0}
Given $\uep{e0}>0$, for all $\lambda>0$ small enough, there is a constant $\boldsymbol{\alpha_1}=\boldsymbol{\alpha_1}(\lambda, \uep{e0})>0$ independent of $m$ such that
\[
\wtP_\mu \left [ \exists ~ i \leq  (1-\uep{e0})m  \text{ s.t. } \exists \text{ seed }\mathbf{x}_i \ : \ \M(\mathbf{x}_i) \cap \A_{(1)} \right ] \leq e^{- \boldsymbol{\alpha_1} m }.
\]
\end{claim}
\begin{proof}
By definition, the event $\M(\mathbf{x}_i) \cap \A_{(1)} $ implies
\begin{equation}\label{eq:MAA2}
\begin{split}
& M_{i,m}^{(\lambda)}(\mathbf{x}_i)< M_{0,m}^{(1)}  \leq (1+\uc{ch}) \wtE_\mu \left [ M_{0,m}^{(1)}\right ].
\end{split}
\end{equation}
Thus, relation \eqref{eq:MAA2} together with \eqref{eq:E-M0m1} imply
\[
\begin{split}
M_{i,m}^{(\lambda)}(\mathbf{x}_i) 
& \leq (1+\uc{ch})C'm
= (1+\uc{ch})C'(m-i)+(1+\uc{ch})C'i\\
& \stackrel{ i\leq  (1-\uep{e0})m  }{\leq} (1+\uc{ch})C'(m-i)+(1+\uc{ch})C' (1-\uep{e0})m  \\
& = (1+\uc{ch})C'(m-i)+ \frac{(1+\uc{ch})C'}{\uep{e0}}\uep{e0} (1-\uep{e0})m   \ 
\stackrel{\uep{e0}m\leq m-i}{< } \ C''( m - i),
\end{split}
\]
for some constant $C''>0$ which depends on $\uep{e0}$ but is independent of $i$ and $m$.
By \eqref{eq:definitionClambda},
\begin{equation}\label{eq:cond12}
\wtE_\mu \left [ M_{i,m}^{(\lambda)}(\mathbf{x}_i) \right ] \geq \frac{\mathcal{C}^\ast }{\lambda}(m-i).
\end{equation}
By putting these facts together we find that 
\begin{equation}\label{eq:cond22}
 M_{i,m}^{(\lambda)}(\mathbf{x}_i)< C''( m - i)<\frac{C''}{\mathcal{C}^\ast} \lambda \wtE_\mu \left [ M_{i,m}^{(\lambda)}(\mathbf{x}_i) \right ].
\end{equation}
Since $\lambda$ is small enough, relations \eqref{eq:cond12} and \eqref{eq:cond22} ensure that there is a constant $\uc{c}>0$ such that uniformly over all $i\leq  (1-\uep{e0})m $,
\[
M_{i,m}^{(\lambda)}(\mathbf{x}_i)<\wtE_\mu \left [  M_{i,m}^{(\lambda)}(\mathbf{x}_i)\right ] - \uc{c}\frac{1}{\lambda}(m-i) .
\]
In particular, we have that
\begin{equation}\label{eq:B_i-to-ABR2}
\wtP_{\den}
\left [ \M(\mathbf{x}_i) \cap \A_{(1)}\right ]  
\leq d^i \max_{\mathbf{x}_i \in \mathcal{L}_i} \wtP_{\den} \left ( \left | M_{i,m}^{(\lambda)}(\mathbf{x}_i)-\wtE_\mu \left [ M_{i,m}^{(\lambda)}(\mathbf{x}_i)\right ] \right |>\frac{ \uc{c}}{\lambda}(m-i) \right ),
\end{equation}
where in the last step we used a union bound over all vertices of $\mathcal{L}_i$.
Thus,
\[
\begin{split}
& \wtP_{\den}  \left ( \exists ~ i \leq \lceil (1-\uep{e0})m \rceil \text{ s.t. } \exists \text{ seed }\mathbf{x}_i \in \mathcal{L}_i \ : \  M_{i,m}^{(\lambda)}(\mathbf{x}_i)<\wtE_\mu \left [  M_{i,m}^{(\lambda)}(\mathbf{x}_i)\right ] - \frac{\uc{c}}{\lambda}(m-i)
\right ) \\
&  \stackrel{\eqref{eq:B_i-to-ABR2} }{\leq }\sum_{i=1}^{\lceil (1-\uep{e0})m \rceil} d^i \max_{\mathbf{x}_i \in \mathcal{L}_i}  \wtP_{\den}\left ( \left | M_{i,m}^{(\lambda)}(\mathbf{x}_i)-\wtE_\mu \left [  M_{i,m}^{(\lambda)} (\mathbf{x}_i)\right ] \right | > \frac{\uc{c}}{\lambda}(m-i) \right )\\
& \stackrel{ \text{Thm.\ \ref{thm:ABR} } }{\leq }
\sum_{i=1}^{\lceil (1-\uep{e0})m \rceil} {C} d^i \, e^{-{\delta }(m-i) \uc{c}/\lambda } \leq e^{- \boldsymbol{\alpha _1}  m},
\end{split}
\]
where the very last inequality follows from the facts that $\lambda$ is small enough 
(in order to have $d^m \, e^{-{\delta }\uep{e0} m \uc{c}/\lambda }<1$ to ensure summability) 
and $i\leq \lceil (1-\uep{e0})m \rceil$.
\end{proof}

The following result is of the same fashion.

\begin{claim}\label{claim:bd3}
Given $\uep{e0}>0$, there is a constant $\boldsymbol{\alpha_2}=\boldsymbol{\alpha_2}(\uep{e0})>0$ independent of $m$ so that
\[
\wtP_\mu \left [   \A_{(1)}^c \right ] \leq e^{-\boldsymbol{\alpha_2} m}.
\]
\end{claim}
\begin{proof}
By definition of $\A_{(1)}$, we obtain
\[
\wtP_\mu \left [ \A_{(1)}^c \right ] \leq 
\wtP_\mu \left [ \left | M_{0,m}^{(1)} -\wtE_\mu \left [ M_{0,m}^{(1)}\right ] \right | \geq \uc{ch} \wtE_\mu \left [ M_{0,m}^{(1)}\right ]  \right ] \stackrel{\text{Thm \ref{thm:ABR} } }{\leq}
 C e^{-\delta \uc{ch} \wtE_\mu \left [  M_{0,m}^{(1)}\right ]}.
\]
The quantity $ M_{0,m}^{(1)}$ is a minimum over all paths starting at the origin as $\Fo$ and arriving to $\mathcal{L}_m$ as $\Fo$, therefore it is easy to see that there is a constant $\mathcal{C}>0$ such that
\[
\wtE_\mu \left [  M_{0,m}^{(1)}\right ]\geq \mathcal{C} m,
\]
implying the statement for a suitable choice of $\boldsymbol{\alpha_2}>0$.
\end{proof}

\begin{lemma}\label{lemma:B_i-unlikely}
Given $\uep{e0}>0$, for all small enough $\lambda>0$ there exists a constant $\delta'=\delta'(\lambda,\uep{e0})>0$ independent of $m$, such that
\[
\wtP_{\den} \left ( \bigcap_{i=1}^{ (1-\uep{e0})m } B_i^c\right )\geq 1-e^{- \delta' m}.
\]
\end{lemma}
\begin{proof} 
For all $i\leq  (1-\uep{e0})m $, for each seed $\mathbf{x}_i\in \mathcal{L}_i$ we have
\begin{equation}\label{eq:tobound}
\wtP_\mu \left [ \M(\mathbf{x}_i) \right ]  = \wtP_\mu \left [ \M(\mathbf{x}_i) \cap \A_{(1)} \right ] + \wtP_\mu \left [ \M(\mathbf{x}_i) \cap \A_{(1)}^c \right ].
\end{equation}
The proof is now a simple consequence of the previous two claims, in fact by \eqref{eq:Bi-in-Mx} 
\[
\begin{split}
\wtP_\mu \left [ \bigcup_{i=1}^{  (1-\uep{e0})m } B_i \right ] 
& \leq \wtP_\mu \left [ \bigcup_{i\leq  (1-\uep{e0})m }  \left \{ \exists \text{ a seed } \mathbf{x}_i \in \mathcal{L}_i \ : \ \M(\mathbf{x}_i) \text{ is realized} \right \} \right ]\\
& \stackrel{\eqref{eq:tobound} \text{ and Claims \ref{claim:bd0},  \ref{claim:bd3} }}{\leq }e^{-\boldsymbol{\alpha _1} m}+ e^{-\boldsymbol{\alpha _2} m}.
\end{split}
\]
Since $ \boldsymbol{\alpha _1}, \boldsymbol{\alpha _2}$ are independent of $m$, the statement follows.
\end{proof}

\subsection{Step 2, subsequent levels}
In this section we will bound the probability of the event
\[
\bigcup_{i= (1-\uep{e0})m  +1}^m B_i.
\]
Since this part is more delicate, we need to control the probabilities of the $B_i$'s in a more refined way.
%
%
\nC{C1}
\nep{e1}
\nep{e1bis}
We proceed by using the BRWCT representation that we defined before.
Recall that by  ``GWT'' we mean the Galton-Watson  tree with offspring distribution given by Bin($d,1-\mu$).
For any individual $u$ in the GWT define
\begin{equation}\label{eq:def-tu}
\bt_u:= \text{birth time of }u,
\end{equation}
where the root $o$ was born at time $0$, that is $\bt_o=0$.
From the definition of BRWCT we have
\[
\bt_u \stackrel{\text{law}}{=}\Upsilon_{\oT_d^m}^{(1,\lambda)}(o \to u).
\]
To make the proof more readable, from now on we set 
\begin{equation}\label{eq:def-ell}
\ell:= \lfloor \uep{e0} m \rfloor.
\end{equation}
Now for each $j\in \{1, \ldots ,  \ell \}$ define
\begin{equation}\label{eq:def-Nj}
\mathcal{N}_j:=\{\text{individuals at }j\text{-th generation of the GWT}\},
\end{equation}
and let $N_j:= |\mathcal{N}_j|$.
Consider a large constant $\uC{C1}>0$ (which will be specified later) and set 
\begin{equation}\label{eq:def-Kj}
\mathcal{K}_j:=\{u\in \mathcal{N}_j \ : \ \bt_u\leq \uC{C1}j\},
\end{equation}
and let $K_j:=|\mathcal{K}_j|$.
In words, $\mathcal{K}_j$ is the subset of $\mathcal{N}_j$ of individuals born within time $\uC{C1}j$.

Recall the technical conditions on $\mu$ from \eqref{eq:techCondMu} and fix two small constants $\uep{e1}>0$ and $\uep{e1bis}>0$.
With this notation in mind, we define a sequence of ``good'' events, namely for $j\in \{1, \ldots ,  \ell \}$ set 
\nc{c2} \nC{C2}
\begin{equation}\label{eq:def-Aj}
A_j:=\left \{K_j\geq \uep{e1bis} \bigl (d\ns \bigr )^{(1-\uep{e1})j} \right \}.
\end{equation}
In order to continue, we need an auxiliary result.
Recall that $\wtP_\mu (\cdot)=\P(\cdot \mid \mathscr{S})$ is the probability measure conditioned on the event of survival of GWT.

\begin{lemma}\label{lemma:new-concentr}
Assume that $d$ and $\mu$ satisfy \eqref{eq:techCondMu} and $ \uep{e1}, \uep{e1bis} $ are as above.
Then there are two positive constants $\uc{c2}=\uc{c2}(\uep{e1})<\frac{\uep{e1}}{2} \log (d(1-\mu))$ and $ \uC{C2}= \uC{C2}(\uep{e1}, \uep{e1bis})$ such that for all $1 \leq j \leq \ell$,
\begin{equation}\label{eq:bound-concj}
\wtP_\mu \left ( \uep{e1bis} \bigl (d\ns \bigr )^{(1-\uep{e1})j} \leq N_j \leq \frac{1}{\uep{e1bis}} \bigl (d\ns \bigr )^{(1+\uep{e1})j}\right )\geq 1-
\uC{C2}e^{-\uc{c2}j}.
\end{equation}
In particular, $ \uC{C2}$ can be made arbitrarily small by choosing 
$ \uep{e1bis}$ small enough.
\end{lemma}
\begin{proof}
For all $j \in \{1, \ldots , \ell\}$ we have
\[
\begin{split}
& \wtP_\mu \left ( \frac{ N_j}{\bigl (d\ns \bigr )^{j}} \geq \frac{1}{\uep{e1bis}} \bigl (d\ns \bigr )^{\uep{e1}j}\right ) 
\leq  \frac{\P_\mu \left ( \frac{ N_j}{\bigl (d\ns \bigr )^{j}} \geq \frac{1}{\uep{e1bis}} \bigl (d\ns \bigr )^{\uep{e1}j} \right ) }{\P_\mu(\mathscr{S})} \\
& \leq \frac{1}{\P_\mu(\mathscr{S})} \frac{\E_\mu \left [ \frac{ N_j}{\bigl (d\ns \bigr )^{j}}\right ]\uep{e1bis} }{\bigl (d\ns \bigr )^{\uep{e1}j}}= \frac{1}{\P_\mu(\mathscr{S})} \frac{\uep{e1bis} }{\bigl (d\ns \bigr )^{\uep{e1}j}} \leq \frac{\uC{C2}e^{-\uc{c2}j}}{2}.
\end{split}
\]
The above inequalities follow from Markov's inequality, from the fact that $\E_\mu N_j=\bigl (d\ns \bigr )^j$ and that $\uc{c2}<\uep{e1}\log (d(1-\mu))$.

To bound the reverse inequality we re-write the sought quantity to apply Markov's inequality,
\begin{equation}\label{eq:bound-to-go}
\begin{split}
& \wtP_\mu \left ( \frac{N_j}{\bigl (d\ns \bigr )^j }\leq  \uep{e1bis} \bigl (d\ns \bigr )^{-\uep{e1}j}  \right ) 
= \wtP_\mu \left ( \frac{\bigl (d\ns \bigr )^j}{N_j}\geq \frac{1}{ \uep{e1bis} }\bigl (d\ns \bigr )^{\uep{e1}j} \right )\\
& \leq \wtE \left [ \frac{\bigl (d\ns \bigr )^j}{N_j}\right ] \uep{e1bis} \bigl (d\ns \bigr )^{-\uep{e1}j} 
\leq \exp \left \{ \left | \log \wtE \left [ \frac{\bigl (d\ns \bigr )^j}{N_j}\right ] \right | \right \} \uep{e1bis} \bigl (d\ns \bigr )^{-\uep{e1}j}.
\end{split}
\end{equation}
Hence our next aim is to bound the expression \ 
$
\exp \left \{ \left | \log \wtE \left [ \frac{\bigl (d\ns \bigr )^j}{N_j}\right ] \right | \right \}.
$
To do so, we employ \cite[Corollary 1]{Ney-Vidyashankar}, which says that under the so-called $L \log L$ condition (which is satisfied here, since the offspring distribution of GWT is binomial) we have that
\[
\lim_n \frac{1}{n} \log \wtE \left [ \frac{1}{N_n}\right ]=\max \left \{ \log(p_1), \, - \log (\E (N_1)) \right \},
\]
where $p_1:=\P(N_1=1)$.
In our specific case this translates into
\[
\lim_n \frac{1}{n} \log \wtE \left [ \frac{1}{N_n}\right ]=\max \left \{ \log \left ( d(1-\mu)\mu^{d-1}\right ), \, - \log (d(1-\mu)) \right \}.
\]
Given that $d$ and $\mu$ must satisfy \eqref{eq:techCondMu}, we obtain that $\log \left ( d(1-\mu)\mu^{d-1}\right )< - \log (d(1-\mu))$, implying that for any fixed $\varepsilon'>0$, for all $n$ sufficiently large we have
\[
\left | \frac{1}{n} \log \wtE \left [ \frac{1}{N_n}\right ]+ \log (\E N_1) \right |< \varepsilon'.
\]
From now on, for convenience, we set $\varepsilon':=\frac{\uep{e1}}{2}\log(\E N_1)$.
Thus, for all $n$ large enough, the above can be rewritten as 
\[
\exp \left \{\left | \log \wtE \left [ \frac{1}{N_n} (\E N_1)^n\right ] \right |\right \}< \exp \{ \varepsilon' n\} \  \stackrel{\varepsilon'=\uep{e1} \log(\E N_1)/2}{=} \ (d(1-\mu))^{\uep{e1} n /2}.
\]
Hence, there is a constant $C'>0$ such that uniformly for all $j\geq 1$ we have
\begin{equation}\label{eq:boundJ}
\exp \left \{\left | \log \wtE \left [ \frac{1}{N_j} (\E N_1)^j\right ] \right |\right \}< C'(d(1-\mu))^{\uep{e1} j /2}.
\end{equation}
Thus, replacing this quantity in \eqref{eq:bound-to-go} we obtain
\[
\begin{split}
& \wtP_\mu \left [ \frac{N_j}{\bigl (d\ns \bigr )^j }\leq  \uep{e1bis} \bigl (d\ns \bigr )^{-\uep{e1}j}  \right ] 
 \stackrel{\eqref{eq:bound-to-go}, \, \eqref{eq:boundJ} }{\leq} C'(d(1-\mu))^{\uep{e1} j /2}\uep{e1bis} \bigl (d\ns \bigr )^{-\uep{e1}j},
\end{split}
\]
which equals $\, C' \uep{e1bis} \bigl (d\ns \bigr )^{-\uep{e1}j/2}$.
Since $\uc{c2}<\frac{\uep{e1}}{2} \log (d(1-\mu))$ then, for all $j \in \{1, 2, \ldots, \ell\}$, 
\[
\wtP_\mu \left ( \frac{N_j}{\bigl (d\ns \bigr )^j }\leq  \uep{e1bis} \bigl (d\ns \bigr )^{-\uep{e1}j}  \right ) 
\leq C' \uep{e1bis} \bigl (d\ns \bigr )^{-\uep{e1}j/2}\leq \frac{\uC{C2}e^{-\uc{c2}j}}{2},
\]
finishing the proof.
In particular, $\uC{C2}$ can be made arbitrarily small by picking $\uep{e1bis}$ small enough.
\end{proof}

Here we recall an auxiliary result that can be found e.g.\ in  \cite{Janson-concentration}, which will allow us to proceed.
The following statement and Lemma \ref{lemma:typical-FPP} are strongly related.
In the result below we denote by $\P$ the probability measure defined on all possible realizations of the sums of independent random variables with the prescribed distribution.
\begin{theorem}\label{thm:aux-janson}{\cite[Thm 5.1(i),(iii)]{Janson-concentration}}
Let $ X:= \sum_{k=1}^K X_k$, where $\{X_k\}_k $ is a collection of independent random variables such that $X_k\sim $ Exp($a_k$) for some constants $a_k>0$ and let $a_\ast:=\min_k a_k$. Then the following bounds hold.
\begin{itemize}
\item[(i)] For all fixed $\delta>0$ we have
$
\P \bigl (X\geq (1+\delta)\E X \bigr )\leq \frac{1}{1+\delta} \exp \bigl \{-a_\ast \E X (\delta - \ln (1+\delta)) \bigr \}
$.
\item[(ii)] For all fixed $0<\delta<1$ we have
$
\P \bigl (X\leq (1-\delta)\E X \bigr )\leq \exp \bigl \{-a_\ast \E X (-\delta - \ln (1-\delta)) \bigr \}
$.
\end{itemize}
\end{theorem}

\begin{claim}\label{claim:NiBorn}
For every $j\in \{1, \ldots , \ell\}$ 
\[
\wtP_{\den} \left [\text{The entire generation }\mathcal{N}_j\text{ is born by time } \uC{C1}j \right ]  \geq  1-\frac{1}{\uC{C1}}e^{- j \uC{C1}/2} .
\]
\end{claim}
\begin{proof}
By Theorem \ref{thm:aux-janson}, with high probability all particles of generation $j$ of GWT are born within time $\uC{C1} j$, for some $\uC{C1}$ large enough.
Clearly, for all $j\in \{1, \ldots , \ell \}$, 
\[
\begin{split}
& \wtP_{\den} \left [\text{Entire generation }\mathcal{N}_j\text{ is born by time } \uC{C1}j \right ] 
\geq  \ 1- \wtP_{\den}\bigl [ \exists \text{ individual in }\mathcal{N}_j\text{ born after time }  \uC{C1}j\bigr ].
\end{split}
\]
Now, let $\{X_k\}_k$ denote an i.i.d.\ sequence of exponential random variables with parameter $1$.
By a union bound we obtain that the above probability is bounded from below by
\[
\begin{split}
1- d^{j}\P \left [ \sum_{k=1}^j X_k \geq \uC{C1}j \right ] 
& \stackrel{\text{Thm \ref{thm:aux-janson}(i)} }{\geq } 1-d^{j}\frac{1}{\uC{C1}}e^{- j \bigl (\uC{C1}-1-\ln(\uC{C1}) \bigr )} \\
& =  1-\frac{1}{\uC{C1}}e^{- j \bigl (\uC{C1}-1-\ln(\uC{C1}) -\ln(d) \bigr )} 
 \ \ \stackrel{ \uC{C1} \text{ large}}{\geq } \ \ 1-\frac{1}{\uC{C1}}e^{- j \uC{C1}/2} .
\end{split}
\]
This concludes the proof.
\end{proof}

The next result proves that for all $j\in \{1, \ldots , \ell \} $, it is likely that at level $j$ there are \emph{a lot} of vertices of type $\Fo$ which have been infected by time $\uC{C1} j$.

\nep{e2}
\begin{lemma}\label{lemma:A^c-unlikely}
Given any small constant $\uep{e2}>0$, we can choose $\uep{e1}>0$ in \eqref{eq:def-Aj} such that 
\[
\wtP_{\den} \left (\exists \, j\in \{1, \ldots, \ell \}\text{ s.t.\ } A_j^c\text{ occurs} \right )\leq \uep{e2}.
\]
\end{lemma}
\begin{proof}
The probability $\wtP_\mu \left  (K_j\geq \uep{e1bis} ( d\ns )^{(1-\uep{e1})j} \right )$ is bounded from below by
\[
\begin{split}
& \wtP_\mu \left [ \left \{ \text{Entire }\mathcal{N}_j\text{ born by time } \uC{C1}j \right \} \cap \left \{ N_j \geq \uep{e1bis} ( d\ns )^{(1-\uep{e1})j} \right \}\right ]\\
& \geq 1- \wtP_\mu \left [ \left \{  \text{Entire } \mathcal{N}_j\text{ born by time } \uC{C1}j \right \}^c \right ] - \P \left [ N_j <  \uep{e1bis} ( d\ns )^{(1-\uep{e1})j} \right ].
\end{split}
\]
The last term of the sum is taken care of by Lemma \ref{lemma:new-concentr}, while for the middle term we rely on Claim \ref{claim:NiBorn}.
Using these facts we obtain that for every $j\in \{1, \ldots , \ell \}$ 
\[
\wtP_{\den}(A_j^c)\leq \frac{1}{\uC{C1}}e^{- j \uC{C1}/2} +\uC{C2}e^{-\uc{c2}j}.
\]
Hence it suffices to pick $\uC{C1}=\uC{C1}(d,\uep{e2})$ large enough and $\uep{e1bis}=\uep{e1bis}(\uep{e2}), ~ \uep{e1}=\uep{e1}(\uep{e2})$ so that
\begin{equation}\label{eq:e2/2}
\sum_{j=1}^{\ell } \frac{1}{\uC{C1}}e^{- j \uC{C1}/2}\leq \frac{\uep{e2}}{2} \quad \text{and} \quad \sum_{j=1}^{ \ell } \uC{C2}e^{-\uc{c2}j}\leq \frac{\uep{e2}}{2},
\end{equation}
in order to obtain $\sum_{j=1}^{\ell} \wtP_{\den} (A_j^c) \leq \uep{e2} $, which is the claim.
\end{proof}

\begin{figure}[h!]
\begin{center}
\includegraphics[width=10cm]{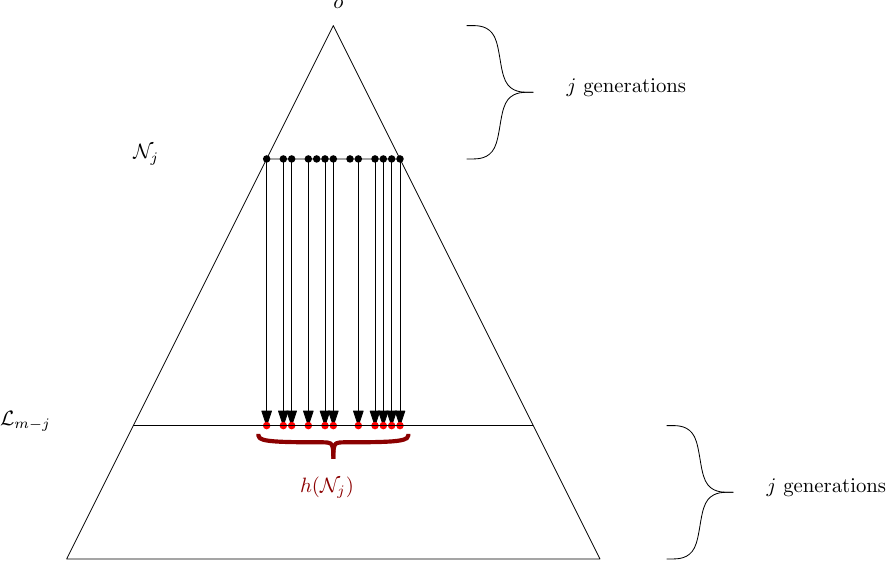}
\end{center}
\caption{Sketch of the construction of $h(\mathcal{N}_j)$.}\label{fig:h(Nj)}
\end{figure}

For a $ B_i $ to be realized, there must be a first seed (at level $\mathcal{L}_i$) initiating the $\Fl$ process which infects level $\mathcal{L}_m$ before any other process does.
From now on we set
\begin{equation}\label{eq:defSeedS}
\s:=\text{first seed initiating  the $\Fl$ process which infects the first vertex on }\mathcal{L}_m.
\end{equation}
In other words, $\cup_{i=m-\ell}^m B_i$ is realized only if there is a seed $\s\in \cup_{j=0}^\ell \mathcal{L}_{m-j} $ such that
\[
\min_{y \in \mathcal{L}_m}\Upsilon^{(1,\lambda)}_{\oT_d^{m}}(o\to y)=\Upsilon^{(1)}_{\oT_d^{m}}(o \to \s )+ \min_{y \in \mathcal{L}_m} \Upsilon^{(\lambda)}_{\oT_d^{m}}(\s\to y) .
\]
\begin{remark}
Note that $j=0$ means that the first infected vertex on $\mathcal{L}_{m}$ is a seed, hence from now on we shall consider $j\geq 1$.
\end{remark}
Recall that for any vertex $y\in V(\oT^m_{d})\setminus \{o\}$ the \emph{father} of $y$ is the unique vertex $y^-$ on the (unique) geodesic line connecting $o$ to $y$ such that 
\[
d_{\T_{d}} (o,y)=d_{\T_{d}} (o,y^-)+1.
\]
For all $j\in \{1, \ldots , \ell \}$ fixed, for every vertex $v\in \mathcal{N}_j $ with birth time $\bt_v\leq \uC{C1}j$, let the $\Fo$ process started at $v$ continue until it reaches level $\mathcal{L}_{m-j}$ for the first time.
More precisely, define the set
\begin{equation}\label{eq:defS(m-j)}
\mathcal{S}(m-j):=\{y\in \mathcal{L}_{m -j} \ : \ y^-\in \mathcal{N}_{m-j-1}\},
\end{equation}
that is, $\mathcal{S}(m-j)$ is the set of vertices $y\in \mathcal{L}_{m -j}$ whose father is of type $\Fo$. 
Thus, for every $v\in \mathcal{N}_{j}$ we look for  a vertex  $h(v) \in \mathcal{S}(m-j)$, that corresponds to the first vertex that has been infected by $\Fo$ and has $v$ as ancestor.
(Note that $h(v)$ can be of type $\Fo$ or a seed.)
Furthermore, it will be convenient to set
\begin{equation}\label{eq:def-h(set)}
h(\mathcal{N}_j):= \cup_{v \in \mathcal{N}_j }h(v),
\end{equation}
keeping in mind that $h(\mathcal{N}_j)\subseteq \mathcal{S}(m-j)$; cf.\ Figure \ref{fig:h(Nj)}.

The next result shows that with high probability all vertices $h(v)$ are born relatively close (in time) to each other.

\nC{C5}
\begin{lemma}\label{lemma:aux}
For any large enough constant $\uC{C5}'$, we can take $\uC{C5}=\uC{C5}(\uC{C5}')$ so large that 
%
\[
\wtP_{\den}  \left [ \sup_{v,v' \in \mathcal{N}_j}\left |  \bt_{h(v)} -  \bt_{h(v')}\right | \leq \bigl (2\uC{C5} + \uC{C1} \bigr ) j \right ]\geq 1- 2e^{-\uC{C5}' j} -\frac{1}{\uC{C1}}e^{- j \uC{C1}/2}.
\]
\end{lemma}
\begin{proof}
For simplicity, denote
\begin{equation}\label{eq:def-tau}
\tau \bigl (v, h(v)\bigr ):= \bt_{h(v)}-\bt_v,
\end{equation}
that is the time needed by $\Fo$ started at $v$ to advance up to infecting  $\mathcal{S}(m-j)$ for the first time.
Note that by definition we have that the variables $\{ \tau \bigl (v, h(v)\bigr )\}_{v\in \mathcal{N}_j}$ are independent. 
%
Furthermore, since the process advances as $\Fo$ for $m-2j$ steps, we also have that $\{ \tau \bigl (v, h(v)\bigr )\}_{v\in \mathcal{N}_j}$ are identically distributed.
In particular, 
\begin{equation}\label{eq:tauLaw}
\text{for all }v\in \mathcal{N}_j, \quad 
\tau \bigl (v, h(v)\bigr ) \stackrel{\text{law}}{=} M^{(1)}_{0,m-2j}.
\end{equation}
We can exploit this fact by fixing a large constant $\uC{C5}>0$ and applying Theorem \ref{thm:ABR}, obtaining
\begin{equation}\label{eq:newstep2}
\wtP_{\den} \left [ \left |  \tau \bigl (v, h(v)\bigr ) - \widetilde{\E}_{\den}  \Bigl [ \tau \bigl (v, h(v)\bigr )\Bigr ]\right | > \uC{C5} j \right ]\leq C e^{-\delta \uC{C5} j}\leq  e^{-\uC{C5}' j},
\end{equation}
for a suitable (and still large, when $\uC{C5}$ is large enough) constant $\uC{C5}'>0$.
Furthermore,
\[
\begin{split}
|  \bt_{h(v)} -  \bt_{h(v')} | & = \left |  \bt_{h(v)} - \bt_v + \bt_v -  \bt_{h(v')}+\bt_{v'} - \bt_{v'}\right | \leq |\bt_{h(v)} - \bt_v - (\bt_{h(v')} - \bt_{v'})| + |\bt_v-\bt_{v'}|\\
& \stackrel{\eqref{eq:def-tau} }{=} \left | \tau \bigl (v, h(v)\bigr )-\tau \bigl (v', h(v')\bigr )\right | + |\bt_v-\bt_{v'}|.
\end{split}
\]
From \eqref{eq:tauLaw} it follows that
$
\widetilde{\E}_{\den} \bigl [ \tau \bigl (v, h(v)\bigr )\bigr ]=\widetilde{\E}_{\den}  \bigl [ \tau \bigl (v', h(v')\bigr )\bigr ]$, 
for all $v,v' \in \mathcal{N}_j$. 
Thus, the above gives
\[
\begin{split}
|  \bt_{h(v)} -  \bt_{h(v')} | & \leq  \left | \tau \bigl (v, h(v)\bigr )-\tau \bigl (v', h(v')\bigr )\right | + |\bt_v-\bt_{v'}|\\
& \leq \left | \tau \bigl (v, h(v)\bigr ) - \widetilde{\E}_{\den}  \bigl [ \tau \bigl (v, h(v)\bigr )\bigr ] \right | + \left |\tau \bigl (v', h(v')\bigr ) - \widetilde{\E}_{\den}  \bigl [ \tau \bigl (v', h(v')\bigr )\bigr ]\right | + |\bt_v-\bt_{v'}|.
\end{split}
\]
By putting this fact together with Claim \ref{claim:NiBorn} and \eqref{eq:newstep2} we obtain that for every $v, v'\in \mathcal{N}_j $,
\[
\begin{split}
& \wtP_{\den}  \Bigl [ \left |  \bt_{h(v)} -  \bt_{h(v')}\right | \leq \bigl (2\uC{C5} + \uC{C1} \bigr ) j \Bigr ]\\
& \geq \wtP_{\den}  \Bigl [ \left | \tau \bigl (v, h(v)\bigr ) - \widetilde{\E}_{\den}  \bigl [ \tau \bigl (v, h(v)\bigr )\bigr ] \right | +\\
& \qquad \quad + \left |\tau \bigl (v', h(v')\bigr ) - \widetilde{\E}_{\den}  \bigl [ \tau \bigl (v', h(v')\bigr )\bigr ]\right | + |\bt_v-\bt_{v'}| \leq \bigl (2\uC{C5} + \uC{C1} \bigr ) j \Bigr ]\\
& \geq 1- 2e^{-\uC{C5}' j} -\frac{1}{\uC{C1}}e^{- j \uC{C1}/2},
\end{split}
\]
yielding the claim.
\end{proof}
\nC{C6}
Using the notation in Lemma \ref{lemma:aux}, define the constant
\begin{equation}\label{eq:def-C6}
\uC{C6}:=2\uC{C5} + \uC{C1},
\end{equation}
furthermore, for all $j\in \{1, \ldots, \ell\} $ define the following event
\begin{equation}\label{eq:event1}
\mathcal{Z}_j:=\left \{ \sup_{v,v' \in \mathcal{N}_j}\left |  \bt_{h(v)} -  \bt_{h(v')}\right | \leq \uC{C6} j \right \}.
\end{equation}
%
%
%
%
For $\gamma\in \{\lambda, 1\}$ recall the definition of $M_{i,n}^{(\gamma)}(x)$ from \eqref{eq:defMx}.
%
For all $j\in \{1, 2, \ldots , \ell \}$, 
whenever a vertex  $\mathbf{x}_{m-j} \in \mathcal{L}_{m-j}$ hosts a seed, all variables $M_{m-j,m}^{(\lambda)} (\mathbf{x}_{m-j}) $ are independent 
and have the same distribution.
Thus introduce $M_{\lambda}(j) $ as an independent copy of $M_{m-j,m}^{(\lambda)} (\mathbf{x}_{m-j})$, that is,
\[
M_{\lambda}(j) \sim M_{m-j,m}^{(\lambda)} (\mathbf{x}_{m-j}).
\]
Similarly, let $ M_{1}(j)$ be an independent copy of 
$M_{m-j,m}^{(1)} (x_{m-j})$ for some (any) 
vertex $x_{m-j}\in \mathcal{N}_{m-j}$, that is,
\[
M_{1}(j)\sim M_{m-j,m}^{(1)} (x_{m-j}).
\]
(If the corresponding GWT is so that the subtree started at $x_{m-j}$ dies out before reaching level $m$, then we set $M_{m-j,m}^{(1)} (x_{m-j}):=+\infty$.)
As a consequence of relation \eqref{eq:definitionClambda} we see that for any large fixed constant $\mathbf{K}$, if $\lambda$ is small enough, we obtain that for all $j\in \{1, 2, \ldots , \ell \}$
\begin{equation}\label{eq:diff-minima}
\E_\mu \left [ M_{\lambda}(j) \right ] > \mathbf{K} \, \wtE_\mu  \left [ M_{1}(j)\right ].
\end{equation}
Recall $\mathcal{C}^\ast$ from \eqref{eq:definitionClambda}.
Since $\lambda$ is small enough, 
\[
\frac{\mathcal{C}^\ast }{\lambda}> 2 \uC{C6},
\]
with $\uC{C6}$ satisfying \eqref{eq:def-C6}.
Subsequently, for $j\in \{1, \ldots , \ell\}$ define
$
\bt_{\text{first}}(m-j):=\min_{x\in \mathcal{L}_{m-j}}\bt_x
$,
and consider the following (random) time interval 
\[
\mathbf{T}_{\lambda}(j):= \left [ \bt_{\text{first}}(m-j), \ \bt_{\text{first}}(m-j) + \frac{\mathcal{C}^\ast j}{2\lambda } \right ] .
\]
%
Our next aim is to use these tools to define auxiliary events that, for each $j\in \{1, \ldots , \ell \}$, contain $B_{m-j}$ and subsequently bound (from above) their probability.

\begin{figure}[h!]
\begin{center}
\includegraphics[width=13.7cm]{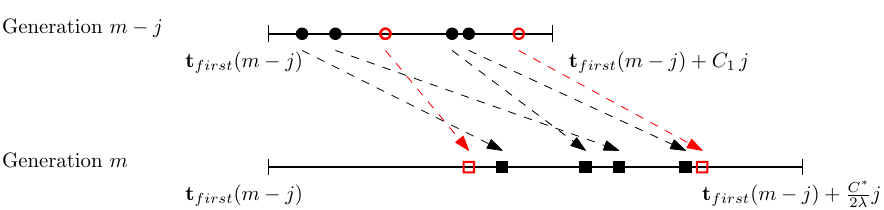}
\end{center}
\caption{Empty (red) dots are seeds and full (black) dots are non-seeds. Here the first individual to be born at the $m$-th level is a descendant of a seed at level $m-j$. Dashed arrows represent the quickest line of descendants for each particle, and the first descendant is represented as a square.
Empty (red) squares are descendants of seeds and full (black) squares are descendants of non-seeds.
Note that among the descendants of the non-seeds we only consider those of type $\Fo$, as further seeds are treated at the corresponding level.}\label{fig:time}
\end{figure}

Recall the definition of $\s$ from \eqref{eq:defSeedS}, and  the definition of $\mathcal{S}(m-j)$   from \eqref{eq:defS(m-j)}.
We now look at all particles in $\mathcal{S}(m-j)$  and order them by their birth time and subsequently, for all $j\in \{1, \ldots , \ell \}$ and for all $1\leq i\leq |\mathcal{S}(m-j) |$ we set
\begin{equation}\label{eq:def-Eij2}
\begin{split}
E_i(j):= & \bigl \{\text{the }i\text{-th particle coincides with }\s \text{ and initiates a }\Fl \\
& \text{ process reaching }\mathcal{L}_{m }\text{ within the time interval }\mathbf{T}_{\lambda}(j)\bigr \}.
\end{split}
\end{equation}
Then define
\begin{equation}\label{eq:def-Ej2}
E(j):= \bigcup_{i=1}^{|\mathcal{S}(m-j) |} E_i(j).
\end{equation}
Recall the definition of $h(\mathcal{N}_j)$ from \eqref{eq:def-h(set)}.
This time consider $i\in h(\mathcal{N}_j)$, and set
\begin{equation}\label{eq:def-Uij2b}
\begin{split}
U_i(j) := & \bigl \{\text{particle }i\text{ is \emph{not} a seed, and the }\Fo \text{ process }\\
& \text{ started from it infects } \mathcal{L}_{m } \text{ within }\mathbf{T}_{\lambda}(j) \bigr \}.
\end{split}
\end{equation}
Subsequently, define
\begin{equation}\label{eq:def-Uj2b}
U(j) := \bigcup_{i\in h(\mathcal{N}_j)} U_i(j).
\end{equation}
Note that the events $E(j)$ and $U(j)$ are \emph{not} mutually exclusive.
Each event $E(j)$ can be seen as follows: for some $i\in \{ 1, \ldots , |\mathcal{S}(m-j) |\}$ the $i$-th particle is a seed $\s$ initiating a $\Fl$ process which is atypically quick.
On the other hand, $U(j)$ states that in the set of all individuals in $h(\mathcal{N}_{j})\subset \mathcal{S}(m-j)\subset \mathcal{L}_{m-j}$, there is one that starts a $\Fo$ process that infects level $\mathcal{L}_m$ within a reasonably large time.
(Recall that for small values of $\lambda$ we have \eqref{eq:diff-minima}.)
Now  for all $1\leq i\leq |\mathcal{S}(m-j) |$ we define 
\begin{equation}\label{eq:def-Fi}
\begin{split}
F_i^{(j)} := E_i(j) \, \cap \, & \Bigl \{ \forall r\in \{1, \ldots , i-1\}, \text{ the process originated from the }r\text{-th particle }\\
& \text{ reaches }\mathcal{L}_{m }  \text{ after the process originated  by the }i\text{-th particle}\Bigr \}.
\end{split}
\end{equation}
Roughly speaking, the event $F_i^{(j)}$ implies that the $i$-th particle of $\mathcal{S}(m-j)$ is a seed initiating a $\Fl$ process that is quicker than what we expect, and at the same time all previously born particles originate processes ($\Fo$ or $\Fl$) that are slower than that.
A graphical representation can be seen in Figure \ref{fig:time}, where the horizontal segments represent time-intervals.

We proceed by showing that $E_i(j)$ is unlikely.
\begin{claim}\label{claim:constantE}
There is large constant $\overline{\mathbf{C}}_\lambda>0$ such that for all fixed $j\in \{1, \ldots , \ell\}$ and uniformly over all $1\leq i \leq |\mathcal{S}(m-j) |$, 
\begin{equation}\label{eq:PEij}
\P_\p  (E_i(j))  \leq e^{-\overline{\mathbf{C}}_\lambda j } .
\end{equation}
Furthermore, $\overline{\mathbf{C}}_\lambda$ can be made large by picking $\lambda$ sufficiently small.
\end{claim}
\begin{proof}
Recall relation \eqref{eq:definitionClambda}; to make the notation less cumbersome, for this proof we set
\[
\mathbf{C}_\lambda:= \frac{\mathcal{C}^\ast}{\lambda}.
\]
From Theorem \ref{thm:ABR} we see that for all $j\in \{1, \ldots , \ell\}$ and uniformly in $1\leq i \leq |\mathcal{S}(m-j) |$,
\[
\P_\mu (E_i(j))  \leq \P_\p \left [ \left |  M_{\lambda}(j) -\E\left ( M_{\lambda}(j)  \right ) \right | > \frac{\mathbf{C}_\lambda j}{2} \right ] \leq C e^{-\delta \mathbf{C}_\lambda j/2 } .
\]
Thus, there is a large constant $\overline{\mathbf{C}}_\lambda= \overline{\mathbf{C}}_\lambda(C,\delta,\mathbf{C}_\lambda) $ so that uniformly over $ i \in \{1, \ldots, |\mathcal{S}(m-j) |\}$
\begin{equation}\label{eq:peij-bis}
\P_\mu (E_i(j))  \leq C e^{-\delta \mathbf{C}_\lambda j/2} \leq e^{-\overline{\mathbf{C}}_\lambda j },
\end{equation}
proving the statement.
\end{proof}

The next result is where we need the assumption that either $d$ is arbitrary and $\mu$ is small, or $\mu>1/2$ is fixed and $d$ is \emph{large enough}.
In the first case we start by defining $\mu'$ so small that
\begin{equation}\label{eq:def-mu-bar-prime}
\P_{\mu'} (\mathscr{S})>1/2,
\end{equation}
and consider only $\mu<\mu'$.

\begin{lemma}\label{lemma:constantU}
For any fixed constant $\mathbf{C}>0$ the following occurs.
Suppose that $d\geq 2$ is arbitrary and $\mu$ is small enough, or $\mu\in(1/2, 1)$ is fixed and $d$ is large enough and satisfies \eqref{eq:techCondMu}.
Then for all fixed $j\in \{1, \ldots , \ell\}$,
\begin{equation}\label{eq:PUij}
 \wtP_\p (U^c(j)) \leq e^{-\mathbf{C}j}.
\end{equation}
\end{lemma}
\begin{proof}
Start by fixing the constant $\mathbf{C}>0$.
Now we want to show that $|h(\mathcal{N}_j)|$ is likely to be large, so that $\wtP_\p (U^c(j))$ satisfies the sought bound.
For this proof denote
\[
h_j:=|h(\mathcal{N}_j)|.
\]
Using Lemma \ref{lemma:new-concentr}, we see that with probability at least $1-\uC{C2}e^{-\uc{c2}j}$, the variable $h_j$ stochastically dominates a Binomial random variable 
\[
Y_j\sim \text{ Bin}\left ( \max \left \{1,\left \lfloor \uep{e1bis} \bigl (d\ns \bigr )^{(1-\uep{e1})j} \right \rfloor \right \}, \P_\mu (\mathscr{S}) \right ).
\]
In particular, since we can choose $\uep{e1bis}$ as small as we please, we can fix it so small that
\begin{equation}\label{eq:FirstC/3}
\uC{C2}e^{-\uc{c2}j}\leq \frac{e^{-\mathbf{C}j}}{3}.
\end{equation}
We proceed with a case distinction, following the two conditions in the statement.

\underline{Case (i)}: $d\geq 2$ is arbitrary and $\mu$ is small enough.
Note that the case $Y_j\sim$ Bin($1, \P_\mu (\mathscr{S}))$ is a Bernoulli random variable.
Recall that $\mu<{\mu}'$ (cf.\ \eqref{eq:def-mu-bar-prime}), then for convenience we pick a large constant $\tilde{C}=\tilde{C}(\mathbf{C})$ so that 
\begin{equation}\label{eq:CtildeCbold}
\frac{1}{\P_\mu (\mathscr{S})} e^{- \tilde{C}j /4}\leq \frac{e^{-\mathbf{C}j}}{3}.
\end{equation}
First we prove the claim for all $j$ such that
\begin{equation}\label{eq:condCtilde}
\tilde{C}j \leq \frac{1}{2} \E_\mu  Y_j = \frac{1}{2} \left  \lfloor \uep{e1bis} \bigl (d\ns \bigr )^{(1-\uep{e1})j}\right \rfloor \P_\mu (\mathscr{S}).
\end{equation}
Here we simply apply a Chernoff bound and the definition of stochastic domination,
\[
\begin{split}
\wtP_\mu \left ( h_j \leq \tilde{C}j \right )
& \leq \frac{\P_\mu \left ( Y_j \leq \frac{1}{2} \E_\mu Y_j \right )}{\P_\mu (\mathscr{S})}\leq \frac{1}{\P_\mu (\mathscr{S})} e^{-\E_\mu Y_j /8} \stackrel{\eqref{eq:condCtilde} }{\leq} \frac{1}{\P_\mu (\mathscr{S})} e^{- \tilde{C}j /4}.
\end{split}
\]
At this point our choice of $\tilde{C}$ as in \eqref{eq:CtildeCbold} gives us that for all $j$ satisfying \eqref{eq:condCtilde},
\begin{equation}\label{eq:hJlarge}
\wtP_\mu \left ( h_j \leq \tilde{C}j \right ) \leq \frac{e^{-\mathbf{C}j}}{3}.
\end{equation}
Whenever $\tilde{C}j > \frac{1}{2} \E_\mu Y_j$ we perform a different approach, using Chebyshev's inequality.
Let 
\begin{equation}\label{eq:defJ0}
j_0=j_0(\uep{e1},\uep{e1bis},d):=\max\left \{ j\in \{1, \ldots , \ell \} \ : \ \tilde{C}j \geq \frac{1}{2} \uep{e1bis} d^{(1-\uep{e1})j/2} \right \}.
\end{equation}
By definition, under our assumptions, all $j$'s for which \eqref{eq:condCtilde} fails are smaller or equal than $j_0$.
Letting
\[
\mathscr{S}_j:=\{\text{GWT survives for }j\text{ generations}\},
\]
then we observe that, given $N_j$, $h_j$ is distributed as a binomial random variable 
$
h_j\sim \text{Bin }(N_j, \P_\mu (\mathscr{S}_j))
$,
thus
\begin{equation}\label{eq:E(h_j)}
\E_\mu h_j = \E_\mu \left [ \E_\mu [h_j  \mid N_j]\right ] = \E_\mu [N_j \P_\mu (\mathscr{S}_j)] = \bigl (d\ns \bigr )^{j} \P_\mu (\mathscr{S}_j).
\end{equation}
This holds for all $j$, and hence for all $j\leq j_0$.
Since $d$ and $\mu$ must satisfy \eqref{eq:techCondMu}, there are two constants $0<\kun<\kov<1$ so that, on the one hand, 
\begin{equation}\label{eq:condKappaUn}
\kun d\ns>1;
\end{equation}
and on the other hand, for all $j\leq j_0$
$
\kov \P_\mu (\mathscr{S}_j) > \kun $.
Then,
\begin{equation}\label{eq:chebyshev}
\begin{split}
\wtP_\mu & \left ( h_j \leq \kun \bigl (d\ns \bigr )^j \right )  \leq 
\wtP_\mu \left ( h_j \leq \kov \bigl (d\ns \bigr )^j \P_\mu (\mathscr{S}_j) \right )\\
& \leq 
\wtP_\mu \left ( | h_j -\E_\mu h_j | \geq (1-\kov)\E_\mu h_j\right )\leq \frac{1}{\P_\mu (\mathscr{S})} \frac{\text{Var}_\mu(h_j)}{\left ( (1-\kov)\E_\mu h_j\right )^2}.
\end{split}
\end{equation}
To evaluate the variance of $h_j$ we reason like above, namely, 
\[
\begin{split}
\E_\mu [(h_j-\E_\mu h_j)^2] & = \E_\mu [\E_\mu [(h_j-\E_\mu h_j)^2 \mid N_j]] = \E [N_j \P_\mu(\mathscr{S}_j) (1-\P_\mu(\mathscr{S}_j))] \\
& =  \bigl (d\ns \bigr )^{j}\P_\mu(\mathscr{S}_j) (1-\P_\mu(\mathscr{S}_j)).
\end{split}
\]
Using these facts in \eqref{eq:chebyshev} we obtain
\[
\wtP_\mu \left ( | h_j -\E_\mu h_j | >(1-\kov)\E_\mu h_j\right )
\leq 
\frac{(1-\P_\mu(\mathscr{S}_j))}{ (1-\kov)^2\bigl (d\ns \bigr )^{j}\P_\mu(\mathscr{S}_j)\P_\mu (\mathscr{S})}.
\]
Now, since for all $j\geq 1$ we have $\P_\mu(\mathscr{S}_j) \geq \P_\mu(\mathscr{S})$, then for all $j\geq 1$ one has
\[
1-\P_\mu(\mathscr{S}_j)\leq 1- \P_\mu (\mathscr{S}).
\]
Hence, Chebyshev's inequality for all $j\in \{1, \ldots , j_0\}$ reduces to 
\[
\wtP_\mu \left ( | h_j -\E_\mu h_j | >(1-\kov)\E_\mu h_j\right )
\leq \frac{(1-\P_\mu(\mathscr{S}))}{ (1-\kov)^2\bigl (d\ns \bigr )^{j}\P_\mu(\mathscr{S})^2}.
\]
Since we are assuming that $\mu<{\mu}'$ (defined in \eqref{eq:def-mu-bar-prime}) then $\P_\mu(\mathscr{S})>1/2$, which yields
\[
\wtP_\mu \left ( | h_j -\E_\mu h_j | >(1-\kov)\E_\mu h_j\right )
\leq \frac{4(1-\P_\mu(\mathscr{S}))}{ (1-\kov)^2\bigl (d\ns \bigr )^{j}}.
\]
Now define $\mu'':=\mu''(\mathbf{C}, \kov)>0$ so that, uniformly over all $j\geq 1 $ one has
\begin{equation}\label{eq:def-mu-doubleprime}
\frac{2}{ (1-\kov)^2\bigl (2 (1-\mu'') \bigr )^{j}} \leq \frac{e^{-\mathbf{C}j}}{3}.
\end{equation}
Hence we can restrict our attention to $\mu<\mu''$ so that, uniformly over all  $j\in \{1, \ldots , j_0\}$ one has 
\begin{equation}\label{eq:cheb2}
\wtP_\mu \left ( h_j \leq \kun \bigl (d\ns \bigr )^{j}  \right )
\leq  \frac{4(1-\P_\mu(\mathscr{S}))}{ (1-\kov)^2\bigl (d\ns \bigr )^{j}} \leq \frac{e^{-\mathbf{C}j}}{3}.
\end{equation}
%
%
%
Now we find an upper bound for $\wtP_\p (U_i^c(j))$, for all $i \in h(\mathcal{N}_j)$.
By definition we have
\[
\begin{split}
& \wtP_\p (U_i^c(j))  \leq \wtP_\p \left [ \left |  M_{1}(j) -\wtE_\p \left [ M_{1}(j) \right ] \right | > \frac{\E ( M_{\lambda}(j) )}{2} - 
\wtE_\p \left [ M_{1}(j) \right ] \right ]+\wtP_{\mu} \left ( \text{particle }i\text{ is  seed}\right )\\
& \stackrel{\eqref{eq:diff-minima} }{\leq }\wtP_\p \left ( \left |  M_{1}(j) -\wtE_\p \left [ M_{1}(j)  \right ] \right | > (\frac{\mathbf{K}}{2}-1) \wtE_\p \left [ M_{1}(j)\right ]  \right ]+\mu \ \stackrel{\text{Thm \ref{thm:ABR}}}{\leq } \ C e^{-\delta \mathbf{K} \wtE_\p \left [ M_{1}(j)  \right ]/3 } +\mu .
\end{split}
\]
At this point fix a value $\mu'''=\mu'''(\lambda,C, \delta, \mathbf{C}, j_0)>0$ so small that
\begin{equation}\label{eq:def-mu-triple}
C e^{-\delta \mathbf{K} \wtE_{\mu'''} \left [ M_{1}(1)  \right ]/3 } +\mu'''  \leq \frac{e^{-\mathbf{C}j_0}}{3}.
\end{equation}
(This quantity is well defined because $j_0$ is independent of $\mu$, cf.\ \eqref{eq:defJ0}.)
Thus, for all $\mu<\mu'''$ and for all $\lambda$ small enough (that is for all $\mathbf{K}$ large enough), one has
\begin{equation}\label{eq:PUic}
\wtP_\p (U_i^c(j)) \leq C e^{-\delta \mathbf{K} \wtE_\p \left [ M_{1}(j)  \right ]/3 } +\mu \leq C e^{-\delta \mathbf{K} \wtE_\p \left [ M_{1}(1)  \right ]/3 } +\mu  \leq \frac{e^{-\mathbf{C}j_0}}{3}.
\end{equation}
As a consequence of this, for all $j\in \{1, \ldots , j_0\}$ and $\mu<\min \{\mu', \mu'', \mu'''\}$ (defined in \eqref{eq:def-mu-bar-prime}, \eqref{eq:def-mu-doubleprime}, \eqref{eq:def-mu-triple} respectively) we get
\[
\begin{split}
\wtP_\p (U^c(j))
& \leq \wtP_\p \left ( \bigcap_{i \in h(\mathcal{N}_j)} U_i^c(j)\right )  \leq  \wtP_\p \left ( U_1^c(j)\right )^{\kun \bigl (d\ns \bigr )^j} + \wtP_\mu \left ( h_j \leq \kun \bigl (d\ns \bigr )^j \right )+\\
&  \stackrel{\eqref{eq:PUic}, ~\eqref{eq:condKappaUn}, ~ \eqref{eq:cheb2} }{\leq } \frac{e^{-\mathbf{C}j_0}}{3} + \frac{e^{-\mathbf{C}j}}{3} \ \stackrel{j\leq j_0 }{\leq } \ e^{-\mathbf{C}j}.
\end{split}
\]
Similarly, for all $j\in \{ j_0, \ldots ,\ell \}$ and $\mu<\min \{\mu', \mu'', \mu'''\}$ we get
\[
\begin{split}
\wtP_\p (U^c(j)) & \leq \wtP_\p \left ( \bigcap_{i \in h(\mathcal{N}_j)} U_i^c(j)\right ) \\ 
& \leq  \wtP_\p \left ( U_1^c(j)\right )^{\tilde{C}j} + \wtP_\mu \left ( h_j \leq \tilde{C}j \right ) + \wtP_\mu \left (N_j\leq  \uep{e1bis} \bigl (d\ns \bigr )^{(1-\uep{e1})j} \right ) \\
& \stackrel{\text{Lemma \ref{lemma:new-concentr} }}{\leq } \wtP_\p \left ( U_1^c(j)\right )^{\tilde{C}j} + \wtP_\mu \left ( h_j \leq \tilde{C}j \right ) + \uC{C2}e^{-\uc{c2}j}\\
& \stackrel{ \eqref{eq:PUic}, \eqref{eq:condKappaUn},  \eqref{eq:hJlarge},\eqref{eq:FirstC/3} }{\leq } \frac{e^{-\mathbf{C}j_0\tilde{C}j}}{3}+\frac{e^{-\mathbf{C}j}}{3} + \frac{e^{-\mathbf{C}j}}{3}  \leq e^{-\mathbf{C}j},
\end{split}
\]
as claimed.

\underline{Case (ii)}: $\mu\in(1/2, 1)$ is fixed and $d$ is large enough.
In this case, when $\mu\in(1/2, 1)$ is fixed, $\tilde{C}>0$ can be chosen again as in \eqref{eq:CtildeCbold}.
Thus, once $\uep{e1}, \uep{e1bis}$ have been set, we can pick $d$ large enough, so that condition \eqref{eq:condCtilde} holds.
Hence the proof for this case follows from a Chernoff bound exactly as performed before, yielding the claim.
This concludes the proof.
\end{proof}

We proceed with the goal of bounding $\wtP_\p (B_{m-j})$.
For all $j$, by applying the inclusion-exclusion principle, we obtain the following bound
\begin{equation}\label{eq:auxPB}
\wtP_\p (B_{m-j})  \leq \wtP_\p \left [ B_{m-j} \cap \left ( \bigcap_{i=1}^{ \exp\{\overline{\mathbf{C}}_\lambda j /2\}}  E_i^c(j) \cap U(j) \right ) \right ] + \wtP_\p \left [ \bigcup_{i=1}^{ \exp\{\overline{\mathbf{C}}_\lambda j /2\}}  E_i(j) \cup U^c(j)  \right ].
\end{equation}
Recall the definition of $\mathcal{Z}_j$ from \eqref{eq:event1}, $A_j$ from \eqref{eq:def-Aj}, and the definition of the events $E(j)$ and $U(j)$ from \eqref{eq:def-Ej2} and \eqref{eq:def-Uj2b}.
Consider the case where all the following assumptions are satisfied, for all $j\in \{1, \ldots , \ell\}$.
\begin{itemize}
\item[(i)] There are at least $\uep{e1bis} \bigl (d\ns \bigr )^{(1-\uep{e1})j}$ individuals of type $\Fo$ alive at generation $j$ which are born by time $\uC{C1}j$ (that is, $A_j$ holds).
\item[(ii)] All their first descendants at generation $m-j$ which are of type $\Fo$ or seeds (if there are any) are born within time $\uC{C6}j$ of each other (that is, $\mathcal{Z}_j$ holds).
\item[(iii)] All particles $i\in \mathcal{S}(m-j)$ such that $i>\exp\{\overline{\mathbf{C}}_\lambda j /2\}$ are either not seeds, or seeds initiating a $\Fl$ process which is not reaching $\mathcal{L}_m$ within the time interval $\mathbf{T}_\lambda(j)$, or yet there is a particle $r$ with $r < i$ originating a process ($\Fo$ or $\Fl$) reaching $\mathcal{L}_m$ before the process originated from the $i$-th particle, that is, $\bigl (F_i^{(j)}\bigr )^c$ holds for all $i>\exp\{\overline{\mathbf{C}}_\lambda j /2\}$.
(Note that if $|\mathcal{S}(m-j)|\leq \exp\{\overline{\mathbf{C}}_\lambda j /2\}$ this assumption is empty.)
\end{itemize}
Then, necessarily at least one of the following events must hold.
Possibly, at level $\mathcal{L}_{m-j}$ there are no seeds that originate a $\Fl$ process infecting $\mathcal{L}_m$ first; or there is at least a particle $i\leq \exp\{\overline{\mathbf{C}}_\lambda j /2\}$ originating a $\Fo$ or an $\Fl$ process reaching $\mathcal{L}_m$ within the time interval $\mathbf{T}_\lambda(j)$.
(This latter possibility is ensured by assumption (iii) above.)
Hence, by considering the complementary events, for every $j\in \{1, \ldots , \ell \}$ 
\begin{equation}\label{eq:inclusionFirst}
B_{m-j} \cap \left ( \bigcap_{i=1}^{ \exp\{\overline{\mathbf{C}}_\lambda j /2\}}  E_i^c(j) \cap U(j) \right ) \ \subseteq \ \left \{  \bigcup_{i > \exp\{\overline{\mathbf{C}}_\lambda j /2\}} F_i^{(j)} \right \} \cup  \mathcal{Z}_j^c \cup  A_j^c .
\end{equation}
Hence, relations \eqref{eq:auxPB} and \eqref{eq:inclusionFirst} yield
\begin{equation}\label{eq:auxConclusion}
\begin{split}
\wtP_\p (B_{m-j}) &  \leq \left [\sum_{i > e^{\overline{\mathbf{C}}_\lambda j /2} } \wtP_\p (F_i^{(j)}) + \wtP_\mu (\mathcal{Z}_j^c) + \wtP_\mu (A_j^c) \right ] + \wtP_\p (U^c(j)) + \sum_{i=1}^{
e^{\overline{\mathbf{C}}_\lambda j /2}} \wtP_\p (E_i(j)) \\
& \stackrel{\eqref{eq:PEij} }{\leq } \sum_{i > \exp\{\overline{\mathbf{C}}_\lambda j /2\}} \wtP_\p (F_i^{(j)}) + \wtP_\mu (\mathcal{Z}_j^c) + \wtP_\mu (A_j^c)  + \wtP_\p (U^c(j)) + e^{\overline{\mathbf{C}}_\lambda j /2} e^{-\overline{\mathbf{C}}_\lambda j }\\
& \stackrel{\eqref{eq:PUij} }{\leq } \sum_{i > \exp\{\overline{\mathbf{C}}_\lambda j /2\}} \wtP_\p (F_i^{(j)}) + \wtP_\mu (\mathcal{Z}_j^c) + \wtP_\mu (A_j^c)  +  e^{-\mathbf{C}j} + e^{-\overline{\mathbf{C}}_\lambda j /2}.
\end{split}
\end{equation}
Our next aim is to show that the quantity $\sum_{i > \exp\{\overline{\mathbf{C}}_\lambda j /2\}} \wtP_\p (F_i^{(j)})$ decreases exponentially fast in $j$, which we do in the next lemma.
\nC{Q5}
\begin{lemma}\label{lemma:auxConclusion2}
With the above definitions, there is a constant $\uC{Q5}>0$ which does not depend on $\lambda$, such that for all $j\in \{1, \ldots , \ell \}$ we have
\[
\sum_{i > \exp\{\overline{\mathbf{C}}_\lambda j /2\}} \wtP_\p (F_i^{(j)}) \leq 
\uC{Q5} e^{-\overline{\mathbf{C}}_\lambda j /2},
\]
\end{lemma}
\begin{proof}
We define two auxiliary events as follows.
For all $k\in \N$ define
\[
\begin{split}
\mathcal{H}_k(i):= &  \bigl \{ \forall r\in \{1, \ldots , i-1\} \text{ with } r\text{-th particle being a non-seed, the process}\\
& \quad \text{ originated from the }r\text{-th particle reaches }\mathcal{L}_{m }  \text{ in time} \geq k\bigr \};\\
\mathcal{Y}_k^{(\lambda)}(i):= &  \bigl \{ i\text{-th particle is a seed initiating a }\Fl\text{ process reaching }\mathcal{L}_{m } \text{ in time}\leq k \bigr \}.
\end{split}
\]
(Note that the above events also depend on $j$, but we omitted this dependence to simplify the notation.)
Recall that $\mathscr{S}$ is the event that GWT survives indefinitely and let
\[
\mathbf{Q}:= \P_\p (\{o\text{ is a non-seed}\} \cap \mathscr{S} ),
\]
then fix any $r\leq i-1$ and let $\{X_l\}_l$ denote an i.i.d.\ sequence of exponential random variables with rate $1$.
We have
\[
\begin{split}
& \P (\{r\text{ is not a seed}\} \cap \{ r \text{ starts a $\Fo$ process infecting }\mathcal{L}_{m }  \text{ in time} \leq  k \})\\
& \geq \P (\{r\text{ is not a seed}\} \cap \{r=\text{root of a GWT  surviving for }j \text{ generations}\} )\P \left [\sum_{l=1}^j X_l \leq k \right ]\\
& \geq \mathbf{Q}\, \P \left (\sum_{l=1}^j X_l \leq k \right ).
\end{split}
\]
In particular, according to the relative sizes of $k$ and $j$, we shall make a case distinction, namely, fix $1\leq j \leq |\mathcal{S}(m-j)|$ and consider the following two situations.
Recall $\mathcal{C}^\ast$ from \eqref{eq:definitionClambda} and for simplicity set $\mathbf{C}_\lambda :=\frac{\mathcal{C}^\ast}{\lambda}$.
\begin{itemize}
\item[(i)] For all $1\leq k\leq \frac{\mathbf{C}_\lambda j}{2}$, we see that 
\[
\begin{split}
\mathbf{Q} \, \P \left (\sum_{l=1}^j X_l \leq k \right ) & = \mathbf{Q} \, \P \left [ \text{Poi}(k)\geq j \right ] \geq \mathbf{Q} \, \P \left [ \text{Poi}(k)= k+j \right ] \\
& \geq \mathbf{Q} e^{-k}\frac{k^{(k+j)}}{(k+j)!} \geq \mathbf{Q} \frac{e^{-k-\frac{1}{12\, (k+j)}}}{\sqrt{2\pi (k+j)}}\left ( \frac{k e}{k+j}\right )^{(k+j)},
\end{split}
\]
in the last inequality we used a standard refinement of Stirling's approximation.


\nC{Q1}
\nC{Q2}
Set $\uC{Q1}:=\mathbf{Q} \frac{e^{-\frac{1}{12}}}{\sqrt{2\pi}}$, then
\[
\begin{split}
& \mathbf{Q} \, \frac{e^{-k-\frac{1}{12\, (k+j)}}}{\sqrt{2\pi (k+j)}}\left ( \frac{k e}{k+j}\right )^{(k+j)}  \ \stackrel{k,j\geq 1}{\geq } \  
\uC{Q1} \frac{ e^{j}}{\sqrt{k+j}}\left ( \frac{k }{k+j}\right )^{(k+j)} \\
&   \stackrel{k\leq \mathbf{C}_\lambda j/2 }{\geq } \ 
\uC{Q1} \frac{ e^{j}}{\sqrt{\frac{\mathbf{C}_\lambda }{2}j+j}}\left ( \frac{k }{\frac{\mathbf{C}_\lambda }{2}j+j}\right )^{(\frac{\mathbf{C}_\lambda }{2}j+j)} 
\geq \uC{Q2} \frac{e^{j/2}}{\sqrt{\mathbf{C}_\lambda }}\left ( \frac{k }{\left ( \frac{\mathbf{C}_\lambda }{2}+1\right )j}\right )^{(\frac{\mathbf{C}_\lambda }{2}+1)j},
\end{split}
\]
for some suitable constant $\uC{Q2}=\uC{Q2}(\uC{Q1})>0$.
To ease the notation, set
\begin{equation}\label{eq:defClambda-ast}
\mathbf{C}_\lambda^\ast :=\frac{\mathbf{C}_\lambda }{2}+1,
\end{equation}
which, for all $1\leq k\leq \frac{\mathbf{C}_\lambda j}{2}$ gives
\begin{equation}\label{eq:PHx<}
\wtP_\p (\mathcal{H}_k(i)) \leq \left ( 1-\frac{\uC{Q2}}{\sqrt{\mathbf{C}_\lambda }} e^{j/2}\left ( \frac{k }{\mathbf{C}_\lambda^\ast j}\right )^{\mathbf{C}_\lambda^\ast j} \right )^i \leq \exp \left \{ -\frac{\uC{Q2}}{\sqrt{\mathbf{C}_\lambda }} \, i \, e^{j/2} \left ( \frac{k }{\mathbf{C}_\lambda^\ast j}\right )^{\mathbf{C}_\lambda^\ast j} \right \}. 
\end{equation}
\item[(ii)] For all $k > \frac{\mathbf{C}_\lambda j}{2}$ we use Theorem \ref{thm:aux-janson}(i) in the following way.
First, we note that for a sequence of i.i.d.\ exponential variables of parameter $1$ we have $\E \bigl (\sum_{l=1}^j X_l \bigr )=j$.
Then, since $k$ is much larger than $j$ (for $\lambda$ small enough) we can write
\[
\begin{split}
& \P \left (\sum_{l=1}^j X_l > k \right ) = \P \left (\sum_{l=1}^j X_l > \frac{k}{j}j \right ) \stackrel{\text{Thm }\ref{thm:aux-janson}(i)}{\leq } \frac{j}{k}e^{-j \left ( \frac{k}{j}-1-\ln \left (\frac{k}{j}\right )\right )}\leq \frac{j}{k}e^{-k/2},
\end{split}
\]
where in the last inequality we have used that $\lambda$ is small enough (and thus $k$ is large).
Hence, for all $k>\frac{\mathbf{C}_\lambda j}{2}$ we obtain
\begin{equation}\label{eq:PHx>}
\P (\mathcal{H}_k(i)) \leq \left ( \frac{j}{k}e^{-k/2}\right )^i 
\leq  \left ( \frac{2}{\mathbf{C}_\lambda} \right )^i  e^{-i \, k/2}.
\end{equation}
\end{itemize}
Using again Theorem \ref{thm:ABR} we obtain that for all $k\leq \frac{\E ( M_{\lambda}(j) )}{2} $
\[
\begin{split}
\P \left ( \mathcal{Y}_k^{(\lambda)}(i) \right ) & = \P \bigl ( M_{\lambda }(j) \leq k \bigr )  \ \leq \ \P \bigl ( \left |  M_{\lambda}(j) -\E\left ( M_{\lambda}(j) \right ) \right | > \E ( M_{\lambda}(j) )-k  \bigr  )\\
& \leq \P \left ( \left |  M_{\lambda}(j) -\E\left ( M_{\lambda}(j) \right ) \right | > \frac{\mathbf{C}_\lambda j}{2}  \right ) \leq C e^{-\delta \mathbf{C}_\lambda j/2} . 
\end{split}
\]
Hence, as we did previously (cf.\ \eqref{eq:peij-bis}) we find that for all $k\leq \frac{\mathbf{C}_\lambda j}{2}$ and all $i$,
\begin{equation}\label{eq:PKx}
\P \left ( \mathcal{Y}_k^{(\lambda)}(i) \right ) \leq e^{-\overline{\mathbf{C}}_\lambda j }.
\end{equation}
From the definition of $F_i^{(j)}$ it follows that 
\[
\wtP_\p \left ( F_i^{(j)}\right ) \leq \wtP_\p \left ( \cup_{k\in \N} \left ( \mathcal{H}_k(i) \cap \mathcal{Y}_k^{(\lambda)}(i)\right ) \right ),
\]
thus, since processes born from different particles (at the same level) are independent we have $\wtP_\p  \left (  \mathcal{H}_k(i) \cap \mathcal{Y}_k^{(\lambda)}(i) \right ) = \wtP_\p  \left (  \mathcal{H}_k(i) \right )\P \left ( \mathcal{Y}_k^{(\lambda)}(i)\right ) $.
Hence,
\[
\begin{split}
\wtP_\p \left ( F_i^{(j)}\right ) & \leq \sum_{k\geq 1} \wtP_\p  \left (  \mathcal{H}_k(i) \cap \mathcal{Y}_k^{(\lambda)}(i) \right ) \leq \sum_{k = 1}^{\mathbf{C}_\lambda j/2} \wtP_\p  \left (  \mathcal{H}_k(i) \right ) \P \left ( \mathcal{Y}_k^{(\lambda)}(i)\right )  + \!\! \sum_{k\geq  \mathbf{C}_\lambda j/2 +1} \!\! \wtP_\p  \left (   \mathcal{H}_k(i)  \right )\\
& \stackrel{\eqref{eq:PKx} }{\leq} e^{-\overline{\mathbf{C}}_\lambda j } \sum_{k = 1}^{\mathbf{C}_\lambda j/2} \wtP_\p  \left (   \mathcal{H}_k(i)  \right )  + \sum_{k\geq  \mathbf{C}_\lambda j/2 +1} \wtP_\p  \left (   \mathcal{H}_k(i)  \right ).
\end{split}
\]
Now we apply relations \eqref{eq:PHx<} and \eqref{eq:PHx>}, which give that the above is bounded by 
\[
\begin{split}
& e^{-\overline{\mathbf{C}}_\lambda j } \frac{\mathbf{C}_\lambda j}{2} \exp \left \{ -\uC{Q2} \, i \, e^{j/2}\left ( \frac{k }{\mathbf{C}_\lambda^\ast j}\right )^{\mathbf{C}_\lambda^\ast j} \right \} 
+
\sum_{k>\mathbf{C}_\lambda j/2} \left ( \frac{2}{\mathbf{C}_\lambda} \right )^i  e^{-i \, k/2}\\
& \leq e^{-\overline{\mathbf{C}}_\lambda j } \mathbf{C}_\lambda j \exp \left \{ -\uC{Q2} \, i \, e^{j/2}\left ( \frac{1 }{\mathbf{C}_\lambda^\ast j}\right )^{\mathbf{C}_\lambda^\ast j} \right \} 
+ 
2 \left ( \frac{2}{\mathbf{C}_\lambda} \right )^i  e^{-i \, \mathbf{C}_\lambda j/4},
\end{split}
\]
where we have used that $\mathbf{C}_\lambda $ is large enough.
Now it follows easily that the above quantity is summable in $i$,
giving at most
\nC{Q4}
\[
\begin{split}
& e^{-\overline{\mathbf{C}}_\lambda j } \mathbf{C}_\lambda j \,  \frac{e^{-\uC{Q2} \, e^{j/2}\left ( \frac{1 }{\mathbf{C}_\lambda^\ast j}\right )^{\mathbf{C}_\lambda^\ast j}}}{1-e^{-\uC{Q2} \, e^{j/2}\left ( \frac{1 }{\mathbf{C}_\lambda^\ast j}\right )^{\mathbf{C}_\lambda^\ast j}}} +  \frac{4}{\mathbf{C}_\lambda}   \frac{e^{-\mathbf{C}_\lambda j/4}}{1-e^{-\mathbf{C}_\lambda j/4}} \\
& \leq  \uC{Q4} e^{-\overline{\mathbf{C}}_\lambda j /2} e^{-\uC{Q2} \, e^{j/2}\left ( \frac{1 }{\mathbf{C}_\lambda^\ast j}\right )^{\mathbf{C}_\lambda^\ast j}} + \frac{8}{\mathbf{C}_\lambda}   e^{-\mathbf{C}_\lambda j/4}\\
& \leq \uC{Q4} e^{-\overline{\mathbf{C}}_\lambda j /2}  +  \frac{8}{\mathbf{C}_\lambda}   e^{-\mathbf{C}_\lambda j/4} \ \leq  \  \uC{Q4} e^{-\overline{\mathbf{C}}_\lambda j /2}  +  8 e^{-\delta \mathbf{C}_\lambda j/4},
\end{split}
\]
for a suitable constant $\uC{Q4}>0$, independent of $\lambda$.
The statement follows from the definition of $\overline{\mathbf{C}}_\lambda$ as in \eqref{eq:peij-bis} and by picking $\uC{Q5}:=\uC{Q4}+8$.
\end{proof}
At this point we are ready to conclude the proof of Proposition \ref{prop:Fo-wins}.

\begin{proof}[Proof of Proposition \ref{prop:Fo-wins}]
For any large (fixed) value of $m$ we have
\[
\begin{split}
& 
\P_{\mu} \Bigl [ 
W\text{ reached first by }\Fl \Bigr ] 
\leq \P_{\mu} \left  [ \bigcup_{j=1}^m B_j \, \cup \, \mathscr{S}^c \right ]  \leq \P_{\mu} \left  [ \bigcup_{j=1}^m B_j \right ] + \P_{\mu} \left  [ \mathscr{S}^c \right ] \\
& \stackrel{\eqref{eq:PS} }{\leq} \ \P_{\mu} \left  [ \bigcup_{j=1}^m B_j \right ] + \f_d(\mu) \ \stackrel{\text{Lemma \ref{lemma:B_i-unlikely}} }{\leq} \ e^{-\delta ' m} + \P_{\mu} \left  [ \bigcup_{j=m-\ell}^m B_j \right ] + \f_d(\mu).
\end{split}
\]
Now, in order to bound $\P_{\mu} \left  [ \bigcup_{j=m-\ell}^m B_j \right ]$ 
we use \ref{eq:auxConclusion}.
By putting together the results obtained in Lemmas \ref{lemma:A^c-unlikely}, \ref{lemma:aux} and \ref{lemma:auxConclusion2} we obtain
\[
\begin{split}
& \P_{\mu} \left  [ \bigcup_{j=m-\ell}^m B_j \right ] 
\leq \uep{e2}+\sum_{j=1}^{\ell} \left [  2e^{-\uC{C5}' j} +\frac{1}{\uC{C1}}e^{- j \uC{C1}/2} + e^{-\mathbf{C}j} + e^{-\overline{\mathbf{C}}_\lambda j /2} + \uC{Q5} e^{-\overline{\mathbf{C}}_\lambda j /4} \right ]\\
& \leq \uep{e2}+\sum_{j=1}^{\ell} 
\left [  2e^{-\uC{C5}' j} +\frac{1}{\uC{C1}}e^{- j \uC{C1}/2} + e^{-\mathbf{C}j}  + (\uC{Q5}+1) e^{-\overline{\mathbf{C}}_\lambda j /4} \right ].
\end{split}
\]
\nC{Q6}
Recall the definitions of $\mathbf{C}$ and $\overline{\mathbf{C}}_\lambda$ from Claim \ref{claim:constantE} and Lemma \ref{lemma:constantU}, and in particular recall that they could be picked as large as we wanted, when $\lambda$ is small and whenever $\mu<\min \{\mu', \mu'', \mu'''\}$ (defined in \eqref{eq:def-mu-bar-prime}, \eqref{eq:def-mu-doubleprime}, \eqref{eq:def-mu-triple} respectively) or $d$ large enough.
Now denote by 
\[
\boldsymbol{\omega}:=\min \{\uC{C5}', \uC{C1}/2, \mathbf{C}, \overline{\mathbf{C}}_\lambda /4 \},
\]
hence
\[
\P_{\mu} \left  [ \bigcup_{i=m-\ell}^m B_i \right ] 
\leq \uep{e2}+\left ( 2 +\frac{1}{\uC{C1}} + 1   + (\uC{Q5}+1) \right ) \sum_{j=1}^{\ell}  e^{-\boldsymbol{\omega} j}.
\]
Set
$
\uC{Q6}:= 5 + \uC{Q5}
$,
since $\lambda $ is chosen to be small enough, then $\boldsymbol{\omega}$ is extremely large, hence
\[
\P_{\mu} \left  [ \bigcup_{i=m-\ell}^m B_i \right ] 
\leq \uep{e2}+2 \uC{Q6}  e^{-\boldsymbol{\omega} } .
\]
Since $m$ was chosen large enough, we let
\[
\bseta_d= \bseta_d(\delta', \mu, d , \lambda, \uep{e2}, \uC{Q6} , \boldsymbol{\omega} ):= e^{-\delta ' m} + \uep{e2}+2 \uC{Q6}  e^{-\boldsymbol{\omega} } .
\]
Moreover, $\uep{e2}$ can be made arbitrarily small and $\boldsymbol{\omega}$ arbitrarily large, provided $\lambda>0$ is small enough and either $\mu$ is small enough, or $d$ is large enough.
This establishes the desired property of $\bseta_d$ in Proposition \ref{prop:Fo-wins}.
To conclude, we note that
\[
\begin{split}
& 
\P_{\mu} \Bigl [ W\text{ is infected by }\Fl \Bigr ] \\
& = 
\P_{\mu} \Bigl [ W\text{ reached first by }\Fo \cap \{W \text{ is seed}\} \Bigr ] + 
\P_{\mu} \Bigl [ W\text{ reached first by }\Fl \Bigr ] \\
& \leq \left ( 1-\bseta_d - \f_d(\mu) \right ) \mu + \left ( \bseta_d + \f_d(\mu) \right ) \leq \mu+ (1-\mu)(\bseta_d + \f_d(\mu)).
\end{split}
\]
Thus, $
\P_{\mu} \Bigl [ W\text{ is infected by }\Fo \Bigr ] \geq (1-\mu)(1-\bseta_d - \f_d(\mu))$,
as claimed.
\end{proof}


\bibliographystyle{amsplain}
\bibliography{Bibliography}

\providecommand{\bysame}{\leavevmode\hbox to3em{\hrulefill}\thinspace}
\providecommand{\MR}{\relax\ifhmode\unskip\space\fi MR }
\providecommand{\MRhref}[2]{%
  \href{http://www.ams.org/mathscinet-getitem?mr=#1}{#2}
}
\providecommand{\href}[2]{#2}
\begin{thebibliography}{10}

\bibitem{AddarioBerry-Reed}
Louigi Addario-Berry and Bruce Reed, \emph{Minima in branching random walks},
  Ann. Probab. \textbf{37} (2009), no.~3, 1044--1079. \MR{2537549}

\bibitem{athreya-ney-BP}
Krishna~B. Athreya and Peter~E. Ney, \emph{Branching processes},
  Springer-Verlag, New York-Heidelberg, 1972, Die Grundlehren der
  mathematischen Wissenschaften, Band 196. \MR{0373040}

\bibitem{Candellero-Stauffer-FppCompetition}
Elisabetta Candellero and Alexandre Stauffer, \emph{Coexistence of competing
  first passage percolation on hyperbolic graphs}, Ann. Inst. Henri
  Poincar\'{e} Probab. Stat. \textbf{57} (2021), no.~4, 2128--2164.
  \MR{4330845}

\bibitem{Dauvergne-Sly}
Duncan Dauvergne and Allan Sly, \emph{Spread of infections in a heterogeneous
  moving population}, {A}rxiv preprint: https://arxiv.org/abs/2105.11947.

\bibitem{Deijfen-nonmonotonic}
Maria Deijfen and Olle H\"{a}ggstr\"{o}m, \emph{Nonmonotonic coexistence
  regions for the two-type {R}ichardson model on graphs}, Electron. J. Probab.
  \textbf{11} (2006), no. 13, 331--344. \MR{2217820}

\bibitem{Finn-PhD}
Thomas Finn, \emph{Topics in random growth models}, 2021, Thesis
  (Ph.D.)--University of Bath.

\bibitem{Finn-Stauffer-FPPHE}
Thomas Finn and Alexandre Stauffer, \emph{Non-equilibrium multi-scale analysis
  and coexistence in competing first passage percolation}, {J.E.M.S.} (2022).

\bibitem{GM2005}
Olivier Garet and R{\'e}gine Marchand, \emph{Coexistence in two-type
  first-passage percolation models}, Ann. Appl. Probab. \textbf{15} (2005),
  no.~1A, 298--330.

\bibitem{GM2008}
\bysame, \emph{{First-passage competition with different speeds: positive
  density for both species is impossible}}, Electronic Journal of Probability
  \textbf{13} (2008), no.~none, 2118 -- 2159.

\bibitem{HP}
Olle H\"aggstr\"om and Robin Pemantle, \emph{First passage percolation and a
  model for competing spatial growth}, Journal of Applied Probability
  \textbf{35} (1998), no.~3, 683--692.

\bibitem{Haggstrom-Pemantle}
Olle H\"{a}ggstr\"{o}m and Robin Pemantle, \emph{Absence of mutual unbounded
  growth for almost all parameter values in the two-type {R}ichardson model},
  Stochastic Process. Appl. \textbf{90} (2000), no.~2, 207--222. \MR{1794536}

\bibitem{Hoffman}
Christopher Hoffman, \emph{Coexistence for richardson type competing spatial
  growth models}, Ann. Appl. Probab. \textbf{15} (2005), no.~1B, 739--747.

\bibitem{Janson-concentration}
Svante Janson, \emph{Tail bounds for sums of geometric and exponential
  variables}, Statist. Probab. Lett. \textbf{135} (2018), 1--6. \MR{3758253}

\bibitem{liggett-SISys}
Thomas~M. Liggett, \emph{Stochastic interacting systems: contact, voter and
  exclusion processes}, Grundlehren der mathematischen Wissenschaften
  [Fundamental Principles of Mathematical Sciences], vol. 324, Springer-Verlag,
  Berlin, 1999. \MR{1717346}

\bibitem{LyonsPeres-book}
Russell Lyons and Yuval Peres, \emph{Probability on trees and networks},
  Cambridge Series in Statistical and Probabilistic Mathematics, vol.~42,
  Cambridge University Press, New York, 2016. \MR{3616205}

\bibitem{Ney-Vidyashankar}
Peter.~E. Ney and Anand.~N. Vidyashankar, \emph{Harmonic moments and large
  deviation rates for supercritical branching processes}, Ann. Appl. Probab.
  \textbf{13} (2003), no.~2, 475--489. \MR{1970272}

\bibitem{Stauffer-Sidoravicius-MDLA}
Vladas Sidoravicius and Alexandre Stauffer, \emph{Multi-particle diffusion
  limited aggregation}, Invent. Math. \textbf{218} (2019), no.~2, 491--571.
  \MR{4011705}

\end{thebibliography}

\end{document}